\documentclass[12pt, oneside]{book}
\usepackage{amssymb, amsmath, amsthm, dsfont}
\usepackage{txfonts, mathrsfs, makeidx, lscape}
\usepackage{enumerate,multicol, graphicx, tikz,}
\usepackage[colorlinks=true,linkcolor=blue, pdfauthor= {Quincy \ Loney}, pdftitle={Decomposition of Level-1 Representations of $D_4^{(1)}$
With Respect to its Subalgebra $G_2^{(1)}$ in the Spinor Construction}, pdfsubject={ Kac-Moody Lie Algebras, Branching Rules} ]{hyperref}
\usepackage[colorlinks=true,linkcolor=blue]{hyperref}
\usepackage[all]{hypcap}
\usepackage[all]{xy}

\makeindex
\makeatletter
\renewcommand*\@makechapterhead[1]{%
  \vspace*{0.75in}%
  {\parindent \z@ \raggedright \normalfont
    \ifnum \c@secnumdepth >\m@ne
        \huge\bfseries \@chapapp\space \thechapter
        \par\nobreak
        \vskip 20\p@
    \fi
    \interlinepenalty\@M
    \Huge \bfseries #1\par\nobreak
    \vskip 40\p@
  }}
\renewcommand*\@makeschapterhead[1]{%
  \vspace*{0.75in}%
  {\parindent \z@ \raggedright
    \normalfont
    \interlinepenalty\@M
    \Huge \bfseries  #1\par\nobreak
    \vskip 40\p@
  }}
\makeatother

\theoremstyle{plain}
\newtheorem{thm}{Theorem}[chapter]
\newtheorem{cor}[thm]{Corollary}
\newtheorem{lem}[thm]{Lemma}
\newtheorem{prop}[thm]{Proposition}

\theoremstyle{definition}
\newtheorem{ex}[thm]{Example}
\newtheorem{dfn}[thm]{Definition}

\newtheorem*{prf}{Proof}
\newtheorem{rmk}[thm]{Remark}
\newtheorem{scon}[thm]{Sugawara Construction}
\newtheorem{ccon}[thm]{Spinor Construction}

\newcommand{\e}{\textrm{End}}
\newcommand{\Vir}{{\bf Vir}}
\def\Cl{\textrm{Cliff}}
\def\CM{\textrm{CM}}
\def\cl{\textrm{Cliff}_{\ell}}
\def\cm{\textrm{CM}_{\ell}}
\def\clz{\textrm{Cliff}(Z)}
\def\cmz{\textrm{CM}(Z)}

\newcommand{\di}{\delta_{i,j}}
\newcommand{\dm}{\delta_{m,-n}}
\newcommand{\ep}{\varepsilon}

\newcommand{\pu}{\phiup}

\DeclareMathOperator{\bv}{\textbf{vac}}
\DeclareMathOperator{\bvz}{\textbf{vac}(Z)}

\def\ds{\displaystyle}
\newcommand{\ca}{^\circledast}
\newcommand{\sub}{\subseteq}
\newcommand{\phm}{\phantom -}
\newcommand{\ph}{\phantom i}
\newcommand{\st}{\phantom{i}|\phantom{i}}

\newcommand{\la}{\langle}
\newcommand{\ra}{\rangle}
\newcommand{\inv}{^{-1}}
\newcommand{\no}{\mbox{\tiny ${\bullet\atop\bullet}$}}
\newcommand{\fno}{\mbox{\tiny ${\circ\atop\circ}$}}
\newcommand{\one}{\tfrac{1}{2}}
\newcommand{\tre}{\tfrac{3}{2}}

\newcommand{\trd}{\tfrac{1}{3}}
\newcommand{\two}{\tfrac{2}{3}}
\newcommand{\sev}{\tfrac{7}{10}}
\newcommand{\for}{\tfrac{1}{4}}
\newcommand{\fb}{\framebox}
\newcommand{\fbf}{\framebox{$44^*$}}
\newcommand{\fbo}{\framebox{$11^*22^*$}}
\newcommand{\fbt}{\framebox{$22^*44^*$}}
\newcommand{\fbs}{\framebox{$1^*234$}}
\newcommand{\fbss}{\framebox{$1^*234^*$}}

\newcommand{\bea}{\begin{eqnarray}}
\newcommand{\eea}{\end{eqnarray}}
\newcommand{\be}{\begin {equation}}
\newcommand{\ee}{\end{equation}}
\newcommand{\bbm}{\begin{bmatrix}}
\newcommand{\ebm}{\end{bmatrix}}

\newcommand{\mb}{\mathfrak b}
\newcommand{\g}{\mathfrak g}
\newcommand{\h}{\mathfrak h}
\renewcommand{\k}{\mathfrak k}
\newcommand{\ghat}{\hat{\mathfrak g}}
\newcommand{\hhat}{\hat{\mathfrak h}}
\newcommand{\khat}{\hat{\mathfrak k}}
\newcommand{\shat}{\hat{\sigma}}
\newcommand{\hV}{\hat{V}}

\DeclareMathOperator{\cee}{\mathcal C}
\DeclareMathOperator{\dee}{\mathcal D}
\DeclareMathOperator{\F}{\mathcal F}
\DeclareMathOperator{\I}{\mathscr I}
\DeclareMathOperator{\W}{\mathcal W}

\newcommand{\D}{\omega_{D_4}}
\newcommand{\B}{\omega_{B_3}}
\newcommand{\G}{\omega_{G_2}}
\newcommand{\DB}{\omega_{(D_4-B_3)}}
\newcommand{\BG}{\omega_{(B_3-G_2)}}
\newcommand{\om}{\Omega}

\def\C{\mathbb{C}}
\def\Z{\mathbb{Z}}

\newcommand{\RR}{\mathbb{R}}

\newcommand{\hz}{\mathbb Z+\tfrac{1}{2}}

\definecolor{green}{rgb}{0,.5,0}

\definecolor{yellow}{rgb}{0,0,0}\definecolor{orange}{rgb}{0,0,0}\definecolor{brown}{rgb}{0,0,0}\definecolor{red}{rgb}{0,0,0}\definecolor{purple}{rgb}{0,0,0}
\begin{document}

\frontmatter

\thispagestyle{empty}

\vbox to 1truein{}

\fontsize{16}{32pt} \selectfont  

\centerline{DECOMPOSITION OF LEVEL-1 REPRESENTATIONS OF $D_4^{(1)}$}
\centerline{WITH RESPECT TO  ITS SUBALGEBRA $G_2^{(1)}$}
\centerline{IN THE SPINOR CONSTRUCTION}
\fontsize{12}{12pt} \selectfont  
\vskip 125pt

\centerline{BY}
\vskip 10pt

\centerline{QUINCY LONEY}
\vskip 10pt

\vskip 150pt

\centerline{DISSERTATION}
\vskip 10pt

\centerline{Submitted in partial fulfillment of the requirements for}
\centerline{the degree of Doctor of Philosophy in Mathematical Sciences}
\centerline{in the Graduate School of}
\centerline{Binghamton University}
\centerline{State University of New York}
\centerline{2012}

\fontsize{12}{24pt} \selectfont 

\chapter*{Abstract}

In \cite{FFR} Feingold, Frenkel and Ries gave a spinor construction
of the vertex operator para-algebra (abelian intertwining algebra)
$ \hV=\hV^0\oplus \hV^1\oplus \hV^2\oplus \hV^3$,
whose summands are four level-1 irreducible representations of the affine
Kac-Moody algebra $D_4^{(1)}$. The triality group $S_3 =\langle  \sigma,\tau\ |\ \sigma^3 = 1 = \tau^2, \tau\sigma\tau = \sigma^{-1} \rangle$ in $Aut(\hV)$ was constructed, preserving $\hV^0$ and permuting the $\hV^i$, for  $i=1,2,3$.
$\hV$ is $\frac{1}{2}\Z$-graded where
$\hV^i_n$ denotes the $n$-graded subspace of $\hV^i$.
Vertex operators $Y(v,z)$ for $v\in\hV^0_1$ represent $D_4^{(1)}$ on $\hV$, while those for which
$\sigma(v) = v$ represent $G_2^{(1)}$.

We investigate branching rules, how $\hV$ decomposes into a direct sum of irreducible $G_2^{(1)}$ representations.
We use a two-step process, first decomposing
with respect to the intermediate subalgebra $B_3^{(1)}$, represented by $Y(v,z)$ for $\tau(v) = v$. 
There are three vertex operators,  $Y(\omega_{D_4},z)$,  $Y(\omega_{B_3},z)$, and $Y(\omega_{G_2},z)$
each representing the Virasoro algebra given by the Sugawara constructions from the three algebras.
The Goddard-Kent-Olive coset construction \cite{GKO} gives two mutually commuting coset Virasoro representations, 
provided by the vertex operators $Y(\DB,z)$ and $Y(\BG,z)$, with central charges $\one$ and $\sev$ respectively. 
The first one commutes with $B_3^{(1)}$, and the second one commutes with $G_2^{(1)}$. 
This gives the space of highest weight vectors for $G_2^{(1)}$ in $\hV$ as tensor products of
irreducible Virasoro modules $L(1/2,h^{1/2})\otimes L(7/10,h^{7/10})$.
This dissertation explicitly constructs these coset Virasoro
operators, and uses them to describe the decomposition of $\hV$ with respect to $G_2^{(1)}$.
This work also provides spinor constructions of the $\sev$ Virasoro modules, and of the two level-1 representations 
of $G_2^{(1)}$ inside $\hV$.

\tableofcontents
\mainmatter
\fontsize{12}{24pt} \selectfont

\chapter{Introduction}


This dissertation solves a specific problem in the representation theory of the infinite dimensional Kac-Moody Lie algebras. Given an irreducible representation $\hV^i$ of such an algebra $\ghat$, and a subalgebra $\ghat_0$, branching rules give the decomposition of $\hV^i$ into a direct sum of irreducible representations of $\ghat_0$. 
While the number of summands in that decomposition may be infinite, it is possible to organize them into a finite number of modules for another Lie algebra that commutes with $\ghat_0$. In many cases, that commuting Lie algebra 
is the Virasoro algebra, which is a vital part of conformal field theory in physics. The specific decompositions 
investigated here are for the four level-1 irreducible modules, $\hV^i$, $0\leq i\leq 3$, of the affine algebra $\ghat$ of type $D_4^{(1)}$ with respect to its affine subalgebra $\ghat_0$ of type $G_2^{(1)}$. We use the spinor construction from Clifford algebras of these four modules given in \cite{FFR} which provides $\hV^0\oplus\hV^1$, the Neveu-Schwarz modules, with the structure of a vertex operator superalgebra (VOSA), provides $\hV^2\oplus\hV^3$, the Ramond modules, with the structure of a module for that VOSA, and provides $\hV = \hV^0\oplus\hV^1\oplus\hV^2\oplus\hV^3$ with the structure of a vertex operator para-algebra (VOPA), also known as an abelian intertwining algebra. An action on $\hV$ of the triality group, the symmetric group $S_{\!3}$, is also provided, and gives $\ghat_0$ as certain fixed points under $S_{\!3}$. The vertex operator algebra (VOA) structure on $\hV^0$ includes vertex operators $Y(v,z)$ for $v\in\hV^0$ which represent the affine algebras $\ghat$ and $\ghat_0$, as well as an intermediate affine algebra, $\mb_1$, of type $B_3^{(1)}$, the fixed points under an order $2$ element of $S_{\!3}$. For each of these three affine algebras, $\ghat\supseteq\mb_1\supseteq\ghat_0$, there is an associated representation of the Virasoro algebra provided by vertex operators $Y(\D,z)$, $Y(\B,z)$ and $Y(\G,z)$. Using the \cite{GKO} coset construction we get vertex operators 
$Y(\DB,z)$ and $Y(\BG,z)$ providing two commuting coset representations of the Virasoro algebra with central charges $\one$ and $\sev$, respectively, both commuting with the operators representing $\ghat_0$. Under these circumstances each irreducible $\ghat$-module, $\hV^i$, decomposes into a direct sum of tensor products of the form
$L(1/2,h^{1/2})\otimes L(7/10,h^{7/10})\otimes W(\om_j)$ where $L(c,h)$ is the irreducible highest weight Virasoro module 
with central charge $c$, the Virasoro operator $L_0$ acts as $h$ on the highest weight vector, and $W(\om_j)$ is 
the irreducible highest weight module for $\ghat_0$ with highest weight $\om_j$. 

The main results of this dissertation are the following explicit decompositions for each level-1 irreducible 
$D_4^{(1)}$-module, $\hV^i$, $0\leq i\leq 3$, as a direct sum of $G_2^{(1)}$-modules, where the infinite number of summands is expressed as a finite sum of tensor products of two kinds of coset Virasoro modules. For the Neveu-Schwarz modules we have

\begin{align*} \hat V^0&= L(1/2,0)\otimes L(7/10,0)\otimes W(\om_0) \\
&\oplus L(1/2,0)\otimes L(7/10,3/5)\otimes W(\om_2) \\
&\oplus L(1/2,1/2)\otimes L(7/10,1/10)\otimes W(\om_2)\\
&\oplus L(1/2,1/2)\otimes L(7/10,3/2)\otimes W(\om_0)
\end{align*}
and
\begin{align*} \hat V^1&= L(1/2,1/2)\otimes L(7/10,0)\otimes W(\om_0) \\
&\oplus L(1/2,1/2)\otimes L(7/10,3/5)\otimes W(\om_2)\\
&\oplus L(1/2,0)\otimes L(7/10,1/10)\otimes W(\om_2)\\
&\oplus L(1/2,0)\otimes L(7/10,3/2)\otimes W(\om_0) 
\end{align*}
and combining them we have
\begin{align*}
 \hat V^0\oplus \hat V^1 &= \Big(L(1/2,0)\oplus L(1/2,1/2)\Big) \otimes \Big(L(7/10,0)\oplus L(7/10,3/2)\Big) \otimes W(\om_0) \\
 &\oplus \Big(L(1/2,0)\oplus L(1/2,1/2)\Big) \otimes \Big(L(7/10,1/10)\oplus L(7/10,3/5)\Big) \otimes W(\om_2).
\end{align*}

For the Ramond modules we have
\begin{align*} \hat V^2&= L(1/2,1/16)\otimes L(7/10,3/80)\otimes W(\om_2) \\
&\oplus L(1/2,1/16)\otimes L(7/10,7/16)\otimes W(\om_0) 
\end{align*}
and
\begin{align*} \hat V^3&= L(1/2,1/16)\otimes L(7/10,3/80)\otimes W(\om_2) \\
&\oplus L(1/2,1/16)\otimes L(7/10,7/16)\otimes W(\om_0) .
\end{align*}
We prove these results by explicitly finding highest weight vectors corresponding to each direct summand above, 
which gives containments for each $\hV^i$. We get the equalities by showing that the principal graded 
characters on both sides are equal. 

The organization of this dissertation is as follows. After giving the necessary background in Lie theory, we give an exposition of the theory of characters (graded dimension) of modules for the algebras we consider, and prove the  character identities needed to complete our decompositions. We then present the algebras and representations in the spinor construction, and construct all of the necessary vertex operators. These are used to find the highest weight vectors and verify their properties giving each of the summands in the above formulas.

\chapter{Background}
\label{cha:background}

This chapter is an attempt to keep the document self contained. We refer the reader to any appropriate text,  for example  \cite{Humphreys} or  \cite{Carter} for a detailed exposition on the material in the first section, and to \cite{Carter}, \cite{Kenji2}, or  \cite{Wan} for the sequel. 

\section{Finite Dimensional Lie Algebras, Automorphisms, and Representations}
\label{sec:finite}

\begin{dfn}\index{Lie algebra}\label{dfn:Lie}
Let $\g$ be a vector space over a field $F$ with a bilinear operation given by
the bracket  $[\cdot,\cdot] : \g \times \g \to \g$. We call $\g$ a \emph{Lie algebra} if
the following properties are satisfied  for all $x,y,z\in \g$:
\begin{enumerate}[(a)]
\item $[x,x] = 0$
\item $[x,[,y,z]] + [y,[z,x]] + [z,[x,y]] =0.$ 
\end{enumerate}
\end{dfn}
An \emph{ideal}\index{ideal} $I\subseteq \g$ is a subspace such that $[I,\g]\subseteq I$, and $\g$
is \emph{simple} \index{Lie algebra!simple} when it has no proper nontrivial ideals. A Lie algebra $\g$ is called \emph{abelian}\index{Lie algebra!abelian} when $[x,y]=0$ for 
all $x,y\in \g$, and an element $x\in\g$ is called \emph{central}\index{central} when $[x,y] =0$ for all $y\in\g$. A Lie algebra \emph{homomorphism} \index{homomorphism} is  a linear map $\phi: \g_1 \to \g_2$ where $\phi([x,y]) = [\phi(x),\phi(y)]$ for all $x,y\in \g_1.$ Lie algebra isomorphisms and automorphisms are defined accordingly. 

\begin{ex}\label{ex:ends}
For any vector space $V$, the vector space $\e(V)$ of endomorphisms forms a Lie algebra with the bracket given by the 
\emph{commutator} $[x,y] =x\circ y - y\circ x $ for all $x,y\in \e(V)$.  
\end{ex}

\begin{dfn}\index{representation}\label{dfn:rep}
A \emph{Lie algebra representation} on a vector space $V$ is a Lie algebra homomorphism
$\phi:\g \to \e(V)$, that is, $\phi[x,y]= \phi(x)\circ\phi(y) - \phi(y)\circ\phi(x)$ for all $x,y\in \g$.  This is equivalent to saying
that $V$ is a \index{module}$\g$-\emph{module}, where the action of $\g$ on $V$ is given by $x\cdot v = \phi(x)(v)$ for $x\in \g$ and 
$v\in V$. A subspace $U\subseteq V$ of a $\g$-module $V$ is called a $\g$-\emph{submodule} 
(\emph{sub-representation}) when $U$ is invariant under the action of $\g$. A $\g$-module $V$ is called \index{irreducible}\emph{irreducible} when it has no proper nontrivial submodules. $V$ is called \index{completely reducible}\emph{completely reducible} when it has a decomposition into
the direct sum of irreducible $\g$-modules. 
\end{dfn}

\begin{ex}\index{representation!adjoint}\label{ex:adj}
The \label{representation!adjoint} \emph{adjoint representation}, $ad: \g\to\e(\g)$, defined by  
$ad_x(y)= [x,y]$ for any $x,y\in\g$, is a Lie algebra representation of $\g$ on itself. 
\end{ex}

\begin{dfn}\index{Cartan!subalgebra}\label{dfn:csa}
We call $\h\sub\g$ a \emph{Cartan subalgebra} (CSA) of $\g$ when $\h$ is a maximal abelian subalgebra of $\g$
such that $\{ad_h\ |\ h\in\h\}$ is simultaneously diagonalizable in $\e(\g)$. The dimension of any CSA, independent of choice, is the \emph{rank} \index{rank} of $\g$. 
\end{dfn}

\begin{dfn}\index{root space}\label{dfn:rootspace}
A \emph{root space}, $\g_\alpha\sub\g$, is a simultaneous eigenspace for the $ad$ action of a CSA $\h$, that
is, $\g_\alpha = \{x\in\g\ |\ [h,x] = \alpha(h)x, \hbox{ for all } h\in\h\}$, where $\alpha\in\h^*$ is a linear
functional in the dual space $\h^*$ such that $\g_\alpha\neq 0$. However, $\g_0 = \h$ is not called a root space.
Denote by $\Phi = \{\alpha\in\h^*\ |\ \alpha\neq 0 \hbox{ and } \g_\alpha\neq 0\}$ the set of\index{roots}  \emph{roots} of $\g$. 
Note that $[\g_{\alpha},\g_{\beta}]\sub\g_{\alpha+\beta}$ for any $\alpha,\beta\in\Phi\cup\{0\}$. 
A decomposition $\ds\g=\h\oplus \bigoplus_{\alpha\in\Phi}\g_{\alpha}$ is called a 
\emph{root space decomposition of $\g$}\index{root space!decomposition}.
\end{dfn}

By the famous Cartan-Killing theorem, the finite dimensional simple Lie algebras over the complex numbers are classified  into four infinite families and five exceptional Lie algebras. Below are these Lie algebras and
their associated Dynkin diagrams\index{Dynkin diagram}. The four families are:
$$
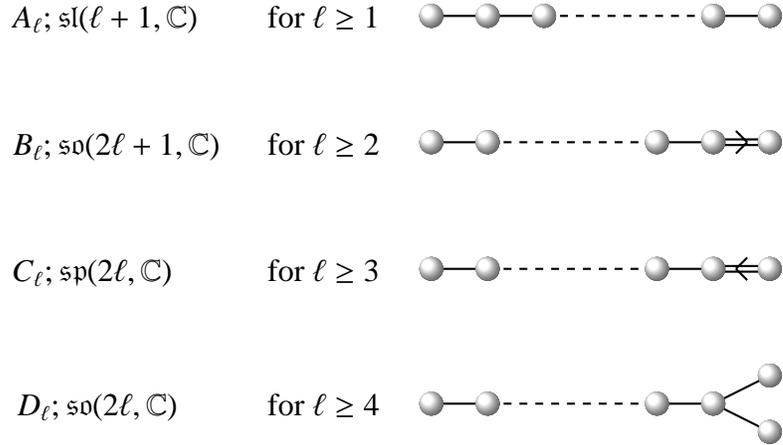
\begin{array}{lcr}
\ph A_\ell;  \mathfrak{sl}(\ell+1,\C)  & \hbox{ for } \ell\geq 1 & 
 \begin{tikzpicture}[scale=.75, thick]
\draw (0,0)--(1,0);
\draw (1,0)--(2,0);
\draw[dashed] (2,0)--(5,0);
\draw (5,0)--(6,0);
\shade[ball color =white](0,0)circle(6pt);
\shade[ball color =white](1,0)circle(6pt);
\shade[ball color =white](2,0)circle(6pt);
\shade[ball color =white](5,0)circle(6pt);
\shade[ball color =white](6,0)circle(6pt);
 \end{tikzpicture}  \\
\\

\ph B_\ell;  \mathfrak{so}(2\ell+1,\C) &\hbox{ for } \ell \geq 2& 
\begin{tikzpicture}[scale=.75, thick]
\draw (0,0)--(1,0);
\draw[dashed] (1,0)--(4,0);
\draw(4,0)--(5,0);
\draw (5,-.05)--(6,-.05);
\draw (5,.05)--(6,.05);
\draw(5.4,.2)--(5.6,0);
\draw(5.4,-.2)--(5.6,0);
\shade[ball color =white](0,0)circle(6pt);
\shade[ball color =white](1,0)circle(6pt);
\shade[ball color =white](4,0)circle(6pt);
\shade[ball color =white](5,0)circle(6pt);
\shade[ball color =white](6,0)circle(6pt);
 \end{tikzpicture} \\ 
\\

\ph C_\ell;  \mathfrak{sp}(2\ell,\C) &\hbox{ for } \ell \geq 3 &
 \begin{tikzpicture}[scale=.75, thick]
\draw (0,0)--(1,0);
\draw[dashed]  (1,0)--(4,0);
\draw(4,0)--(5,0);
\draw (5,-.05)--(6,-.05);
\draw (5,.05)--(6,.05);
\draw(5.4,0)--(5.6,.2);
\draw(5.4,0)--(5.6,-.2);
\shade[ball color =white](0,0)circle(6pt);
\shade[ball color =white](1,0)circle(6pt);
\shade[ball color =white](4,0)circle(6pt);
\shade[ball color =white](5,0)circle(6pt);
\shade[ball color =white](6,0)circle(6pt);
 \end{tikzpicture} \\ \\

\begin{array}{l} D_\ell;  \mathfrak{so}(2\ell,\C)  \\
\vspace{1pt} \end{array} &
\begin{array}{c} \hbox{ for } \ell \geq 4 \\
\vspace{1pt} \end{array}
& 
 \begin{tikzpicture}[scale=.75, thick]
\draw (0,-4)--(1,-4);
\draw[dashed]  (1,-4)--(4,-4);
\draw(4,-4)--(5,-4);
\draw(5,-4)--(6,-3.5);
\draw(5,-4)--(6,-4.5);
\shade[ball color =white](0,-4)circle(6pt);
\shade[ball color =white](1,-4)circle(6pt);
\shade[ball color =white](4,-4)circle(6pt);
\shade[ball color =white](5,-4)circle(6pt);
\shade[ball color =white](6,-3.5)circle(6pt);
\shade[ball color =white](6,-4.5)circle(6pt);
 \end{tikzpicture} 
 
\end{array}
$$

\noindent the special linear algebras, orthogonal algebras, and symplectic algebras;
and the five exceptional algebras are:
$$
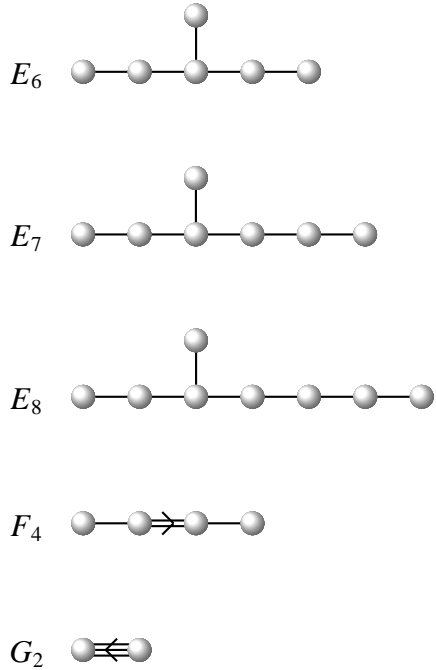
\begin{array}{ll}

 E_6&
 \begin{tikzpicture}[scale=.75, thick]
\draw (0,0)--(1,0);
\draw (1,0)--(2,0);
\draw (2,0)--(3,0);
\draw (2,0)--(2,1);
\draw (3,0)--(4,0);
\shade[ball color =white](0,0)circle(6pt);
\shade[ball color =white](1,0)circle(6pt);
\shade[ball color =white](2,0)circle(6pt);
\shade[ball color =white](2,1)circle(6pt);
\shade[ball color =white](3,0)circle(6pt);
\shade[ball color =white](4,0)circle(6pt);
 \end{tikzpicture} \\ \\ 

E_7&
 \begin{tikzpicture}[scale=.75, thick]
\draw (3,-1)--(4,-1);
\draw (4,-1)--(5,-1);
\draw (5,-1)--(6,-1);
\draw (5,-1)--(5,0);
\draw (6,-1)--(7,-1);
\draw (7,-1)--(8,-1);
\shade[ball color =white](3,-1)circle(6pt);
\shade[ball color =white](4,-1)circle(6pt);
\shade[ball color =white](5,-1)circle(6pt);
\shade[ball color =white](5,0)circle(6pt);
\shade[ball color =white](6,-1)circle(6pt);
\shade[ball color =white](7,-1)circle(6pt);
\shade[ball color =white](8,-1)circle(6pt);
 \end{tikzpicture} \\ \\
 
 E_8&
  \begin{tikzpicture}[scale=.75, thick]
\draw (0,-2)--(1,-2);
\draw (1,-2)--(2,-2);
\draw (2,-2)--(3,-2);
\draw (2,-2)--(2,-1);
\draw (3,-2)--(4,-2);
\draw (4,-2)--(5,-2);
\draw (5,-2)--(6,-2);
\shade[ball color =white](0,-2)circle(6pt);
\shade[ball color =white](1,-2)circle(6pt);
\shade[ball color =white](2,-2)circle(6pt);
\shade[ball color =white](2,-1)circle(6pt);
\shade[ball color =white](3,-2)circle(6pt);
\shade[ball color =white](4,-2)circle(6pt);
\shade[ball color =white](5,-2)circle(6pt);
\shade[ball color =white](6,-2)circle(6pt);
 \end{tikzpicture} \\ \\

 F_4&
 \begin{tikzpicture}[scale=.75, thick]
\draw (0,0)--(1,0);

\draw (2,0)--(3,0);
\draw (1,-.05)--(2,-.05);
\draw (1,.05)--(2,.05);
\draw(1.4,.2)--(1.6,0);
\draw(1.4,-.2)--(1.6,0);
\shade[ball color =white](0,0)circle(6pt);
\shade[ball color =white](1,0)circle(6pt);
\shade[ball color =white](2,0)circle(6pt);
\shade[ball color =white](3,0)circle(6pt);
 \end{tikzpicture} \\ \\

 G_2&
 \begin{tikzpicture}[scale=.75, thick]
\draw(5,0)--(6,0);
\draw (5,-.1)--(6,-.1);
\draw (5,.1)--(6,.1);
\shade[ball color =white](5,0)circle(6pt);
\shade[ball color =white](6,0)circle(6pt);
\draw(5.4,0)--(5.6,.2);
\draw(5.4,0)--(5.6,-.2);
 \end{tikzpicture}

\end{array}
$$

When $\ell=4$ the Dynkin diagram for $D_4$ has an $S_3$ symmetry called \index{triality}\emph{triality}. Let
$\sigma$ be the diagram automorphism which rotates the outer nodes counterclockwise
and fixes the middle node. Let $\tau$ be the order-two diagram automorphism which fixes the
middle and leftmost nodes while permuting the two rightmost nodes. Let $G$ be the permutation group generated by
$\sigma$ and $\tau$. These diagram automorphisms correspond to Lie algebra outer automorphisms which we denote by the same symbols. Notice that $\tau$, $\sigma\tau$ and $\sigma^2\tau$ are the order-two automorphisms which fix the middle node and one other node. For each of those order-two elements, the fixed points of $D_4$ under it form 
a Lie subalgebra of type $B_3$, $D_4^{\tau}\cong B_3$, $D_4^{\sigma\tau}\cong B_3$ and $D_4^{\sigma^2\tau}\cong B_3$. The fixed points of $D_4$ under $G$ form a Lie subalgebra of type $G_2$ inside each of those copies of $B_3$, but in fact, that subalgebra is the fixed points under just $\sigma$, $D_4^{\sigma}\cong G_2$.  

\begin{dfn}\label{dfn:rootsystem} Let $\mathbb R^\ell$ be the Euclidean space with an inner product $(\cdot,\cdot)$.
A subset $\Phi\sub \mathbb R^\ell$ is called a finite \emph{root system}\index{root system} if the following four
properties are satisfied.
\begin{enumerate}[(a)] 
\item $\Phi$ is finite, spans $\RR^\ell$ and does not contain $0$.
\item For any $\alpha\in\Phi$, $k \alpha\in\Phi$ if and only if $k = \pm 1$.
\item For each $\alpha\in \Phi$, the reflection $r_{\alpha}$  through the
hyperplane perpendicular to $\alpha$ leaves $\Phi$ invariant. The \index{reflection}reflection is given by $r_{\alpha}(\lambda) =\lambda-\la\lambda,\alpha\ra\alpha $ where $\ds\la\lambda, \alpha\ra =2\frac{(\lambda, \alpha)}{(\alpha, \alpha)}$. 
\item For $\alpha,\beta\in \Phi$, $\la\alpha,\beta\ra \in \Z$. 
\end{enumerate}
\end{dfn}

\begin{dfn}\label{dfn:roots}
A subset $\Delta = \{\alpha_1,\cdots,\alpha_\ell\}$ of a root system $\Phi$ is called
a \emph{set of simple roots} \index{roots!simple} if span$(\Delta)$ = span$(\Phi)$, so $\Delta$ is a basis of 
$\mathbb R^\ell$, 
and for each root $\alpha\in\Phi$, either $\alpha$ or $-\alpha$ can be written as a non-negative integral linear
combination of the \emph{simple} roots, $\sum_{i=1}^\ell k_i \alpha_i$ for $0\leq k_i\in \Z$.
The \index{roots!positive}\emph{positive roots} relative to a choice of $\Delta$ are 
$\Phi^+ = \{\alpha\in\Phi\ |\ \alpha = \sum_{i=1}^\ell k_i \alpha_i, 0\leq k_i\in \Z\}$, the  \index{roots!negative}\emph{negative roots} are 
$\Phi^- = - \Phi^+$, and $\Phi =\Phi^+ \cup \Phi^-$. The \index{lattice!root}\emph{root lattice} of $\Phi$ is
the integral span of the simple roots, $Q_\Phi = \sum_{i=1}^\ell \Z\alpha_i$ and we write 
$Q_\Phi^+ = \{ \sum_{i=1}^\ell k_i \alpha_i\ |\  0\leq k_i\in \Z\}$. 
There is a \index{roots!partial order}\emph{partial order} on $\RR^\ell$ defined by 
$\mu\leq\lambda$ when $\lambda - \mu \in Q_\Phi^+$. 

\end{dfn}

\begin{dfn}\label{dfn:Dynkin Diagram} Let $\Phi$ be a root system and $\Delta$ a set of simple roots in $\Phi$. Define the 
associated \index{Dynkin diagram}\emph{Dynkin diagram} to be the graph whose vertices (nodes) correspond to the simple roots, and whose edges are determined as follows. The number of edges between nodes $\alpha_i$ and $\alpha_j$ for $i\neq j$ is 
$\la\alpha_i,\alpha_j\ra \la\alpha_j,\alpha_i\ra$, an integer between $0$ and $3$. If this integer is $1$ there is an undirected edge between those nodes. If this integer is $2$ or $3$ then $(\alpha_i,\alpha_i) \neq (\alpha_j,\alpha_j)$ and the edges are directed by an arrow aimed from the longer root to the shorter one. The information encoded in the Dynkin diagram is equivalent to the \index{Cartan!matrix}\emph{Cartan matrix} $A = [ \la\alpha_i,\alpha_j\ra ]$. 
\end{dfn}

\begin{dfn}\index{Weyl!group}\label{dfn:Weyl}
Let $\Phi$ be a root system. The \emph{Weyl group}, $\W$, is the
group generated by the \index{reflection!simple}\emph{simple reflections} $r_{\alpha_i}$, 
$\W = \la r_{\alpha_i} \ | \ \alpha_i\in \Delta\ra$.
\end{dfn}

\begin{rmk} In the proof of the classification of finite dimensional simple Lie algebras over the complex numbers, the existence of a Cartan
subalgebra, $\h$, with $dim(\h) = \ell$, gives a finite set of nonzero linear functionals, $\Phi\subset \h^*$, and a root space decomposition 
of the Lie algebra. It is shown that there is a positive definite inner product on the real span of $\Phi$, which is therefore isomorphic to 
$\RR^\ell$, and $\Phi$ satisfies the axioms of a root system given above. It is shown that $\Phi$ has a set of simple roots, and the 
Weyl group generated by the simple reflections acting on the real span of $\Phi$ can also be understood as acting on $\h^*$, which is the
complex span of $\Phi$. 
\end{rmk}

\begin{dfn}\label{dfn:weight} Let $V$ be a finite dimensional $\g$-module such that the CSA $\h$ acts
simultaneously diagonalizably on $V$, with simultaneous eigenspaces (\index{weight space}\emph{weight spaces})
$V_{\mu} = \newline \{v\in V \ | \ h\cdot v = \mu(h) v, \hbox{ for all }h\in\h\}$ where $\mu\in\h^*$. 
Call $\mu$ a \emph{weight}\index{weight} of $V$ if $V_\mu \neq 0$ and write 
$\Pi(V) = \{\mu\in\h^*\ |\ V_\mu \neq 0\}$ for the set of all weights of $V$. 
We say that a nonzero vector $v\in V$ is a \emph{highest weight vector} \index{highest weight!vector} (HWV)
when $v\in V_\lambda$ for some $\lambda\in\Pi(V)$ and $x\cdot v = 0$ for all $x\in\g_\alpha$ with $\alpha\in\Phi^+$.
If $V$ is irreducible then $V$ contains a highest weight vector unique up to scalar multiples, and $V$ is uniquely determined by the weight $\lambda$ of that HWV, so we label that module by $V^\lambda$ and call it a  \emph{highest weight representation}.
\index{highest weight!representation} In that case, we write $\Pi(V^\lambda) = \Pi^\lambda$ and for any weight 
$\mu\in\Pi^\lambda$ we have $\mu \leq \lambda$. The \index{weight space!decomposition}\emph{weight space decomposition} of $V^\lambda$ is its direct sum
decomposition into weight spaces, $\ds V^{\lambda} = \bigoplus_{\mu\in\Pi^{\lambda}} V^{\lambda}_{\mu}$. 
For $\g$ finite dimensional simple, $\g$ itself is an irreducible $\g$-module by the $ad$-action, and its highest weight is its highest root, denoted by $\theta$.
\end{dfn}

\begin{dfn}\index{character}\label{def:character}
For a finite dimensional irreducible highest weight $\g$-module $V^{\lambda}$, the \emph{character}\index{character}
of $V^\lambda$ is the following element of $\Z[\h^*] = \{\sum_{\mu\in S} n_\mu e^\mu\st n_\mu\in\Z,S\sub\h^*, |S|<\infty\}$, the integral group ring over $\h^*$, 
$ch(V^{\lambda}) = \sum_{\mu\in\Pi^{\lambda}} dim(V_{\mu}^{\lambda}) e^{\mu}$. 
The character of a direct sum of such $\g$-modules is the sum of the characters of the summands. The character of a tensor
product of such $\g$-modules is the product of the characters of the tensor factors. 
\end{dfn}

\begin{dfn}\index{lattice!weight}\label{def:lattice} The \emph{integral weight lattice} associated with a root system $\Phi$ is
$P_\Phi = \newline \{\lambda\in\h^*\ |\ \la\lambda,\alpha_i\ra\in\Z, \alpha_i\in\Delta\}$. It can be seen that 
$P_\Phi = \sum_{i=1}^\ell \Z\lambda_i$ where $\lambda_i$ is the \index{weight!fundamental}\emph{fundamental weight} defined by 
$\la\lambda_i, \alpha_j\ra = \delta_{i,j}$ for $1\leq i,j\leq \ell$, and $Q_\Phi \subseteq P_\Phi$. 
The set of \emph{dominant integral weights} \index{weight!dominant integral} is 
$P_\Phi ^+ = \{\lambda\in\h^*\ |\ 0\leq\la\lambda,\alpha_i\ra\in\Z, \alpha_i\in\Delta\}$ and it is clear that
$P_\Phi ^+ = \{ \sum_{i=1}^\ell n_i\lambda_i \ | \ 0\leq n_i\in \Z\}$. For all $\lambda \in P^+$ there is a finite dimensional 
irreducible representation $V^{\lambda}$ of $\g$ with highest weight $\lambda$ and all of its weights, 
$\Pi^\lambda = \{w(\mu)\ |\ \mu\in P_\Phi ^+, \mu\leq\lambda, w\in\W\}$, are in
$P_\Phi$. Also, for each finite dimensional irreducible $\g$-module $V$ there is a $\lambda\in P^+$ such that 
$V = V^\lambda$, and $\Pi^\lambda$ is invariant under the action of the Weyl group, $\W$.
\end{dfn}

\begin{thm}\index{Weyl!character formula}\label{def:weylchar} For any $\lambda\in P_\Phi^+$ we have the following Weyl character formula for the character of irreducible finite dimensional highest weight $\g$-module $V^\lambda$:
$$ ch(V^{\lambda})=\frac{S_{\lambda+\rho}}{S_{\rho}}\quad \hbox{ where } \quad
S_{\mu}= \sum_{w\in\W} sgn(w) e^{w\mu-\rho}\quad\hbox{ and }\quad \rho = \sum_{i=1}^\ell \lambda_i.$$ 
The \index{Weyl!denominator formula} Weyl denominator formula is:
\be\label{welyden}S_{\rho} = \prod_{\alpha\in\Phi^+} (1-e^{-\alpha}).\ee
Note that setting $u_i = e^{-\alpha_i}$ for $1\leq i\leq \ell$, both $S_\rho$ and $e^{-\lambda}\ ch(V^{\lambda})$ are in the 
group ring $\Z[u_1,\cdots,u_\ell]$. 
\end{thm}

\begin{dfn}\index{Clifford algebra}\label{def:CliffordAlg} Let $V$ be a vector space over a field $F$ with a symmetric bilinear form $(\cdot,\cdot)$. The \emph{Clifford algebra} $\Cl = \Cl(V, (\cdot,\cdot) )$ is the associative algebra over $F$ generated by $V$ with unit element $1$ and
relations $ab+ba = (a,b)1$ for any $a,b\in V$. If $\I$ is any left ideal of $\Cl$ then $\CM = \Cl/\I$ is a left $\Cl$-module\index{module!Clifford}, that is, a vector space on which $\Cl$ acts by left multiplication as endomorphisms satisfying the relations in $\Cl$, thus giving a representation of $\Cl$ on $\CM$. 

Let $V^{\pm}\cong\C^\ell$, $V = V^+\oplus V^-$, and let $(\cdot,\cdot)$ be non-degenerate such that each subspace $V^{\pm}$ is isotropic. For $a,b\in V$ let $\fno ab\fno = \one(ab-ba)\in\Cl$ and let $\g = span\{\fno ab\fno\st a,b\in V\}$. Then the relations in $\Cl$ imply that $\g$ is closed under commutator, thus forming a Lie algebra which can be shown to be of type $D_\ell$. The relations in $\Cl$ also imply that for all $a,b,c\in V$, the commutator $[\fno ab\fno,c]\in V$, showing that $V$ is a $\g$-module. Let $\I$ be the left ideal of $\Cl$ generated by $V^+$, so that $\CM = \Cl/\I = (\wedge V^-)\cdot \bv$ where $\bv = 1+\I$, and $\CM = \CM^0\oplus \CM^1$ is the decomposition according to \emph{parity}, the $\Z_2$-grading of $\wedge V^-$. Although $\CM$ is an irreducible $\Cl$-module, both $\CM^0$ and $\CM^1$ are invariant and irreducible under the left multiplication by $\g$. These explicit constructions give the adjoint representation of dimension $2\ell^2 - \ell$ and three other irreducible representations of $D_\ell$, the \emph{natural}\index{representation!natural} representation of dimension $2\ell$, and the two \emph{semi-spinor}\index{representation!spinor} representations, each of dimension $2^{\ell-1}$. We will be using only the case of $\ell=4$, and details will be provided later. 
\end{dfn}

\section{Infinite Dimensional Lie Algebras and Representations}
\label{sec:infinite}

\begin{dfn}\index{Lie algebra!Kac-Moody}\label{def:KacMoody}
For a finite dimensional simple Lie algebra $\g$, define the untwisted \emph{affine Kac-Moody Lie algebra} 
$\ghat = \g\otimes \C[t,t\inv]\oplus\C c\oplus\C d$, where we use the notation $x(m) = x\otimes t^m $ for $x\in\g$ and $m\in\Z$, and the Lie brackets are 
\be\label{KMb}[x(m),y(n)] =  [x,y](m+n)+m\dm(x,y)c,\ \  [c,\hat\g] =  0,\hbox{ and } [d,x(m)] = mx(m)\ee
for $x,y\in\g$ and $m,n\in\Z$, and where $(\cdot,\cdot)$ is a non-degenerate invariant symmetric bilinear form on $\g$ with a certain standard normalization. We identify $\g$ with $\g\otimes t^0\sub\ghat$ so $\g$ is a Lie subalgebra of $\ghat$.
\end{dfn}

This affine algebra $\ghat$ has \index{Cartan!subalgebra!affine}CSA $H = \h\oplus \C c \oplus \C d$, a maximal abelian subalgebra whose $ad$ action is diagonalizable, and \index{root space!decomposition!affine}root space decomposition 
$\ds\ghat = H \oplus \bigoplus_{\alpha\in\hat\Phi} \ghat_{\alpha}$ where the \index{root system!affine}
\emph{affine root system} is denoted  $\hat\Phi$, and we write $dim(\ghat_{\alpha}) = mult(\alpha)$. In order to
describe the affine root system in the dual of the Cartan subalgebra, $H^*$, first extend the simple roots $\alpha_i\in\Delta$ to be 
the elements in $H^*$ which are zero on $c$ and $d$. Then define the elements $c^*, d^*=\delta\in H^*$ by 
$c^*(\h) = 0 = d^*(\h)$, $c^*(c) = 1 = d^*(d)$ and $c^*(d) = 0 = d^*(c)$. Then 
$\hat\Phi = \{\alpha + m\delta\ |\ \alpha\in\Phi, m\in\Z\} \cup \{n\delta\ |\ 0\neq n\in\Z\}$ 
and the simple roots of $\ghat$ are $\hat\Delta = \{\alpha_0, \alpha_1,\cdots,\alpha_\ell\}$ where the new  \index{roots!affine}affine simple root, 
$\alpha_0 = \delta - \theta$ and $\theta$ is the highest root of $\Phi$. 
The dual Cartan, $H^*$, has a basis $\{\alpha_1, ... \alpha_\ell, c^*, \delta \}$. The inner product on $H$ is the 
extension of the inner product on $\h$ such that $(\h,c) = 0 = (\h,d)$, $(c,c) = 0 = (d,d)$ and $(c,d) = 1$.
The inner product on $H^*$ is the extension of the inner product on $\h^*$ such that 
$(\h^*,c^*) = 0 = (\h^*,\delta)$, $(c^*,c^*) = 0 = (\delta,\delta)$ and $(c^*,\delta) = 1$. 
The partial order on $\h^*$ extends to a partial order on $H^*$ by $\mu\leq\lambda$ when 
$\lambda - \mu = \sum_{i=0}^\ell a_i \alpha_i$ for some $0\leq a_i\in\Z$.

\begin{ex}
The three untwisted affine Kac-Moody Lie algebras that we will be concerned with are $D_4^{(1)}, B_3^{(1)}$ and $G_2^{(1)}$. Below are their \label{Dynkin Diagram!affine} Dynkin diagrams, respectively.  The darker node corresponds to the affine root $\alpha_0$. 

\begin{center}
\begin{tikzpicture}
\draw(0:0cm)--(45:1cm);
\draw(0:0cm)--(135:1cm);
\draw(0:0cm)--(225:1cm);
\draw(0:0cm)--(315:1cm);
\shade[ball color =white](0:0cm)circle(6pt);
\shade[ball color =white](45:1cm)circle(6pt);
\shade[ball color =black](135:1cm)circle(6pt);
\shade[ball color =white](225:1cm)circle(6pt);
\shade[ball color =white](315:1cm)circle(6pt);
 \end{tikzpicture}
\qquad 
\begin{tikzpicture}

\draw(4,.5)--(5,0);
\draw(4,-.5)--(5,0);
\draw (5,-.05)--(6,-.05);
\draw (5,.05)--(6,.05);
\draw(5.4,.2)--(5.6,0);
\draw(5.4,-.2)--(5.6,0);
\shade[ball color =white](5,0)circle(6pt);
\shade[ball color =white](6,0)circle(6pt);
\shade[ball color =black](4,.5)circle(6pt);
\shade[ball color =white](4,-.5)circle(6pt);
\end{tikzpicture}
\qquad 
\begin{tikzpicture}
\draw(5,0)--(4,0);
\draw(5,0)--(6,0);
\draw (5,-.1)--(6,-.1);
\draw (5,.1)--(6,.1);
\shade[ball color =black](4,0)circle(6pt);
\shade[ball color =white](5,0)circle(6pt);
\shade[ball color =white](6,0)circle(6pt);
\draw(5.4,.2)--(5.6,0);
\draw(5.4,-.2)--(5.6,0);
 \end{tikzpicture}
 \qquad 

\end{center}

\end{ex}

\begin{dfn}
The \index{Weyl!group!affine}\emph{affine Weyl group}, $\hat\W$, is the group generated by the simple reflections 
$r_{\alpha_i}$ on $H^*$, $\hat\W = \la r_{\alpha_i}\ |\ \alpha_i\in\hat\Delta\ra$.
\end{dfn}
\begin{dfn}
The \index{weight!fundamental}\emph{fundamental weights} of $\ghat$ are $\Lambda_0=c^*$, and $\Lambda_i = n_ic^*+\lambda_i$, $1\leq i\leq \ell$. 
The integral weight lattice and dominant weights are defined similarly as in Definition \ref{def:lattice}, 
$\hat P = \sum_{i=0}^\ell\Z\Lambda_i$ and $\hat P^+ = \{ \sum_{i=0}^\ell m_i\Lambda_i \ | \ 0\leq m_i\in \Z\}$. 
\end{dfn}
\begin{dfn}\index{level}\label{def:level}
For each $\Lambda = \sum_{i=0}^\ell m_i\Lambda_i \in \hat P^+$ there exists an irreducible representation of $\ghat$, $V^{\Lambda}$. The \emph{level} of the representation, $\Lambda(c) = m_0 + \sum_{i=1}^\ell m_i n_i$, is the scalar by which the central element $c$ acts on $V^{\Lambda}$. The module has a \index{weight space!decomposition}weight space decomposition 
$\ds V^{\Lambda}=\bigoplus_{\mu\in \Pi^{\Lambda}} V_{\mu}^{\Lambda} $,  where we have as before, 
$V^\Lambda_\mu = \{v\in V^\Lambda \ | \ h\cdot v = \mu(h) v, \hbox{ for all }h\in H\}$ where $\mu\in H^*$ and 
$\Pi^{\Lambda} = \{ \mu \in \hat P \ | \ V_{\mu}^{\Lambda}\neq 0\} = \{ w(\mu)\in \hat P \ | \ \mu\in\hat P^+,
\mu\leq\Lambda, w\in\hat\W \}\sub \{\mu\in \hat P\st \mu\leq\Lambda\}$. 
\end{dfn}

\begin{dfn}\index{graded dimension}\label{def:graded dimension} Generalizing the character of a finite dimensional module, the \emph{graded dimension} 
\be\label{gradedim} gr(V^{\Lambda}) = e^{-\Lambda} \sum_{\mu\in\Pi^{\Lambda}} dim(V_{\mu}^{\Lambda})e^{\mu}\in \Z[[u_i\st 0\leq i\leq \ell]]\ee
is a formal power series of $\ell+1$ variables, $u_i = e^{-\alpha_i}$, $0\leq i\leq \ell$, since
$\Pi^{\Lambda}\sub \{\mu\in \hat P\st \mu\leq\Lambda\}$. The character \index{character} of the module is the shifted power series $ch(V^{\Lambda}) =e^{\Lambda}gr(V^{\Lambda})$.
\end{dfn}

\begin{thm}\index{Weyl-Kac!character formula}\label{def:weylkacchar} For any $\Lambda\in {\hat P}^+$ we have the following Weyl-Kac character formula for the character of irreducible finite dimensional highest weight $\ghat$-module $V^\Lambda$:
$$ ch(V^{\Lambda})=\frac{S_{\Lambda+\hat\rho}}{S_{\hat\rho}}\quad \hbox{ where } \quad
S_{\mu}= \sum_{w\in\hat\W} sgn(w) e^{w\mu-\hat\rho}\quad\hbox{ and }\quad \hat\rho = \sum_{i=0}^\ell \Lambda_i.$$ 
The \index{Weyl-Kac!denominator formula} Weyl-Kac denominator formula is:
\be\label{weylkacden}S_{\hat\rho} = \prod_{\alpha\in\hat\Phi^+} (1-e^{-\alpha})^{mult(\alpha)}. \ee
Note that setting $u_i = e^{-\alpha_i}$ for $0\leq i\leq \ell$, both $S_{\hat\rho}$ and 
$e^{-\Lambda}\ ch(V^{\Lambda}) = gr(V^{\Lambda})$ are in the power series ring $\Z[[u_0,\cdots,u_\ell]]$. 
\end{thm}

\begin{dfn}\index{Witt algebra}\label{def:Witt}
The \emph{Witt algebra} $\mathcal D$ is the infinite dimensional Lie algebra with basis $\{d_m\ | \ m\in\Z\}$ and
brackets $[d_m,d_n] = (n-m) d_{m+n}$, for $m,n\in \Z$. This Lie algebra has a representation on the ring of Laurent polynomials, $\C[t,t\inv]$, where the action is given by the formula $d_m=t^{m+1}\tfrac{d}{dt}$.  
This action extends to an action on $\ghat$ by $d_m\cdot x(n) = nx(m+n)$ so that $d_0 = d$ as in Definition \ref{def:KacMoody}. \index{Virasoro!algebra} The \emph{Virasoro algebra}, $\Vir$, is the central extension of $\mathcal D$, with basis $\{L_m, c_{Vir}\ | \ m\in\Z\}$ , and brackets \begin{align}\label{virbrackets} [\Vir,c_{Vir}] & = 0 \hbox{ and } \notag  \\
 [L_m, L_n] & =  (m-n)L_{m+n}+\tfrac{1}{12}(m^3-m)\dm c_{Vir} \mbox{ for } m,n\in\Z.\end{align}
\end{dfn}

\begin{rmk} \label{rmk:Lisminusd} Note that $[-d_m,-d_n] = (m-n)(-d_{m+n})$ so the projection from $\Vir$ to 
$\mathcal D$ sends $L_m$ to $-d_m$. 
\end{rmk}

\begin{thm}\index{Virasoro!representations}\label{def:Virasororepresentations} For each $(c,h)\in\C^2$ there is an irreducible representation of the Virasoro algebra (an irreducible $\Vir$-module), denoted by $\Vir(c,h)$ or $L(c,h)$, such that the central element $c_{Vir}$ acts by the scalar $c$, called the \emph{central charge} of the module, and such that there is a highest weight vector $v^+$ satisfying $L_m(v^+) = 0$ for $m>0$ and $L_0(v^+) = h v^+$. For certain values of $c$ and $h$ the $\Vir$-module $\Vir(c,h)$\index{module!Virasoro} admits a positive definite Hermitian form such that 
$(L_m(u),v) = (u, L_{-m}(v))$ for any $u,v\in \Vir(c,h)$, in which case the module is called \emph{unitary}. A very important class of these are the \emph{discrete series} for which the central charge $0 < c < 1$ can only have one of the discrete set of values $c = 1 - \frac{6}{s(s+1)}$ for $3 \leq s\in\Z$, and for which there are only a finite set of $h$ given by the special formula $h = h^{m,n} = \frac{[(s+1)m-sn]^2 - 1}{4s(s+1)}$, where $1\leq m\leq n < s+1$. The graded dimension for each of these $\Vir$-modules, as well as those of a larger class of $\Vir$-modules called the \emph{minimal models}, is given by a theorem of Feigin-Fuchs \cite{Feigin}, and will be used in this dissertation.
\end{thm}

\begin{scon}\index{Sugawara construction}\label{con:Sugawara} Let $V^\Lambda$ be a irreducible highest weight $\ghat$-module. 
Let $\{u_i\ |\ 1\leq i\leq dim(\g)\}$ be a basis of the finite dimensional simple Lie algebra $\g$, and let $\{u^i\ |\ 1\leq i\leq dim(\g)\}$ be the dual basis with respect to the normalized invariant symmetric bilinear form on $\g$ used in the bracket formula in Definition \ref{def:KacMoody}. 
Let $Cox^\vee(\g) = 1 + (\theta,\rho)$ be the \emph{dual Coxeter number} of $\g$ \cite{Kac}. Then the following 
\emph{Sugawara operators} represent $\Vir$ on $V^\Lambda$:
\be\label{sugstart}L_m = \frac{1}{2(Cox^\vee(\g) + \Lambda(c))} \ \sum_{k\in\Z} \no u_i(-k) u^i(m+k)\no \qquad\hbox{for}\qquad m\in\Z\ee
where the colons, $\no \ \ \no$ indicate \emph{bosonic normal ordering}\index{normal order!bosonic} of these operators on $V^\Lambda$, that is, 
$$\no u_i(-k) u^i(m+k)\no \ =  u_i(-k) u^i(m+k)\quad \hbox{ if }\quad -k \leq m+k$$ 
but 
$$\no u_i(-k) u^i(m+k)\no \  =  u^i(m+k) u_i(-k) \quad \hbox{ if }\quad m+k < -k.$$
In this $\Vir$ representation the central charge has value
\be\label{charge}\index{central charge} c =  \frac{dim(\g)\cdot \Lambda(c)}{Cox^\vee(\g) + \Lambda(c)} .\ee
\end{scon}

\begin{ccon}\index{spinor construction!affine}\label{con:AffineSpinorConstruction} 
Let $Z=\Z+\epsilon\in \left\{\Z, \hz\right\}$. Let $A = A^+\oplus A^-$ where $A^\pm\cong\C^{\ell}$ are maximal isotropic subspaces with respect to a non-degenerate symmetric bilinear form on $A$. Choose \emph{canonical bases} 
$\{a_1,\cdots,a_\ell\}$ and $\{a^*_1,\cdots,a^*_\ell\}$ for $A^+$ and $A^-$, respectively, such that 
$(a_i,a_j) = 0 = (a^*_i,a^*_j)$ and $(a_i,a^*_j) = \delta_{i,j}$. 
Define $A(Z)= A\otimes t^\epsilon\C[t,t\inv]$, with elements written $a(n)= a\otimes t^n$, where $a\in A$ and $n\in Z$, equipped with the symmetric bilinear form 
$$(a(m),b(n))=(a,b)\dm \ph \textrm{ for } \ph a,b\in A, \ph m,n\in Z .$$
Let $\cl(Z) = \Cl(A(Z),(\cdot,\cdot))$ be the infinite dimensional Clifford algebra \index{Clifford algebra!infinite dimensional} whose elements satisfy the  relations 
$$a(m)b(n)+b(n)a(m)=\big(a(m),b(n)\big)1.$$ \index{Clifford algebra!main relation}
When $\epsilon = \frac{1}{2}$, define the maximally isotropic subspaces
\begin{align*}
A(Z)^+ :=& \mbox{ span }\{a(n)\in A(Z)\st n>0\}\\
A(Z)^- :=& \mbox{ span }\{a(n)\in A(Z)\st n<0\}
\end{align*}
and when $\epsilon = 0$, define
\begin{align*}
A(Z)^+ :=& \mbox{ span }\{a(n)\in A(Z)\st n>0\}\cup  \mbox{ span }\{a(0)\st a\in A^+\}\\
A(Z)^- :=& \mbox{ span }\{a(n)\in A(Z)\st n<0\}\cup  \mbox{ span }\{a(0)\st a\in A^-\}.
\end{align*}

Let $\I(Z)$ be the left ideal in $\cl(Z)$ generated by $A(Z)^+$ so the quotient 
\be \cm(Z)=\cl(Z)/\I(Z)\ee 
is an irreducible left $\cl(Z)$-module with \emph{vacuum vector}\index{vacuum vector} 
\be \bvz=1+\I(Z) \ee
such that $A(Z)^+\cdot \bvz=0$. We have
\be \cm(Z) =  \mbox{ span } \{b_1(-n_1)\cdots b_r(-n_r)\bvz\st b_i(-n_i)\in A(Z)^-\}, \ee 
and is $\frac{1}{2}\mathbb Z$ graded, where the grading is given by the absolute value of the sum of the mode numbers, $\ds\sum_{i=1}^{r} n_i$. There is a direct sum decomposition $\cm(Z)=\cm^0(Z)\oplus \cm^1(Z)$ into even and odd subspaces, depending on the parity of $r$.
Define the \emph{fermionic normal ordering}\index{normal order!fermionic} \be\fno\label{afno} a(m)b(n)\fno = \left\{\begin{array}{ll}
a(m)b(n)& \textrm{ for } m<n\\
1/2\Big(a(m)b(n)-b(n)a(m)\Big)& \textrm{ for } m=n\\
-b(n)a(m)& \textrm{ for } m>n\\
\end{array}
\right.\ee
so that $\fno a(m)b(n)\fno = - \fno b(n)a(m)\fno$. 

When $Z=\hz$ the modules $\cm(Z)=\cm^0(Z)\oplus \cm^1(Z)$ are called Neveu-Schwarz\index{module!Neveu-Schwarz}, and when $Z=\Z$ they are called Ramond\index{module!Ramond}. The decomposition of these two Clifford modules provide the spinor construction of four level-$1$ representations for the affine Kac-Moody Lie algebra $\ghat$ of type $D_\ell^{(1)}$. Adopt the notation,
\be \label{dellmodules}
\hV^0 = \cm^0(\hz), \quad \hV^1 = \cm^1(\hz),\quad \hV^2 = \cm^0(\Z), \quad \hV^3 = \cm^1(\Z),
\ee
and let  $\hV^i_n$  denote the $n^{th}$ graded subspace of $\hV^i$.
The highest weights of these four modules are 
\be \label{dellhws} 
\Lambda_0, \quad \Lambda_1 - \one\delta,\quad \Lambda_4  - \tfrac{\ell}{8} \delta, \quad \Lambda_3  - \tfrac{\ell}{8} \delta, \ee
respectively, corresponding to the four endpoints of the Dynkin diagram, and with highest weight vectors 
\be \label{d4hwvs}
\bv = \bv(\hz), \quad a_1(-\one)\bv, \quad \bv' = \bv(\Z), \quad a^*_\ell(0)\bv'.
\ee
The best way to express the operators representing $\ghat$ and $\Vir$ on these modules is through generating functions. 
Let $w$ be a formal variable and for $a\in A$ write the generating function of Clifford operators
\begin{align}\label{gena}a(w):=& \sum_{m\in Z}a(m)w^{-m-1/2}\end{align} 
and its derivative generating function
\begin{align}\label{genader}a^{(1)}(w):=& \sum_{m\in Z}(-m-\one)a(m)w^{-m-3/2}.\end{align} 
Using the fermionic normal ordering, we can define the following generating functions of operators
\begin{align}\label{genab} 
\fno a(w)b(w)\fno &= \sum_{m\in Z}\sum_{n\in Z}\fno a(m)b(n) \fno w^{-n-m-1}  \notag \\
&= \sum_{k\in \mathbb Z}\left(\sum_{n\in Z}\fno a(k-n)b(n)\fno \right)w^{-k-1}\end{align}
and
\begin{align}\label{genaderb} 
\fno a^{(1)}(w)b(w)\fno &= \sum_{m\in Z}\sum_{n\in Z}(-m-\one)\fno a(m)b(n) \fno w^{-n-m-2}  \notag \\
&= \sum_{k\in \Z}\left(\sum_{n\in Z}(n-k-\one)\fno a(k-n)b(n)\fno \right)w^{-k-2}.\end{align}
 For each $k\in \Z$, we will use the following notations for the inner summations in \eqref{genab} and
\eqref{genaderb}, respectively,
\be\label{normk} \ds\sum_{n\in Z}\fno a(k-n)b(n)\fno  \ph = \ph\fno a(w)b(w)\fno _k
 \ee
and 
\be\label{normderk} \ds\sum_{n\in Z}(n-k-\one)\fno a(k-n)b(n)\fno  \ph = \ph\fno a^{(1)}(w)b(w)\fno _k.
 \ee
Because of the normal ordering, both $\fno a(w)b(w)\fno_k $ and $\fno a^{(1)}(w)b(w)\fno_k$ are well-defined operators on $\cmz$ for all $k\in\Z$, providing a representation of $\ghat$ and $\Vir$
on $\hV = \hV^0 \oplus \hV^1 \oplus \hV^2 \oplus \hV^3$ as follows. Let the operators $L^Z_k$ be defined by their generating function
\begin{align}\label{genVirdef} 
L^Z(w) &= \sum_{k\in \Z} L^Z_k w^{-k-2}\notag \\
&= \frac{1+(-1)^{2\epsilon}}{16} \ell w^{-2}\ + \frac{1}{2}\ \sum_{i=1}^\ell \fno a_i^{(1)}(w)a_i^*(w) + a_i^{*(1)}(w)a_i(w)\fno .\end{align}
This means that for $0\neq k\in\Z$, 
\be\label{Virnonzero} L^Z_k =  \sum_{i=1}^\ell \sum_{n\in Z} (n-\one k)\fno a_i(k-n)a_i^*(n)\fno \ee
and
\be\label{VirzeroNS} L^{\hz}_0 =  \sum_{i=1}^\ell \sum_{n\in \hz} n \fno a_i(-n)a_i^*(n)\fno \ee
but 
\be\label{VirzeroR} L^{\Z}_0 = \frac{\ell}{8} + \sum_{i=1}^\ell \sum_{n\in \Z} n \fno a_i(-n)a_i^*(n)\fno .\ee
\end{ccon}

\begin{prop}\label{prop:rep}\cite{FFR} Let $\ghat$ be the untwisted affine Kac-Moody Lie algebra of type $D_\ell^{(1)}$. For $Z\in \Z$ or $\hz$ the operators \eqref{normk} represent the elements
$\fno ab\fno(k)\in\ghat$ for $k\in\Z$, the identity operator represents $c\in\ghat$, the operator $- L^Z_0$ given in \eqref{VirzeroNS} and \eqref{VirzeroR} represents $d \in\ghat$, and the operators \eqref{Virnonzero}, \eqref{VirzeroNS} and \eqref{VirzeroR}
represent $L_k\in\Vir$ on $\cm(Z)$ with $c_{Vir}$ acting as the scalar $\ell$. $\cm(Z)=\cm^0(Z)\oplus \cm^1(Z)$ is the decomposition of $\cm(Z)$ into two irreducible $\ghat$-modules, so that \eqref{dellmodules} are the four level-1 $\ghat$-modules with highest weights given in \eqref{dellhws}. 
\end{prop}

\begin{rmk}
 \label{rmk:VO} These operators are examples of \emph{vertex operators}\index{vertex operators}, 
\be\label{VO} Y(v,w) = \sum_{m\in\one\Z} Y_m(v) w^{-m-wt(v)} \ee
where $\ds v\in V = \bigoplus_{r\in\one\Z} V_r$ is a $\one\Z$-graded vector space, $v\in V_r$ is said to have weight
$r = wt(v)$, $w$ is a formal variable and $Y_m(v)$ is an operator on $V$. 
That vector space is $\cm(\hz)$ where the grading is given by the eigenvalues of $L_0^{\hz}$. In particular, 
\eqref{gena} - \eqref{genaderb} with $Z = \hz$ are the vertex operators
\begin{align}
\label{VOa} Y\Big(a(-\one)\bv,w\Big) &= a(w) \\
\label{VOader} Y\Big(a(-\tre)\bv,w\Big) &= a^{(1)}(w) \\
\label{VOab} Y\Big(a(-\one)b(-\one)\bv,w\Big) &= \fno a(w)b(w)\fno \\
\label{VOaderb} Y\Big(a(-\tre)b(-\one)\bv,w\Big) &= \fno a^{(1)}(w)b(w)\fno 
\end{align}
and with \emph{conformal vector} \index{conformal vector} 
\be\label{omegadell}
\omega_{D_\ell} = \one \sum_{i=1}^\ell \Big(a_i(-\tre) a_i^*(-\one) + a_i^*(-\tre)a_i(-\one)\Big)\bv
\ee
\eqref{genVirdef} with $Z = \hz$ is the vertex operator
\be\label{VOVir}
Y(\omega_{D_\ell},w) = L^{\hz}(w) .\ee
Furthermore, the operators in \eqref{normk} and \eqref{normderk} are 
\begin{align}
 \label{VOabk} Y_k\Big(a(-\one)b(-\one)\bv\Big) &= \fno a(w)b(w)\fno_k \\
 \label{VOaderbk} Y_k\Big(a(-\tre)b(-\one)\bv\Big) &= \fno a^{(1)}(w)b(w)\fno_k .
\end{align}
In \cite{FFR} it is shown that $\cm(\hz)$ has the structure of a \emph{vertex operator superalgebra} (VOSA)
with even part $\hV^0$ and odd part $\hV^1$. In order to give $\cm(\Z)$ the structure of a module for that VOSA, 
we first define unmodified vertex operators ${\overline Y}(v,w)$ as in the definitions 
\eqref{VOa} - \eqref{VOaderbk} with $Z = \Z$. 
Define the Ramond module vertex operators to be
\begin{align}\label{RamondVO} Y(v,w) = {\overline Y}(\exp(\Delta(w))v,w)  \end{align}  
where \begin{align} \Delta(w)= \sum_{1 \leq i \leq \ell} \sum_{0\leq m, n\in\Z }C_{m,n}a_i(m+\one)a_i^*(n+\one) 
w^{-m-n-1}\end{align} 
is a generating function of operators on $\cm(\hz)$ with coefficients  
\begin{align}C_{m,n} = \tfrac{1}{2}\frac{m-n}{(m+n+1)} {-\frac{1}{2}\choose \phm m} {-\frac{1}{2}\choose\phm n}.\end{align}  
These vertex operators provide the Ramond module $\cm(\Z)$ with the structure of a vertex operator superalgebra module, as shown in \cite{FFR}. For $v\in\hV^i_n$ with $i\in\{0,1\}$ and $n\in\{0,1/2,1,3/2\}$, 
$\exp(\Delta(w))v = v$ so $Y(v,w) = {\overline Y}(v,w)$. For $n\geq 2$, we may have $\exp(\Delta(w))v \neq v$, 
and, in particular, for $v = \omega_{D_\ell}$ we have 
\be\label{Deltadell} \exp(\Delta(w))\omega_{D_\ell} = \omega_{D_\ell} + \tfrac{1}{8} \ell w^{-2} \bv\ee
(note correction to a misprint in \cite{FFR}, equation (3.62)) so 
\be\label{RamondVirasoroVO} Y(\omega_{D_\ell},w) = {\overline Y}(\omega_{D_\ell},w) 
+ \tfrac{1}{8} \ell w^{-2} {\overline Y}(\bv,w) 
= {\overline Y}(\omega_{D_\ell},w) + \tfrac{1}{8} \ell w^{-2} I = L^{\Z}(w) \ee
where $I$ is the identity operator on $\cm(\Z)$. The effect of this ``correction'' is to add the scalar
$\tfrac{1}{8} \ell$ to the unmodified operator ${\bar L}^\Z_0$. We will not need to deal with cases where $n > 2$. 
The relations among vertex operators established in \cite{FFR} will be used here without proof. 

The Sugawara construction for $\ghat$ on each irreducible representation $\hV^i$, $0\leq i\leq 3$, simplifies to
the vertex operator $Y(\omega_{D_\ell},w)$, which is actually a Sugawara construction for the Heisenberg subalgebra $\hhat$ of $\ghat$. For other subalgebras of $\ghat$ such simplification does not occur, so the full Sugawara construction will be needed in later chapters, where we will focus on the case $\ell = 4$. For $\ghat$ of type 
$D_4^{(1)}$ we have triality, and $\ghat$ contains subalgebras of types $B_3^{(1)}$ and $G_2^{(1)}$. The associated full Sugawara constructions play a vital role in this investigation.
\end{rmk}

To close this introductory chapter on background information we present the elegant theorem of \cite{GKO} on the construction of \emph{commuting coset Virasoro algebras} \index{Virasoro!coset}. 
\begin{thm}\cite{GKO}\label{thm:cosetvir}
Let  $\k$ be a Lie subalgebra of a finite dimensional Lie algebra $\g$, let $V^\Lambda$ be an irreducible representation of $\ghat$, 
and let $\Vir^{\k}$ and $\Vir^{\g}$ be the representations of the Virasoro algebra on $V^\Lambda$ provided by the Sugawara construction for
$\khat$ and $\ghat$, respectively, given by the operators $\{L_m^{\k},c_{Vir}^{\k}\ |\ m\in\Z\}$ and $\{L_m^{\g},c_{Vir}^{\g}\ |\ m\in\Z\}$.  Then the differences $L_m^{\g} - L_m^{\k}$ also provide a representation of the Virasoro algebra on $V^\Lambda$ with central charge 
$c_{Vir}^{\g-\k} = c_{Vir}^{\g}-c_{Vir}^{\k}$. Moreover, this \emph{coset Virasoro representation}, $\Vir^{\g-\k}$, commutes with $\khat$ and with 
$\Vir^{\k}$, that is, $[L_m^{\g} - L_m^{\k}, \khat] = 0$ and $[L_m^{\g} - L_m^{\k}, L_n^{\k}] = 0$, for all $m,n\in\Z$. 
\end{thm}

We will use this result to explicitly find two such coset Virasoro representations, one as we step down from $D_4^{(1)}$ to $B_3^{(1)}$, 
$\Vir^{D_4-B_3}$, and another,  $\Vir^{B_3-G_2}$,  from $B_3^{(1)}$ to $G_2^{(1)}$.  These coset Virasoros have central charges of 
$\one$ and $\tfrac{7}{10}$, respectively, and their representations on $\hV$ are the key to the branching rules that we wish to expose.

\chapter {Graded Dimensions for Clifford, Virasoro, and $G_2^{(1)}$-Modules }\label{cha:Chars}
\section{Clifford Modules}\label{sec:Cliff}

Infinite dimensional Clifford modules play a central role in the spinor construction of level-$1$ modules for the affine Kac-Moody
Lie algebras of type $D_\ell^{(1)}$. In this investigation we will only be concerned with the case $\ell=4$, where triality is present. In Chapter \ref{cha:constructs} 
we will see how the four irreducible $D_4^{(1)}$ modules are constructed as the even and odd subspaces of two kinds of Clifford modules, the Neveu-Schwarz and the Ramond, 
\be \label{D4modules}
\hV^0 = \CM_4^0(\hz), \quad \hV^1 = \CM_4^1(\hz),\quad \hV^2 = \CM_4^0(\Z), \quad \hV^3 = \CM_4^1(\Z).
\ee
with highest weight vectors 
\be \label{D4hwvs}
\bv = \bv(\hz), \quad a_1(-\one)\bv, \quad \bv' = \bv(\Z), \quad a^*_4(0)\bv'.
\ee 
The highest weights of these four vectors, and their modules, are the fundamental weights 
$\Lambda_0$, $\Lambda_1 - \one\delta$, $\Lambda_4  - \one\delta$, $\Lambda_3  - \one\delta$, respectively,
corresponding to the four endpoints of the Dynkin diagram, but the last three are shifted by $-\one\delta$ because 
the derivation $-d$ is represented by the Virasoro operator $L_0$. 

The homogeneous graded dimension\index{graded dimension!homogeneous} of each Clifford module is easily written in a product form as follows. 
\be \label{NSgrdim}
gr(\CM_4(\hz)) = \prod_{i=1}^4 \prod_{0<n\in\hz} (1 + e^{\epsilon_i} q^n) (1 + e^{-\epsilon_i} q^n) ,
\ee
\be \label{Rgrdim}
gr(\CM_4(\Z)) =  \left(\prod_{i=1}^4 (1 + e^{-\epsilon_i})\right)  
\prod_{i=1}^4 \prod_{0<n\in\Z} (1 + e^{\epsilon_i} q^n) (1 + e^{-\epsilon_i} q^n) ,
\ee
where $\epsilon_i$ is the weight of $a_i$ in the natural representation of $D_4$, $q = e^{-\delta}$ and $\delta$ is the 
minimal positive imaginary root of $D_4^{(1)}$. We will use the following notations: 
\be \label{D4abbrev}
v_i = e^{\epsilon_i}, \qquad u_i = e^{-\alpha_i}, \quad 1\leq i\leq 4, \qquad u_0 = e^{-\alpha_0} = e^{\theta - \delta} 
\ee
where $\theta = \alpha_1+ 2\alpha_2 + \alpha_3 + \alpha_4 = \epsilon_1 + \epsilon_2$ is the highest root of $D_4$. The homogeneous characters of these Clifford modules are then just shifted according to their highest weights, 
\be \label{NSRchars}
ch(\CM_4(\hz)) = e^{\Lambda_0}\ gr(\CM_4(\hz)) ,\qquad ch(\CM_4(\Z)) = e^{\Lambda_4-\delta/2}\ gr(\CM_4(\Z)).
\ee
In order to relate these to each other it is necessary to know that 
\be \label{Topshift}
\Lambda_4 - \Lambda_0 = \lambda_4 = \one(\epsilon_1 + \epsilon_2 + \epsilon_3 + \epsilon_4).
\ee
In the principal specialization \index{specialization!principal} of the graded dimension\index{graded dimension!principal} and character we set all $u_i = u$ for $0\leq i\leq 4$, which means
$q = u_0 e^{-\theta} =  u_0 u_1 u_2^2 u_3 u_4 = u^6$ and from the relations between the simple roots and the $\epsilon_i$,
we have $v_i = u^{i-4}$ for $1\leq i\leq 4$. This gives the following results for the principal graded dimension: 
\bea \label{NSprgrdim}
gr_{pr}(\CM_4(\hz)) &=&  \prod_{0\leq m\in\Z}  \prod_{i=1}^4 (1 + v_i q^{m+\one}) (1 + v_i^{-1} q^{m+\one}) \notag \\
&=& \prod_{0\leq m\in\Z}  \prod_{i=1}^4 (1 + u^{i-4} u^{6m+3}) (1 + u^{4-i} q^{6m+3}) \notag \\
&=& \prod_{0\leq m\in\Z} (1+u^{6m}) (1+u^{6m+1}) (1+u^{6m+2}) (1+u^{6m+3}) \cdot \notag \\
&&\qquad            (1+u^{6m+6}) (1+u^{6m+5}) (1+u^{6m+4}) (1+u^{6m+3})  \notag \\
&=& 2 \prod_{1\leq n\in\Z} (1 + u^n) (1 + u^{3n}) \notag \\
&=& 2 \frac{\phi(u^2) \phi(u^6)}{\phi(u) \phi(u^3)}.
\eea
For the Ramond Clifford module the principal graded dimension is exactly the same:
\bea \label{Rprgrdim}
gr_{pr}(\CM_4(\Z)) &=& \left(\prod_{i=1}^4 (1 + v_i^{-1})\right)  
   \prod_{i=1}^4 \prod_{1\leq n\in\Z} (1 + v_i q^n) (1 + v_i^{-1} q^n) \notag \\
&=& \left(\prod_{i=1}^4 (1 + u^{4-i})\right)  
   \prod_{i=1}^4 \prod_{1\leq n\in\Z} (1 + u^{i-4} u^{6n}) (1 + u^{4-i} u^{6n}) \notag \\
&=&(1+u^3)(1+u^2)(1+u)(2) \prod_{1\leq n} (1+u^{6n-3}) (1+u^{6n-2}) (1+u^{6n-1})  \cdot  \notag \\
&&\qquad (1+u^{6n}) (1+u^{6n+3}) (1+u^{6n+2}) (1+u^{6n+1}) (1+u^{6n}) \notag \\
&=& 2 \prod_{1\leq n\in\Z} (1 + u^n) (1 + u^{3n}) \notag \\
&=& 2 \frac{\phi(u^2) \phi(u^6)}{\phi(u) \phi(u^3)}.
\eea

From (\ref{D4modules}) we see that
\be 
gr_{pr}(\hV^0) + gr_{pr}(\hV^1) = gr_{pr}(\CM_4(\hz)) 
\ee
and
\be
 gr_{pr}(\hV^2) + gr_{pr}(\hV^3) = gr_{pr}(\CM_4(\Z)) 
\ee
and from the triality symmetry we also must have
\be \label{triality1}
gr_{pr}(\hV^1) = gr_{pr}(\hV^2) = gr_{pr}(\hV^3) 
\ee
which implies 
\be \label{triality2}
gr_{pr}(\hV^0) = gr_{pr}(\hV^1) = gr_{pr}(\hV^2) = gr_{pr}(\hV^3) = \frac{\phi(u^2) \phi(u^6)}{\phi(u) \phi(u^3)}.
\ee

Although we do not need it here, there is another interesting specialization of the homogeneous character which, when
combined with the triality symmetry among $\hV^j$ for $1\leq j\leq 3$, gives a famous identity of Jacobi. 
If we use the specialization $v_i = 1$ for $1\leq j\leq 4$, then we get a ``horizontal graded dimension"\index{graded dimension!horizontal}
\be \label{NShorgrdim}
gr_{hor}(\CM_4(\hz)) = \prod_{0<n\in\hz} (1 + q^n)^8 ,
\ee
\be \label{Rhorgrdim}
gr_{hor}(\CM_4(\Z)) =  16  \prod_{0<n\in\Z} (1 + q^n)^8 .
\ee
The horizontal characters of the two Neveu-Schwarz modules, after replacing $q$ by $q^2$ to avoid dealing with half-integral
powers, and after multiplication by $e^{-\Lambda_0}$ to start the series at $1$, are then
\be \label{V0horchar}
e^{-\Lambda_0}\ ch_{hor}(\hV^0)(q^2) 
= \frac{1}{2} \left( \prod_{0\leq n\in\Z} (1 + q^{2n+1})^8 + \prod_{0\leq n\in\Z} (1 - q^{2n+1})^8 \right) ,
\ee
\be \label{V1horchar}
e^{-\Lambda_0}\ ch_{hor}(\hV^1)(q^2) 
= \frac{1}{2} \left( \prod_{0\leq n\in\Z} (1 + q^{2n+1})^8 - \prod_{0\leq n\in\Z} (1 - q^{2n+1})^8 \right) ,
\ee
but the shifted horizontal characters of the two Ramond modules, including the $(q^2)^{1/2}$ factor from the highest weight, are equal, each being half of the total, 
\be \label{V2V3horchar}
e^{-\Lambda_0}\ ch_{hor}(\hV^2)(q^2) = e^{-\Lambda_0}\ ch_{hor}(\hV^3)(q^2) = 8q \prod_{1\leq n\in\Z} (1 + q^{2n})^8  .
\ee
The equality $ch_{hor}(\hV^1)(q^2) = ch_{hor}(\hV^2)(q^2)$ then gives the Jacobi identity \index{Jacobi identity}
\be \label{JacobiAbsIden}
\prod_{0\leq n\in\Z} (1 + q^{2n+1})^8 - \prod_{0\leq n\in\Z} (1 - q^{2n+1})^8 = 16q \prod_{1\leq n\in\Z} (1 + q^{2n})^8.
\ee

\section{Virasoro Modules}\label{sec:vir}

The Virasoro algebra, $\Vir$, is the infinite dimensional Lie algebra with basis $\{L_m, C\ |\ m\in\Z\}$ and brackets
$$[L_m,L_n] = (m-n) L_{m+n} + \frac{m^3 - m}{12} \delta_{m,-n} C\qquad [L_m, C] = 0.$$
Representations of the Virasoro algebra, in particular the highest weight modules, are closely connected to the theory of affine Kac-Moody Lie algebras, and to the theory of vertex operator algebras. The class of ``minimal models", and their graded dimensions (characters) play an especially important role in conformal field theory, string theory and vertex algebra theory. Each irreducible highest weight $\Vir$-module, $\Vir(c,h)$ is uniquely determined by an ordered pair of complex numbers, the central charge $c$ being the eigenvalue of the central element $C$, and $h$ being the eigenvalue of $L_0$ on the highest weight vector. The module $\Vir(c,h)$ decomposes into a direct sum of finite dimensional eigenspaces for $L_0$, and the graded dimension (also called the character, although it does not actually characterize the module uniquely) is the formal power series
$$gr(c,h) = \sum_{n\geq 0} dim(\Vir(c,h)_{h-n}) q^n.$$
For some purposes it is important to multiply this series by a (possibly) fractional power of $q$, so one defines the ``character"
$$\chi(c,h) = q^{h-c/24} gr(c,h).$$
The minimal models \index{minimal models} have particular values of the central charge \index{central charge} $c < 1$, and for each one a finite number of $h$ values, as follows. For $2\leq s,t\in\Z$ relatively prime and $1\leq m < s$, $1\leq n < t$, let
\be \label{cst}
c_{s,t}=1-\frac{6(s-t)^2}{st}, \ \ h_{s,t}^{m,n}=\frac{(mt-ns)^2-(s-t)^2}{4 st}.
\ee
Then Feigin-Fuchs proved that the character $\chi_{s,t}^{m,n}(q)$ of the module $\Vir(c_{s,t},h_{s,t}^{m,n})$ is given by 
\be \label{VirChar}
\chi_{s,t}^{m,n}(q) 
= \frac{q^{(h_{s,t}^{m,n}-c_{s,t}/24)}}{\phi(q)} \cdot \sum_{k \in \Z} q^{st k^2}(q^{k(mt-ns)}-q^{(mt+ns)k+mn})
\ee 
where $\phi(q) = \prod_{i= 1}^\infty (1 - q^i)$. 
Some of these characters have a product form, for example, the following two characters with $c_{2,5} = \frac{-22}{5}$, 
\bea \label{rr1}
&& \chi_{2,5}^{1,1}(q)=q^{11/60} \prod_{n =0}^\infty \frac{1}{(1-q^{5n+2})(1-q^{5n+3})}, \\
\label{rr2}
&& \chi_{2,5}^{1,2}(q)=q^{-1/60} \prod_{n = 0}^\infty \frac{1}{(1-q^{5n+1})(1-q^{5n+4})}
\eea
contain the product sides of the famous \index{Rogers-Ramanujan series} Rogers-Ramanujan series, which appear 
in many papers on the representation theory of Virasoro and affine Lie algebras.

One way to see this is through the Jacobi Triple Product Identity \index{Jacobi identity! triple product}(JTPI), which is the Weyl-Kac denominator formula\index{Weyl-Kac!denominator formula}
for the affine algebra $A_1^{(1)}$, 
\be \label{JTPI}
\prod_{n\geq 1} (1 - u^n v^n) (1 - u^n v^{n-1}) (1 - u^{n-1} v^n) = \sum_{k\in\Z} (-1)^k u^{k(k+1)/2} v^{k(k-1)/2}
\ee
whose specializations yield the following identities
\be \label{JTPSpec1}
\phi(q) = \prod_{n\geq 1} (1 - q^n) = \sum_{k\in\Z} (-1)^k q^{k(3k-1)/2} \quad \hbox{with } u=q,\ \ v=q^2,
\ee
\be \label{JTPSpec2}
\prod_{n\geq 1} (1 - q^{5n}) (1 - q^{5n-1}) (1 - q^{5n-4})  = \sum_{k\in\Z} (-1)^k q^{k(5k+3)/2}  \quad \hbox{with } u=q^4,\ \ v=q,
\ee
\be \label{JTPSpec3}
\prod_{n\geq 1} (1 - q^{5n}) (1 - q^{5n-2}) (1 - q^{5n-3})  = \sum_{k\in\Z} (-1)^k q^{k(5k+1)/2}  \quad \hbox{with } u=q^3,\ \ v=q^2.
\ee
We record here for future reference the following notations and formulas for these series:
\be \label{RR1}
a(q) = \prod_{n\geq 1} \frac{1}{(1 - q^{5n-2}) (1 - q^{5n-3})} = \frac{1}{\phi(q)}\ \sum_{k\in\Z} (-1)^k q^{k(5k+3)/2} ,
\ee
\be \label{RR2}
b(q) = \prod_{n\geq 1} \frac{1}{(1 - q^{5n-1}) (1 - q^{5n-4})} = \frac{1}{\phi(q)}\ \sum_{k\in\Z} (-1)^k q^{k(5k+1)/2} .
\ee
Note that the summation in (\ref{VirChar}) for $(s,t) = (2,5)$ and $(m,n) = (1,1)$ is
\be \label{VirChar2511}
\sum_{k \in \Z} q^{10 k^2}(q^{3k}-q^{7k+1}) = \sum_{j\in\Z} (-1)^j q^{j(5j+3)/2}
\ee 
and for $(m,n) = (1,2)$ it is
\be \label{VirChar2512}
\sum_{k \in \Z} q^{10 k^2}(q^{k}-q^{9k+2}) = \sum_{j\in\Z} (-1)^j q^{j(5j+1)/2}
\ee 
using $j = 2k$ to get the terms with coefficient $1$ and $j = -2k-1$ to get those with coefficient $-1$. 
These yield (\ref{rr1}) and (\ref{rr2}) by using (\ref{RR1}) and (\ref{RR2}). 

We will be concerned with the characters of Virasoro modules with central charge $c = c_{3,4} = \frac{1}{2}$,
which has three distinct $h$ values, $h_{3,4}^{1,1} = 0$, $h_{3,4}^{1,3} = \frac{1}{2}$ and $h_{3,4}^{1,2} = \frac{1}{16}$.
These three characters are:
\be \label{VirChar3411}
\chi_{3,4}^{1,1}(q) 
= \frac{q^{(-\frac{1}{48})}}{\phi(q)} \cdot \sum_{k \in \Z} q^{12 k^2}(q^{k}-q^{7k+1})
\ee 
\be \label{VirChar3413}
\chi_{3,4}^{1,3}(q) 
= \frac{q^{(\frac{1}{2}-\frac{1}{48})}}{\phi(q)} \cdot \sum_{k \in \Z} q^{12 k^2}(q^{-5k}-q^{13k+3})
\ee 
\be \label{VirChar3412}
\chi_{3,4}^{1,2}(q) 
= \frac{q^{(\frac{1}{16}-\frac{1}{48})}}{\phi(q)} \cdot \sum_{k \in \Z} q^{12 k^2}(q^{-2k}-q^{10k+2}).
\ee 
Note that, replacing $q$ by $q^2$ and then multiplying by $q^{\frac{1}{24}}$ we find the power series 
\bea \label{VirChar34diff}
&&[\chi_{3,4}^{1,1}(q^2) - \chi_{3,4}^{1,3}(q^2)] q^{\frac{1}{24}} \notag \\
&=& \frac{1}{\phi(q^2)} \cdot \sum_{k \in \Z} (q^{24 k^2+2k} - q^{24 k^2+14k+2} - q^{24 k^2-10k+1} + q^{24 k^2+26k+7} )\\
&=& \frac{1}{\phi(q^2)} \cdot \sum_{n \in \Z} (-1)^n q^{n(3n+1)/2} \notag \\
&=& \frac{\phi(q)}{\phi(q^2)} = \prod_{m\geq 1} (1 - q^{2m-1}) \notag
\eea
using $n = 4k$, $n = 4k+1$, $n = 4k-1$ and $n = 4k+2$ to combine the four kinds of terms, and then using (\ref{JTPSpec1}) 
to replace the summation by $\phi(q)$. This is a ``virtual" graded dimension, that is, a linear combination of graded dimensions,
with a shifted power of $q$ to account for the ``conformal dimension", $q^{2h}$. Since the even powers of $q$ all came from 
$\chi_{3,4}^{1,1}(q^2)$ and the odd powers of $q$ all came from $\chi_{3,4}^{1,3}(q^2)$, we see that 
\be \label{VirChar34sum}
[\chi_{3,4}^{1,1}(q^2) + \chi_{3,4}^{1,3}(q^2)] q^{\frac{1}{24}} 
= \frac{\phi(-q)}{\phi(q^2)} = \prod_{m\geq 1} (1 + q^{2m-1}) = \frac{\phi(q^2)^2}{\phi(q)\phi(q^4)}.
\ee
For later use we define the notation 
\be 
v(q) = \frac{\phi(q)}{\phi(q^2)} = \prod_{m\geq 1} (1 - q^{2m-1}) 
\ee
so that 
\be 
v(-q) = \frac{\phi(-q)}{\phi(q^2)} = \prod_{m\geq 1} (1 + q^{2m-1}) 
\ee
and we may write
\bea
\chi_{3,4}^{1,1}(q^2) q^{\frac{1}{24}} &=& \one (v(-q) + v(q)) \label{ch3411} \\
\chi_{3,4}^{1,3}(q^2) q^{\frac{1}{24}} &=& \one (v(-q) - v(q))  \label{ch3413} .
\eea

Turning our attention now to the $h = \frac{1}{16}$ module, if we first replace $k$ by $-k$ in the summation in (\ref{VirChar3412}),
without replacing $q$ by its square, we can write
\bea
\chi_{3,4}^{1,2}(q) q^{\frac{-1}{24}}  \notag
&=& \frac{1}{\phi(q)} \cdot \sum_{k \in \Z} q^{12 k^2}(q^{2k}-q^{-10k+2}) \\
&=& \frac{1}{\phi(q)} \cdot \sum_{n \in \Z} (-1)^n q^{n(3n+1)} \notag \\
&=& \frac{\phi(q^2)}{\phi(q)} = \prod_{m\geq 1} (1 + q^{m}) = \frac{1}{v(q)} \label{VirChar3412Ram}
\eea 
using $n = 2k$ and $n = 2k-1$ to combine the two kinds of terms in the summation, and using (\ref{JTPSpec1}) with $q$
replaced by $q^2$. 

We will also be concerned with the characters of Virasoro modules with central charge $c = c_{4,5} = \frac{7}{10}$,
which has six distinct $h$ values, $h_{4,5}^{1,1} = 0$, $h_{4,5}^{1,2} = \frac{1}{10}$, $h_{4,5}^{1,3} = \frac{3}{5}$,
$h_{4,5}^{1,4} = \frac{3}{2}$, $h_{4,5}^{2,1} = \frac{7}{16}$, $h_{4,5}^{2,2} = \frac{3}{80}$.
These six characters are:
\be \label{VirChar4511}
\chi_{4,5}^{1,1}(q) 
= \frac{q^{(-\frac{7}{240})}}{\phi(q)} \cdot \sum_{k \in \Z} q^{20 k^2}(q^{k}-q^{9k+1})
\ee 
\be \label{VirChar4512}
\chi_{4,5}^{1,2}(q) 
= \frac{q^{(\frac{1}{10}-\frac{7}{240})}}{\phi(q)} \cdot \sum_{k \in \Z} q^{20 k^2}(q^{3k}-q^{13k+2})
\ee 
\be \label{VirChar4513}
\chi_{4,5}^{1,3}(q) 
= \frac{q^{(\frac{3}{5}-\frac{7}{240})}}{\phi(q)} \cdot \sum_{k \in \Z} q^{20 k^2}(q^{7k}-q^{17k+3})
\ee 
\be \label{VirChar4514}
\chi_{4,5}^{1,4}(q) 
= \frac{q^{(\frac{3}{2}-\frac{7}{240})}}{\phi(q)} \cdot \sum_{k \in \Z} q^{20 k^2}(q^{11k}-q^{21k+4})
\ee 
\be \label{VirChar4521}
\chi_{4,5}^{2,1}(q) 
= \frac{q^{(\frac{7}{16}-\frac{7}{240})}}{\phi(q)} \cdot \sum_{k \in \Z} q^{20 k^2}(q^{6k}-q^{14k+2})
\ee 
\be \label{VirChar4522}
\chi_{4,5}^{2,2}(q) 
= \frac{q^{(\frac{3}{80}-\frac{7}{240})}}{\phi(q)} \cdot \sum_{k \in \Z} q^{20 k^2}(q^{2k}-q^{18k+4}).
\ee 
As above, we find that the first four combine in pairs after replacing $q$ by $q^2$ and multiplying by $q^{\frac{7}{120}}$,
but the last two simplify by themselves. We have
\bea \label{VirChar45diff1114}
&&[\chi_{4,5}^{1,1}(q^2) - \chi_{4,5}^{1,4}(q^2)] q^{\frac{7}{120}} \notag \\
&=& \frac{1}{\phi(q^2)} \cdot \sum_{k \in \Z} (q^{40 k^2+2k} - q^{40 k^2+18k+2} - q^{40 k^2+22k+3} + q^{40 k^2+42k+11} )\\
&=& \frac{1}{\phi(q^2)} \cdot \sum_{n \in \Z} (-1)^n q^{n(5n+1)/2} \notag \\
&=& \frac{\phi(q)b(q)}{\phi(q^2)} = v(q) b(q) \notag
\eea
using $n = 4k$, $n = -4k-1$, $n = 4k+1$ and $n = 4k+2$ to combine the four kinds of terms in the summation, and using
(\ref{RR2}). But this also gives 
\be \label{VirChar45sum1114}
[\chi_{4,5}^{1,1}(q^2) + \chi_{4,5}^{1,4}(q^2)] q^{\frac{7}{120}}
= \frac{\phi(-q)b(-q)}{\phi(q^2)} = v(-q) b(-q) .
\ee
We may then write
\bea
\chi_{4,5}^{1,1}(q^2) q^{\frac{7}{120}} &=& \one(v(-q) b(-q) + v(q) b(q)) \label{ch4511} \\
\chi_{4,5}^{1,4}(q^2) q^{\frac{7}{120}} &=& \one(v(-q) b(-q) - v(q) b(q)) \label{ch4514} .
\eea
Similarly, we find 
\be \label{VirChar45diff1213}
[\chi_{4,5}^{1,2}(q^2) - \chi_{4,5}^{1,3}(q^2)] q^{\frac{-17}{120}}
= \frac{\phi(q)a(q)}{\phi(q^2)} = v(q) a(q) ,
\ee
\be \label{VirChar45sum1213}
[\chi_{4,5}^{1,2}(q^2) + \chi_{4,5}^{1,3}(q^2)] q^{\frac{-17}{120}}
= \frac{\phi(-q)a(-q)}{\phi(q^2)} = v(-q) a(-q) .
\ee
So we may also write
\bea
\chi_{4,5}^{1,2}(q^2) q^{\frac{-17}{120}} &=& \one(v(-q) a(-q) + v(q) a(q)) \label{ch4512} \\
\chi_{4,5}^{1,3}(q^2) q^{\frac{-17}{120}} &=& \one(v(-q) a(-q) - v(q) a(q)) \label{ch4513}.
\eea

Finally, for the last two characters we have
\be \label{VirChar4521a}
\chi_{4,5}^{2,1}(q) q^{\frac{-49}{120}}
= \frac{1}{\phi(q)} \cdot \sum_{k \in \Z} (-1)^n q^{n(5n+3)} =  \frac{\phi(q^2)a(q^2)}{\phi(q)} = \frac{a(q^2)}{v(q)},
\ee 
\be \label{VirChar4522a}
\chi_{4,5}^{2,2}(q) q^{\frac{-1}{120}}
= \frac{1}{\phi(q)} \cdot \sum_{k \in \Z} (-1)^n q^{n(5n+1)} =  \frac{\phi(q^2)b(q^2)}{\phi(q)} = \frac{b(q^2)}{v(q)}.
\ee

\section{$G_2^{(1)}$-Modules}\label{sec:G2}

We will be concerned with the characters (graded dimensions) of two level-$1$ irreducible highest weight modules for the affine 
Kac-Moody Lie algebra $G_2^{(1)}$. These modules are denoted by $W(\Omega_0)$ and $W(\Omega_2)$, 
whose highest weights, $\Omega_0$ and $\Omega_2$ are the fundamental weights for $G_2^{(1)}$ corresponding to the two endpoints of the Dynkin diagram. 

The Weyl-Kac character formula for any Kac-Moody algebra, $\g$, gives the homogeneous graded dimension of any highest weight module, $V^\Lambda$, with highest weight $\Lambda$, in the form of a ratio, 
\be \label{WeylKacChar}
\chi(V^\Lambda) = \frac{S_{\Lambda+\rho}}{S_\rho} \qquad\hbox{where}\qquad S_\mu = \sum_{w\in\W} sgn(w) e^{w\mu - \rho}
\ee
is the alternating summation over the Weyl group $\W$, and $\rho$ is the sum of the fundamental weights of the algebra. 
The Weyl-Kac denominator formula gives a product form for $S_\rho$, 
\be \label{WeylKacDenom}
S_\rho = \prod_{\alpha>0} (1 - e^{-\alpha})^{mult(\alpha)}
\ee
where the product is over the positive roots of the algebra, and $mult(\alpha)$ is the dimension of the $\alpha$ root space. 
All positive roots are non-negative integral linear combinations of the simple roots, $\alpha_0$, $\alpha_1$, ..., $\alpha_\ell$, 
and for $\Lambda$ any dominant integral weight, we have that $\Lambda + \rho - w(\Lambda + \rho)$ is a non-negative integral
linear combination of the simple roots. So if we define variables $u_i = e^{-\alpha_i}$ for $0\leq i\leq\ell$, then the normalized
character (graded dimension),
\be \label{grdim}
gr(V^\Lambda) = \chi(V^\Lambda) e^{-\Lambda} 
\ee
is in the power series ring $\Z[u_0,u_1,\cdots,u_\ell]$. 

It is a remarkable result of Lepowsky that the principal specialization of the numerator, $S_{\Lambda+\rho} \ e^{-\Lambda}$,  
where one sets all variables $u_i = u$, equals a corresponding non-principal specialization of the denominator $S^\vee_\rho$ 
for the dual Kac-Moody Lie algebra, $\g^\vee$, where one sets $u_i = u^{s_i}$ for $s_i = (\Lambda+\rho)(h_i)$. We denote such 
a specialization by ${S\!pec}_{(s_0,s_1,\cdots,s_\ell)}(S^\vee_\rho)$, and note that the principal specialization of the denominator of $\g$ is $S\!pec_{(1,1,\cdots,1)}(S_\rho)$. 
That means the principal specialization of $gr(V^\Lambda)$ is a power series in one variable having a product form
\be \label{grprdim}
gr_{pr}(V^\Lambda) = \frac{{S\!pec}_{(s_0,s_1,\cdots,s_\ell)}(S^\vee_\rho)}{S\!pec_{(1,1,\cdots,1)}(S_\rho)}.
\ee
Furthermore, the character of a direct sum of modules is the sum of the characters of the summands, and the character of a tensor product of modules is the product of the characters of the factors, so the same is true of principal specializations. 
Using the theory of $Z$-algebras developed by Lepowsky and Wilson, Mandia computed the principal graded dimensions of 
the two $G_2^{(1)}$ modules we need. 

\newpage

\begin{thm} \label{Mandia} \cite{Mandia} Letting
\be \label{G2Fock}
\F(u) = \frac{\phi(u^2) \phi(u^3)}{\phi(u) \phi(u^6)} = \prod_{n\geq 1} \frac{1}{(1 - u^{6n-5}) (1 - u^{6n-1})} = \frac{v(u^3)}{v(u)},
\ee
then we have
\be \label{G2modules}
gr_{pr}(W(\Omega_0)) = \F(u) a(u^3) \qquad\hbox{and}\qquad gr_{pr}(W(\Omega_2)) = \F(u) b(u^3).
\ee
\end{thm}

\chapter{Decompositions of $D_4^{(1)}$-Modules into a direct sum of $G_2^{(1)}$-Modules }\label{cha:punchline}
\section{Notations and Simplifications}\label{sec:simp}

Let
\be \label{AandB}
A(t) = \sum_{m\geq 0} A_m t^m \qquad\hbox{and}\qquad B(t) = \sum_{m\geq 0} B_m t^m
\ee
be the power series whose coefficients give the branching rules for the decomposition of $\hV^0 \oplus \hV^1 = CM_4(\hz)$, 
the Neveu-Schwarz representation of $D_4^{(1)}$, with respect to $G_2^{(1)}$, that is, with $t = u^3$,
\bea 
\hV^0 &=& \left( \sum_{m\geq 0} A_{2m} t^{2m} \right) W(\Omega_0) + 
\left( \sum_{m\geq 0} B_{2m+1} t^{2m+1} \right) W(\Omega_2) \label{NSDecomp0} \\
\hV^1 &=& \left( \sum_{m\geq 0} A_{2m+1} t^{2m+1} \right) W(\Omega_0) + 
\left( \sum_{m\geq 0} B_{2m} t^{2m} \right) W(\Omega_2) \label{NSDecomp1}  .
\eea
Because of the triality symmetry among $\hV^1$, $\hV^2$ and $\hV^3$, we must also have 
\bea 
\hV^2 &=& \left( \sum_{m\geq 0} A_{2m+1} t^{2m+1} \right) W(\Omega_0) + 
\left( \sum_{m\geq 0} B_{2m} t^{2m} \right) W(\Omega_2) \label{RDecomp2} \\
\hV^3 &=& \left( \sum_{m\geq 0} A_{2m+1} t^{2m+1} \right) W(\Omega_0) + 
\left( \sum_{m\geq 0} B_{2m} t^{2m} \right) W(\Omega_2)\label{RDecomp3} .
\eea

Using the principal graded dimensions of these modules as computed above, the left hand sides of all four of the last expressions 
become
\be \label{Allprgrdims}
\frac{\phi(u^2) \phi(u^6)}{\phi(u) \phi(u^3)}
\ee
and on the right hand sides the principal graded dimensions of the two $G_2^{(1)}$ modules are given in (\ref{G2modules}) by
$\F(u) a(u^3)$ and $\F(u) b(u^3)$, where $\F(u)$ is given in (\ref{G2Fock}). Dividing by $\F(u)$, we find the new left hand sides
are all 
\be \label{AllLHSs}
\frac{\phi(u^6)^2}{\phi(u^3)^2} = \frac{\phi(t^2)^2}{\phi(t)^2} 
\ee
giving 
\bea \label{AllSimplified}
\frac{\phi(t^2)^2}{\phi(t)^2} &=& \left( \sum_{m\geq 0} A_{2m} t^{2m} \right) a(t) + 
\left( \sum_{m\geq 0} B_{2m+1} t^{2m+1} \right) b(t)  \\
\frac{\phi(t^2)^2}{\phi(t)^2} &=& \left( \sum_{m\geq 0} A_{2m+1} t^{2m+1} \right) a(t) + 
\left( \sum_{m\geq 0} B_{2m} t^{2m} \right) b(t)   .
\eea
Adding these two equations gives
\be \label{sumeq}
2 \frac{\phi(t^2)^2}{\phi(t)^2} = A(t) a(t) + B(t) b(t)
\ee
but subtracting them gives 
\be \label{diffeq}
0 = A(-t) a(t) - B(-t) b(t)
\ee
which, after replacing $t$ by $-t$, is equivalent to 
\be \label{betterdiffeq}
0 = A(t) a(-t) - B(t) b(-t).
\ee
Now express (\ref{sumeq}) and (\ref{betterdiffeq}) as a single matrix equation, using the notation 
\be \label{cdef}
c(t) = 2 \frac{\phi(t^2)^2}{\phi(t)^2} = 2 \frac{1}{v(t)^2}
\ee
to get
\be \label{matrixeq}
\bbm a(t) &b(t) \\ a(-t) &-b(-t)\ebm \ \bbm A(t) \\ B(t) \ebm = \bbm c(t) \\ 0 \ebm .
\ee
Let 
\be \label{Dmatrix}
\dee(t) = \bbm a(t) &b(t) \\ a(-t) &-b(-t)\ebm
\ee
and 
\be \label{Ddet}
D(t) = det (\dee(t)) = - a(t) b(-t) - a(-t) b(t) .
\ee
Since $a(t)$ and $b(t)$ are invertible power series with constant term $1$, so is $a(t) b(-t)$. But then $D(t)$ is a power 
series with constant term $-2$ and $D(-t) = D(t)$, so it contains only even powers of $t$. Calculations made with the computer
program ``Mathematica" gave strong evidence for the following result. 
\begin{lem} We have $D(t) = - c(t^2)$. \index{Ramanujan identities}
\end{lem}
\begin{prf} This identity is actually one of the 40 identities in Ramanujan's famous notebook, proven by Watson \cite{Watson} 
in 1933. A modern treatment can be found in the Memoirs paper by Berndt and coauthors \cite{Berndt}. 
It is interesting to learn that this famous identity can now be given a Lie algebra representation meaning. We are very 
grateful to George Andrews for drawing our attention to the source of this identity, and to Bruce Berndt for helpful comments.\qed
\end{prf}
Now knowing this Lemma is true, we see that the matrix $\dee(t)$ is invertible, 
\be \label{Dmatrixinverse}
\dee(t)^{-1} = D(t)^{-1} \bbm -b(-t) &-b(t) \\ -a(-t) &a(t)\ebm
\ee
and we get
\be \label{AandBsolution}
\bbm A(t) \\ B(t) \ebm = \dee(t)^{-1} \bbm c(t) \\ 0 \ebm = D(t)^{-1} \bbm -b(-t) &-b(t) \\ -a(-t) &a(t)\ebm \bbm c(t) \\ 0 \ebm
= \frac{c(t)}{c(t^2)} \ \bbm b(-t) \\ a(-t) \ebm
\ee
\begin{thm} We have
\be \label{AandBformulas}
A(t) = \frac{b(-t) c(t)}{c(t^2)} = b(-t) v(-t)^2 \quad\hbox{and}\quad B(t) = \frac{a(-t) c(t)}{c(t^2)} = a(-t) v(-t)^2 
\ee
with
\be
 \frac{c(t)}{c(t^2)} = \prod_{n\geq 1} (1 + t^{2n-1})^2 = v(-t)^2 . 
\ee
\end{thm}
It now remains to show that this information allows us to write the coefficient summations in 
(\ref{NSDecomp0}), (\ref{NSDecomp1}), (\ref{RDecomp2}) and (\ref{RDecomp3}) as appropriate products of Virasoro characters. 

To verify the claimed decompositions of $\hV^0$ and $\hV^1$ as sums of tensor products of certain modules for the two
coset Virasoro algebras and the two $G_2^{(1)}$-modules, we need to check the following identities:
\bea 
\one(A(t) + A(-t)) &=& \chi_{3,4}^{1,1}(t^2) \chi_{4,5}^{1,1}(t^2) + \chi_{3,4}^{1,3}(t^2) \chi_{4,5}^{1,4}(t^2) \label{V0Omega0} \\
\one(A(t) - A(-t)) &=& \chi_{3,4}^{1,1}(t^2) \chi_{4,5}^{1,4}(t^2) + \chi_{3,4}^{1,3}(t^2) \chi_{4,5}^{1,1}(t^2) \label{V1Omega0} \\
\one(B(t) - B(-t)) &=& \chi_{3,4}^{1,1}(t^2) \chi_{4,5}^{1,3}(t^2) + \chi_{3,4}^{1,3}(t^2) \chi_{4,5}^{1,2}(t^2) \label{V0Omega2} \\
\one(B(t) + B(-t)) &=& \chi_{3,4}^{1,1}(t^2) \chi_{4,5}^{1,2}(t^2) + \chi_{3,4}^{1,3}(t^2) \chi_{4,5}^{1,3}(t^2) \label{V1Omega2} .
\eea
To verify the claimed decompositions of $\hV^2$ and $\hV^3$, we need to check the following identities:
\bea
\one(A(t) - A(-t)) &=& t\ \chi_{3,4}^{1,2}(t^2) \chi_{4,5}^{2,1}(t^2)  \label{V2or3Omega0} \\
\one(B(t) + B(-t)) &=&  \chi_{3,4}^{1,2}(t^2) \chi_{4,5}^{2,2}(t^2)  \label{V2or3Omega2}
\eea
where we have suppressed the fractional powers of $t$ which are usually associated to these characters because we really
mean to write the graded dimensions of these Virasoro modules, but with certain shifts preserved. 

\begin{thm} Using the notations 
\begin{align}
v_\pm(t) &=  \one \big(v(-t)\pm v(t)\big) \\
[v(t)b(t)]_{\pm} &= \one \big(v(-t)b(-t)\pm v(t)b(t)\big) \\
[v(t)a(t)]_{\pm} &= \one \big(v(-t)a(-t)\pm v(t)a(t)\big)
\end{align}
we have the following power series identities:
\bea 
\one(A(t) + A(-t)) &=& v_+(t) [v(t)b(t)]_+ + v_-(t) [v(t)b(t)]_-  \\ 
\one(A(t) - A(-t)) &=& v_+(t) [v(t)b(t)]_-  +  v_-(t) [v(t)b(t)]_+  \\ 
\one(B(t) - B(-t)) &=& v_+(t) [v(t)a(t)]_-  +  v_-(t) [v(t)a(t)]_+ \\ 
\one(B(t) + B(-t)) &=&v_+(t) [v(t)a(t)]_+  + v_-(t) [v(t)a(t)]_- 
\eea
\be 
\one(A(t) - A(-t)) = t \chi_{3,4}^{1,2}(t^2)  \chi_{4,5}^{2,1}(t^2) 
= t  \frac{\phi(t^4)}{\phi(t^2)} \frac{\phi(t^4)a(t^4)}{\phi(t^2)} = t \one c(t^2) a(t^4) ,
\ee
\be
\one(B(t) + B(-t)) = \chi_{3,4}^{1,2}(t^2) \chi_{4,5}^{2,2}(t^2) 
=  \frac{\phi(t^4)}{\phi(t^2)}  \frac{\phi(t^4)b(t^4)}{\phi(t^2)} = \one c(t^2) b(t^4)
\ee
\end{thm}

\begin{prf}
The six power series identities correspond to \eqref{V0Omega0} - \eqref{V2or3Omega2}. 
Computations using Mathematica verify all six identities up to a high power of $t$. Algebraic proofs for 
all six identities are given below. 

The first two identities are certainly equivalent to their corresponding sum and difference. But the sum is just
$A(t) = v(-t)^2 b(-t)$, which was established in (\ref{AandBformulas}), and the difference is just 
$A(-t) = v(t)^2 b(t)$, which is the same identity with $-t$ in place of $t$. 
The second two identities also are equivalent to their sum and difference, but we get for the sum
$B(t) = v(-t)^2 a(-t)$ which was established in (\ref{AandBformulas}), and the difference is just 
$B(-t) = v(t)^2 a(t)$, which is the same identity with $-t$ in place of $t$.

The last two identities are equivalent to 
\bea
b(-t) v(-t)^2 - b(t) v(t)^2 &=& t c(t^2) a(t^4) \qquad\hbox{and} \\
a(-t) v(-t)^2 + a(t) v(t^2) &=& c(t^2) b(t^4) 
\eea
which can be written as one matrix identity
\be
\bmatrix b(-t) &b(t) \\a(-t) &-a(t)\endbmatrix \bmatrix v(-t)^2\\-v(t^2)\endbmatrix = c(t^2) \bmatrix ta(t^4) \\ b(t^4)\endbmatrix.
\ee
Dividing both sides by $-D(t) = c(t^2)$ gives 
\be
\dee(t)^{-1} \bmatrix v(-t)^2\\-v(t^2)\endbmatrix = 
D(t)^{-1} \bmatrix -b(-t) &-b(t) \\-a(-t) &a(t)\endbmatrix \bmatrix v(-t)^2\\-v(t^2)\endbmatrix = \bmatrix ta(t^4) \\ b(t^4)\endbmatrix
\ee
which is equivalent to 
\be
\dee(t)\bmatrix ta(t^4) \\ b(t^4)\endbmatrix =  \bmatrix v(-t)^2\\-v(t^2)\endbmatrix.
\ee
The two equations this gives, 
\bea
t a(t) a(t^4) + b(t) b(t^4) &=& v(-t)^2 \qquad\hbox{and} \\
t a(-t) a(t^4) - b(-t) b(t^4) &=& -v(t)^2, 
\eea
are equivalent to each other under the substitution of $-t$ for $t$. 
The first one, where $a(t) = H(t)$ and $b(t) = G(t)$, is one of the forty Ramanujan identities, 
\be
G(t) G(t^4) + t H(t) H(t^4) = \frac{\phi(t^2)^4}{\phi(t^4)^2 \phi(t)^2} \index{Ramanujan identities}
\ee
proven by Watson \cite{Watson} in 1933. (Also see \cite{Berndt}.)\qed
\end{prf}

\chapter{Detailed Constructions} \label{cha:constructs}

This chapter embellishes in greater detail the spinor construction\index{spinor construction} of the finite dimensional Lie algebras of type $D_4, B_3$, and $G_2$ as well as the untwisted affine Kac-Moody Lie algebra  of type $D_4^{(1)}$. We also provide the details for constructing the conformal vectors $\D, \B, \G, \DB$ and $\BG$ using the Sugawara construction. 

\section{Spinor Construction of $D_4$ and $D_4$-Modules }\label{sec:finspinconstruct}

Let $A = A^+\oplus A^-$ where $A^\pm\cong\C^{4}$ are maximal isotropic subspaces with respect to a non-degenerate symmetric bilinear form on $A$. Choose canonical bases 
$\mathcal A^+ = \{a_1,\cdots,a_4\}$ and $\mathcal A^- = \{a^*_1,\cdots,a^*_4\}$ for $A^+$ and $A^-$, respectively, such that 
\be\label{cliffpairing}  (a_i,a^*_j)=\di  \quad \hbox{and} \quad  (a_i,a_j)=(a^*_i,a^*_j)=0   \quad \hbox{for all} \quad 1\leq i,j \leq 4,\ee
and let $\mathcal A = \mathcal A^+ \cup\mathcal A^-$. 
Define the \emph{tensor spaces}\index{tensor!space} $A^{\otimes n}$ for $n\geq 0$ in the following way:
$$
A^{\otimes 0}:=\C,\qquad
A^{\otimes 1}:=A \simeq \C \otimes A,\qquad
 \cdots  \qquad
A^{\otimes n}:=A^{\otimes (n-1)}\otimes A \textrm{ for } n>1
$$
Denote the (associative) \index{tensor!algebra}\emph {tensor algebra} $\ds\mathscr T=\mathscr T(A):= \bigoplus_{n\geq 0} A^{\otimes n} $ as the infinite direct sum of the tensor spaces $A^{\otimes n}$. Let $\mathscr J$ be the ideal in $\mathscr T$ generated by all of the elements of the form $v_1\otimes v_2+v_2\otimes v_1-(v_1,v_2)1$, where $1\in A^{\otimes 0}=\C$.
Define the (associative) Clifford algebra\index{Clifford algebra} as the quotient 
$$\Cl = \Cl(A, (\cdot,\cdot)) = \mathscr T/\mathscr J .$$
When the context is clear, we will write elements $v+\mathscr J$ of Cliff as $v$, and we will write products  $(v+\mathscr J)(w+\mathscr J)$ as $vw$. The unit $1+\mathscr J$ will be written as 1, when there is no chance for ambiguity.

Quite often we will apply the \emph {main Clifford relation}  \index{Clifford algebra!main relation} 
$v_1v_2 + v_2v_1 = (v_1,v_2)1$ to the basis elements in $\mathcal A$. In particular we have  
\be\label{maincliff} a_ia^*_j+a_j^*a_i= \di 1, \qquad
a_ia_j=-a_ja_i \qquad \textrm{ and }  \qquad a_i^*a_j^*=-a_j^*a_i^*  \ee
which implies
\be\label{cliffsquare}(a_i)^2=(a^*_i)^2=0.\ee
Let $\I$ be the left ideal of $\Cl$ generated by $A^+$, so that as in Definition \ref{def:CliffordAlg}, define the left $\Cl$-module $\CM = \Cl/\I = \Cl\cdot\bv = (\wedge A^-)\cdot \bv$ \index{module!Clifford} where $\bv = 1+\I$ is called the vacuum vector\index{vacuum vector}. It follows that $A^+\cdot\bv = 0$ and that 
\be\label{cliffmodule}\CM = \mbox{ span} \left\{\left(a^*_1\right)^{k_1} \left(a^*_2\right)^{k_2} \left(a^*_3\right)^{k_3} \left(a^*_4\right)^{k_4}\bv \left|\right. k_i \in\{0,1\}\right\}.\ee
 Define its ``even" and ``odd" \index{module!Clifford!even, odd} subspaces $\CM^j$ for $j\in\{0,1\}$, as the span above with the additional condition that: $$\sum_{i=1}^4k_i\equiv j\ (mod2).$$

 For $j\in\{0,1\}$, $dim(\CM^j)=8$ and   $\CM = \CM^{0} \oplus \CM^{1}$ as vector spaces where
 $\CM^0$ is spanned by only  even numbers of elements of $\mathcal A^-$ applied on the left of $\bv$ and
 $\CM^1$  is spanned by only odd  numbers of elements of $\mathcal A^-$ applied to $\bv$.

There is an anti-automorphism defined on $\Cl$, $\alpha: \Cl \rightarrow\Cl$, such that for $v_1,...,v_n \in A$, 
$$\alpha(v_1v_2\cdots v_{n-1}v_n) \quad = \quad v_nv_{n-1}\cdots v_2v_1.$$
It is easy to show that  $\alpha(\bv)= \bv$ and also that for any $u\in$  Cliff , $(\bv)u(\bv)$ is a scalar multiple of $\bv$. Define a non-degenerate, symmetric bilinear form $(\cdot,\cdot)$ on $\CM$ such that for $b,c\in \CM$,  $(b,c)$ is the unique scalar for which $\alpha(b)c=(b,c)\bv.$

\newpage

\begin{prop}\label{prop:pair}\cite{FFR} For $a\in A$ and $ b,c\in \CM$ we have 
\begin{enumerate}[\rm(a)]
\item $(ab,c) = (b,ac).$
\item $(ab,ac) = \one(a,a)(b,c).$
\item $(\CM^0,\CM^1)=0.$
\end{enumerate}
\end{prop}
Denote the \index{Chevalley algebra} \emph{Chevalley algebra}  by  $\mathcal C= A \oplus \CM^0 \oplus \CM^1$ and extend the forms on $A$ and $\CM$ to a non-degenerate symmetric bilinear form on $\mathcal C$ where $(A,\CM)=0$, making the direct sum for $\mathcal C$ above an  orthogonal direct sum of vector spaces.
For $u\in \mathcal C$ write $ u=a+b+c$ where
$a\in A, b\in \CM^0, c\in \CM^1$.
Define a cubic form $\phi:\mathcal C\rightarrow \C$ by $\phi(u)= (ab,c)$ and a trilinear form $\Phi: \mathcal C\times \mathcal C\times \mathcal C\rightarrow \C$ \quad  obtained by polarization of $\phi$,
$$\Phi(u_1,u_2,u_3) = \phi(u_1+u_2+u_3) - \phi(u_2+u_3) - \phi(u_1+u_3) - \phi(u_1+u_2) + \phi(u_1) + \phi(u_2) + \phi(u_3).$$
Define the \index{$\circ$ product} operation\ $\circ$\ :$\mathcal C \times \mathcal C\rightarrow \mathcal C$ so that for $u_1,u_2\in \mathcal C$, $u_1\circ u_2$ is the unique element of $\mathcal C$ such that for all $u\in \mathcal C$ , \be\label{circprod} \Phi(u_1,u_2,u) = (u_1\circ u_2, u). \ee

\begin{prop}\label{prop:circprod}\cite{FFR} The $\circ$ product has the following properties:
\begin{enumerate}[\rm (a)]
\item The operation $\quad\circ\quad$ is commutative.
\item $A\circ A =\CM^0\circ \CM^0=\CM^1\circ \CM^1=0$.
\item $a\circ b=ab$\quad for \quad  $a\in A, \quad b\in \CM^0\oplus \CM^1$.
\item $\CM^0\circ \CM^1\subseteq A$.
\item $(u_1\circ u_2, u_3 ) = (u_1, u_2 \circ u_3)$ for all  $u_1, u_2, u_3 \in \cee.$
\end{enumerate}
\end{prop}

 The following will play a large role throughout. We will only present the details necessary for this work, however for a more comprehensive discussion one may consult \cite{FFR} or \cite{Chevalley}. Define 
 \bea\label{e0} e_1 := \ph &\ph a_4+a_4^* \ph\in\ph A \quad \quad \quad \mbox{ and } \notag\\  e_2:= \ph &\ph a_1^*a_4^*\bv-a_2^*a_3^*\bv \ph \in \ph \CM^0 \nonumber \\  e_3:= \ph &\ph -a_1^*\bv -a_2^*a_3^*a_4^*\bv  \ph \in \ph \CM^1.\eea
 It is easy to see that $e_3=e_1e_2 $,  and also that  \be  (e_1,e_1)=2=(e_2,e_2)=(e_3,e_3) .  \label{e1}\ee Further, define the operators \begin{align*} \rho(e_1): &\mathcal C\rightarrow \mathcal C \quad \mbox{ and} \\ \rho(e_2): &\mathcal C\rightarrow \mathcal C\quad \mbox{ by } \end{align*}
\bea\label{rhoaction}
&\rho(e_1)b =& e_1\circ b \quad \hbox{for} \quad b\in \CM^0\oplus \CM^1 \notag\\
&\rho(e_1)a =& (a,e_1)e_1-a \quad \hbox{for} \quad a\in A. \notag\\
&\rho(e_2)c =& e_2\circ c \quad \hbox{for} \quad c\in A\oplus \CM^1 \notag\\
&\rho(e_2)b =& (b,e_2)e_2-b \quad \hbox{for} \quad b\in \CM^0.
\eea

Denote the following  automorphisms:  \bea \label{rho}\tau  := \ph & \ph \rho(e_1)\nonumber \\ \sigma \ph :=  \ph &\rho(e_1)\rho(e_2).\eea
Then $\sigma^3=\tau^2=1$, \quad and $\quad \tau\sigma\tau=\sigma^{-1}. $
The operators \quad $\sigma$ \quad and \quad $\tau$ \quad are automorphisms of $(\mathcal C,\circ)$, and they  generate the \emph{triality group}\index{triality}, which is isomorphic to the symmetric group $S_3$. One may use $\sigma$ to polarize the $\CM^j$ into maximally isotropic subspaces $\CM^{j\pm}$ which are the images of $A^{\pm}$ under $\sigma$ and $\sigma^2$, that is, $\sigma(A^{\pm})=$ \CM$^{0\pm}$ and  $\sigma^2(A^{\pm})=$ \CM$^{1\pm}$. As vector spaces $A\cong \CM^0\cong \CM^1$ so we can consider the Clifford algebras:
$\mbox{Cliff}(\CM^0, (\cdot,\cdot))$ and $\mbox{Cliff}(\CM^1, (\cdot,\cdot)).$ We will suppress the bilinear forms when  there is no chance for ambiguity.
Using $\sigma$  it is clear that $$\mbox{ Cliff}\cong \mbox{ Cliff}(\CM^0) \cong \mbox{ Cliff}(\CM^1).$$ 

For any $u,v\in A$ we define the \emph{fermionic normal ordering} \index{normal order!fermionic} $\fno uv \fno$ by \be\label{fno} \fno uv\fno \ph := \ph \tfrac{1}{2}(uv-vu).\ee Immediately we see that $\fno uv\fno = -\fno vu\fno = uv-\tfrac{1}{2}(u,v)1.$ Consider the subspace of Cliff given by $\g := \mbox{ span} \{\fno uv\fno \st u, v \in A \}$. Note that for $\fno u_1v_1\fno, \fno u_2v_2\fno \in \g$ we have the commutator
 \be \left[\fno u_1v_1\fno, \fno u_2v_2\fno\right] = (v _1,u_2)\fno u_1v_2\fno + (u_1,v_2)\fno v_1u_2\fno - (u_1,u_2)\fno v_1v_2\fno - (v_1,v_2)\fno u_1u_2\fno\label{adjoint}\ee
 which can be computed from the definition. This shows that $\g$ is a Lie algebra closed under this bracket. There is a
 canonical action of $\g$ on $(\mathcal C, \circ)$ , given by \bea
\fno uv\fno\cdot a & = (v,a)u-(u,a)v \mbox{ for } a\in A\label{natural}\\
\fno uv\fno\cdot b & = \tfrac{1}{2}(u\circ(v\circ b))-v\circ(u\circ b)) \mbox{ for } b \in \CM^0\label{spinoreven} \\
\fno uv\fno\cdot c &=\tfrac{1}{2}(u\circ(v\circ c))-v\circ(u\circ c)) \mbox{ for } c \in \CM^1 \label{spinorodd} \eea which preserves each of the vector spaces $A$, $\CM^0$, and  $\CM^1$. Hence $\g$ is a Lie algebra of operators on $(\mathcal C, \circ)$. Similarly, we may repeat \eqref{fno} - \eqref{spinorodd} with $(A, \CM^0, \CM^1)$ replaced by $(\CM^0, \CM^1, A)$ or by $(\CM^1, A, \CM^0)$, yielding two additional Lie algebras $\sigma(\g)\subseteq \mbox{ Cliff}(\CM^0)$ and $\sigma^2(\g)\subseteq \mbox{ Cliff}(\CM^1)$, both acting on $(\mathcal C, \circ)$. This involves identifying $\CM^1\oplus A$ as the left Clifford module for Cliff$(\CM^0)$ and $A\oplus \CM^0$ as the left Clifford module for Cliff$(\CM^1)$. We use the notation \bea
\sigma\fno uv\fno & := \fno(\sigma u)(\sigma v)\fno \in \mbox{ Cliff}(\CM^0) &\hbox{ and } \notag\\
\sigma^2\fno uv\fno &  := \fno(\sigma^2 u)(\sigma^2 v)\fno  \in \hbox{ Cliff}(\CM^1) \notag
\eea for $u,v\in A$ to denote the elements in $\sigma(\g)$ and $\sigma^2(\g)$ respectively. As shown in \cite{FFR} these three Lie algebras of operators on $(\mathcal C, \circ)$ are identical, so that $\sigma$ acts on one Lie algebra $\g$ as an automorphism of order three. In section \ref{sub:circle} we will further discuss the action of $\g$ on  $(\mathcal C, \circ)$ and also give the identification $\g=\sigma(\g)=\sigma^2(\g)$ implied by \cite{FFR} via the ordered basis $\mathcal A$ of $A$.
To simplify future calculations we will occasionally adopt the following notation
\bea\label{finmods}
V^{(1)}:= A \phantom{M^0} & \quad V^{(2)}:= \CM^0 & \qquad  V^{(3)}:= \CM^1 \nonumber\\
\g^{(0)} :=  \g\phantom{\sigma^2)} & \quad \g^{(1)} := \sigma(\g)\phantom{^2} & \qquad \g^{(2)} := \sigma^2(\g). \label{upper}
\eea

Note that the triality automorphism $\sigma$ not only induces the Lie algebra isomorphism above, but also permutes the summands of the \index{Chevalley algebra} Chevalley algebra
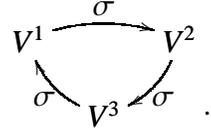
  \begin{xy}
{\ar@{->}@/^/ (0,10)*+{V^1}; (20,10)*+{V^2}}?*!/_2mm/{\sigma}; 
{\ar@{->}@/^/ (20,10)*+{V^2}; (10,0)*+{V^3}}?*!/_2mm/{\sigma}; 
{\ar@{->}@/^/ (10,0)*+{V^3}; (0,10)*+{V^1}}?*!/_2mm/{\sigma}
\end{xy}.

\begin{rmk} \label{rmk:LievsOps} Although $\g$ has been constructed above as a Lie algebra of operators, we may also use the terminology given in Section \ref{sec:finite} to discuss aspects of $\g$ such as root spaces, Cartan subalgebra, representations and their weights. 
\end{rmk}

Denote by $E_{i,j}$ the $8\times8$ complex matrix with every entry equal to zero except for the $(i,j)^{th}$ entry, which shall be 1. The standard description of the Lie algebra $\mathfrak{so}(8)$, as can be found in \cite{Humphreys}, is the algebra of $8\times8$ complex matrices of the block form $\left[\begin{array}{cc}
m&p\\q&n\\ \end{array}\right]$, where $\ph p=-p^T,\ph q=-q^T, \ph n=-m^T$; with basis
\begin{align*}
&\left\{E_{i,j+4}-E_{j, i+4} \st 1\leq i<j\leq 4\right\}\ph\cup\\
&\left\{E_{i+4,j}-E_{j+4, i} \st 1\leq i<j\leq 4\right\}\ph\cup\\
&\left\{E_{i,j}-E_{j+4, i+4} \st 1\leq i,j \leq 4\right\}.
\end{align*}
The diagonal matrices forming its Cartan subalgebra (CSA) are the matrices in the last set with $i=j$.
There is a  map from $\g$ to  $\mathfrak{so}(8)$ given by
\begin{align}
\fno a_ia_j\fno  &\longleftrightarrow E_{i,j+4}-E_{j,i+4} \notag \\
\fno a_i^*a_j^*\fno  &\longleftrightarrow E_{i+4,j}-E_{j+4,i} \notag \\
\fno a_ia_j^*\fno  &\longleftrightarrow E_{i,j}-E_{j+4,i+4}\label{Es}.
\end{align} Using \eqref{adjoint}, one may verify that this map is a Lie algebra isomorphism $\g\cong \mathfrak{so}(8)$.
The action of $\g$ on itself in \eqref{adjoint} corresponds to the adjoint representation of $\mathfrak{so}(8)$. We have also seen in \eqref{natural}-\eqref{spinorodd} that the actions of $\g$ on $A, \CM^0$, and $\CM^1$ yield the natural and the two irreducible spinor representations of $\mathfrak{so}(8)$. We note here that there is a non-degenerate invariant symmetric bilinear form on $\g$ given by   
\be(\fno u_1v_1\fno,\fno u_2v_2\fno) = (u_1,v_2)(v_1,u_2)-(u_1,u_2)(v_1,v_2)\label{spinkill}.\ee

Let $h_i=\fno a_ia_i^*\fno$, $1\leq i \leq 4$, and note that $\{h_1,h_2,h_3,h_4\}$ is an orthonormal basis, with respect to \eqref{spinkill}, for the canonical CSA, $\h$, corresponding to the diagonal matrices in \eqref{Es} when $i=j$. 
 Defining the linear functionals $\ep_i\in \h^*$ by $\ep_i(h_j)=\di$  for $1\leq i \leq 4$, we see that the dual basis $\{\ep_1,\ep_2,\ep_3,\ep_4\}$ of $\h^*$ is orthonormal with respect to the induced bilinear form on $\h^*$. We will use the standard notation for the root system of $\mathfrak{so}(8)$, of type $D_4$ given as
\be\label{Dsystem} \Phi_{D_4} = \{\pm\ep_i\pm\ep_j\st 1\leq i<j\leq 4\} \ee with the simple roots given as \be\label{Dsimple} \Delta_{D_4} = \{\alpha_1 = \ep_1-\ep_2, \ph\ph \alpha_2 = \ep_2-\ep_3, \ph\ph \alpha_3 = \ep_3-\ep_4, \ph\ph \alpha_4 = \ep_3+\ep_4\}. \ee
In fact, using the bracket formula \eqref{adjoint} we see the following correspondence of root vectors in $\g$ with their roots
\bea
\fno a_ia_j\fno  &\longleftrightarrow \ep_i+\ep_j & 1\leq i<j\leq 4\nonumber\\
\fno a_i^*a_j^*\fno  &\longleftrightarrow -\ep_i-\ep_j & 1\leq i<j\leq 4\nonumber \\
\fno a_ia_j^*\fno  &\longleftrightarrow \ep_i -\ep_j & 1\leq i\neq  j\leq 4\label{roots}.
\eea
The fundamental weights of type $D_4$ are then 
\be\label{D4fundamentalweights} \lambda_1 = \ep_1,\quad \lambda_2 = \ep_1+\ep_2,\quad 
\lambda_3 = \one( \ep_1+\ep_2+ \ep_3 - \ep_4),\quad \lambda_4 = \one( \ep_1+\ep_2+ \ep_3 + \ep_4) .\ee
Using the formula \eqref{natural} we see the following correspondence of weight vectors in $V$ with their weights
\bea
&a_i \longleftrightarrow \ep_i & 1\leq i\leq 4\nonumber\\
&a_i^* \longleftrightarrow -\ep_i & 1\leq i\leq 4.\label{naturalweights}
\eea
Using the formulas \eqref{spinoreven}-\eqref{spinorodd} we see the following correspondence of weight vectors in $\CM^0$ and $\CM^1$ with their weights
\bea
\bv &\longleftrightarrow& \tfrac{1}{2}(\ep_1+\ep_2+\ep_3+\ep_4)\nonumber\\
a_i^*\bv &\longleftrightarrow& \tfrac{1}{2}(\ep_1+\ep_2+\ep_3+\ep_4)-\ep_i\nonumber\\
a_i^*a_j^*\bv &\longleftrightarrow& \tfrac{1}{2}(\ep_1+\ep_2+\ep_3+\ep_4)-\ep_i-\ep_j\nonumber\\
a_i^*a_j^*a_k^*\bv &\longleftrightarrow& \tfrac{1}{2}(\ep_1+\ep_2+\ep_3+\ep_4)-\ep_i-\ep_j-\ep_k\nonumber\\
a_1^*a_2^*a_3^*a_4^*\bv &\longleftrightarrow& -\tfrac{1}{2}(\ep_1+\ep_2+\ep_3+\ep_4)\label{spinorweights}.
\eea
Thus the  set of weights for the natural representation on $V$ is $\{\pm\ep_i\st 1\leq i \leq 4 \}$ and its highest weight is $\ep_1=\lambda_1$. 
The set of weights for the spinor representation on $\CM^0\oplus \CM^1$ is $\left\{\one(\pm\ep_1\pm\ep_2\pm\ep_3\pm\ep_4)\right\}$. Specifically, for $j\in \{0,1\}$, the weights of $\CM^j$ are those where the number of minus signs is congruent to j mod 2. The highest weight in $\CM^0$ is
$\one( \ep_1+\ep_2+ \ep_3 + \ep_4) = \lambda_4$ and the the highest weight in $\CM^1$ is $\one( \ep_1+\ep_2+ \ep_3 - \ep_4) = \lambda_3$

Recall that we defined $\sigma$ and $\tau$ as diagram automorphisms. The diagram below shows their action on the simple roots of $D_4$ as a Dynkin diagram automorphism.     
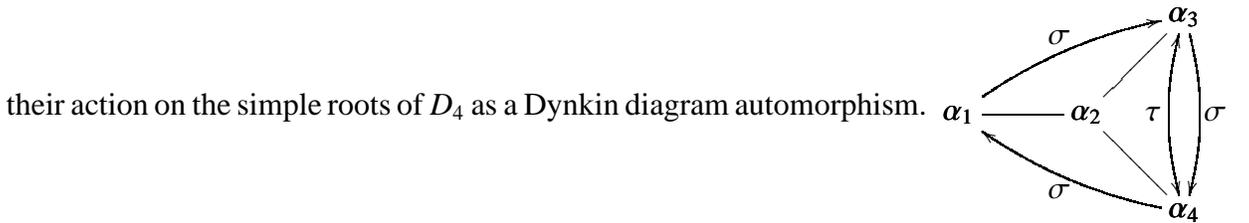
\begin{xy}
{\ar@{-} (0,0)*+{\alpha_1}; (17,0)*+{\alpha_2}}; 
{\ar@{-}(17,0)*+{\alpha_2}; (30,13)*+{\alpha_3}}; 
{\ar@{-}(17,0)*+{\alpha_2}; (30,-13)*+{\alpha_4}}; 
{\ar@{->}@/^/ (0,0)*+{\alpha_1}; (30,13)*+{\alpha_3}}?*!/_2mm/{\sigma}; 
{\ar@{->}@/^/ (30,13)*+{\alpha_3}; (30,-13)*+{\alpha_4}}?*!/_2mm/{\sigma}; 
{\ar@{->}@/^/ (30,-13)*+{\alpha_4}; (0,0)*+{\alpha_1}}?*!/_2mm/{\sigma}; 
{\ar@{<->}@/^/ (30,-13)*+{\alpha_4}; (30,13)*+{\alpha_3}}?*!/_2mm/{\tau}; 
\end{xy}

\section{A Spinor Description of $G_2$ as a Subalgebra of $D_4$} \label{sec:G2}

In this section we use the spinor description of $D_4$ and the automorphism $\sigma$ to give a spinor description of  $G_2$. 

\subsection*{Invariance of the Cartan Subalgebra}\label{sub:CartanG2}

In section \ref{sec:finspinconstruct} we choose a CSA
$\h\subset \g$ with ordered orthonormal basis $\{h_1,...,h_4 \}$. With analogous notation to \eqref{upper}, there are CSA's $\h^{(1)}:= \sigma(\h)\subset \g^{(1)}$ and
 $\h^{(2)}:= \sigma^2(\h)\subset \g^{(2)}$. In fact, we have the following:

\begin{thm}
As algebras of operators on $(\mathcal C,\circ)$,  $\h^{(0)} :=\h = \h^{(1)} = \h^{(2)}$.
\end{thm}

\begin{prf}
Since $\g\cong \mathfrak{so}(8)$ is a simple Lie algebra, the irreducible representation $\phi_{V^{(0)}}: \g \rightarrow \e(V)$ of $\g$ on $V$ is non-trivial and therefore faithful. The representations $\phi_{V^{(1)}}: \g \rightarrow \e(\CM^0)$ and $\phi_{V^{(2)}}: \g \rightarrow \e(\CM^1)$ are also faithful and irreducible. For all $x\in\g$, these three representations satisfy the relations
\begin{align}
\sigma\inv\circ\phi_{V^{(i)}}(x)\circ\sigma=\phi_{V^{(i-1)}}(\sigma\inv x) &\ph &\quad \mbox{ where}\quad i,i-1 \in \{0,1,2\}\mod3.\label{repi}
\end{align}
This is summarized by the commutative diagram below.
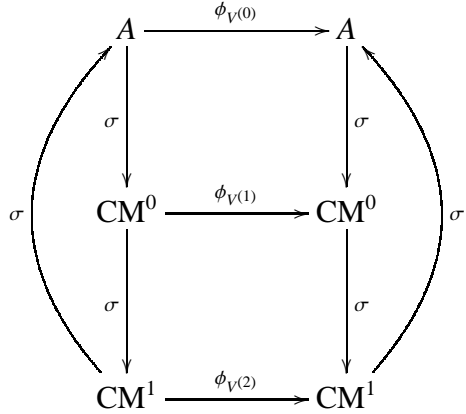
\[\xymatrix{  A \ar[rr]^{\phi_{V^{(0)}}} \ar[dd]_{\sigma}  &  & A \ar[dd]^{\sigma}\\  & & \\
 \CM^0 \ar[rr]^{\phi_{V^{(1)}}}\ar[dd]_{\sigma}  &  & \CM^0
\ar[dd]^{\sigma}  \\
& & \\
\CM^1 \ar[rr]^{\phi_{V^{(2)}}}\ar@/^3pc/[uuuu]^{\sigma} &  & \CM^1 \ar@/_3pc/[uuuu]_{\sigma}\\
    }\]

Consider $\sigma(h_1) \in \sigma(\g)$. It can be shown that $\sigma(h_1) - \tfrac{1}{2}(h_1+h_2+h_3+h_4)$ acts as the zero of $\g$ on $V$, and hence by the above, on all of $(\mathcal C,\circ)$. Thus $\sigma(h_1) = \tfrac{1}{2}(h_1+h_2+h_3+h_4)$ as an operator on $(\mathcal C,\circ)$; $\sigma(h_1)\in \h$. Similar calculations can be done for $h_2,h_3,h_4$ with the result that $\sigma(\h) \subseteq \h$. Applying sigma proves the theorem. $ \h \subseteq \sigma(\h) \subseteq \sigma^2(\h) \subseteq \h$. 
\qed\end{prf}
The remaining identifications of the operators in the CSA (and actually in all of $\g$) are provided in Lemma \ref{lem:identification}.

\subsection*{Roots of $G_2$}\label{sub:g2roots}
Recall the description of the root system for the finite dimensional rank 4 Lie algebra of type $D_4$ given in \eqref{Dsystem};
$\Phi_{D_4} = \{\pm\ep_i\pm\ep_j\st 1\leq i<j\leq 4\}$, the simple roots \eqref{Dsimple}
$\Delta_{D_4} = \{\alpha_i \st 1\leq i \leq 4\}$, and the correspondence of root vectors with the roots given in \eqref{roots}. Restricting the automorphism $\sigma$ to $\h$ gives the automorphism of the previous theorem $\sigma: \h\rightarrow \h$ and induces a dual map $\sigma^*:\h^*\rightarrow\h^*$ such that
\be\label{sigstar}
\sigma^*(\alpha_1)=\alpha_4 ,\qquad
\sigma^*(\alpha_2)=\alpha_2 ,\qquad
\sigma^*(\alpha_3)=\alpha_1 ,\qquad
\sigma^*(\alpha_4)=\alpha_3.
\ee
Actually, since $\sigma^*$ permutes the simple roots of $\g$ it is also the case that $\sigma$ permutes the root spaces of $\g$. That is, for $\mu \in \Phi_{D_4}$ \be \sigma: \g_{\mu}\to \g_{\sigma^*(\mu)}.\label{rootperm}\ee
One may also compute the action of $\sigma^*$ on the individual $\ep_i$ via the appropriate use of \eqref{sigstar}.

The the fixed points of $\g$ under $\sigma$ form a rank 2 finite dimensional simple exceptional Lie algebra of type $G_2$, with CSA  denoted by $\mathfrak k$. We can restrict both $\sigma$ and $\sigma^*$,  to $\mathfrak k$ and its dual space $\mathfrak k^*$, to find the following root system of the  $G_2$
\be\label{Gsystem} \Phi_{G_2} = \pm\{ \beta_2,\ph\ph \beta_1,\ph\ph \beta_2+\beta_1,\ph\ph 2\beta_2+\beta_1,\ph\ph 3\beta_2+\beta_1,\ph\ph 3\beta_2+2\beta_1 \} \ee with the short and long simple roots given respectively as \be\label{Gsimple} \Delta_{G_2} = \{\beta_2 =\tfrac{1}{3}( \ep_1-\ep_2+2\ep_3), \ph\ph \beta_1 = \ep_2-\ep_3\} \ee
which has the following Dynkin diagram 
\quad
\begin{tikzpicture}
\draw(5,0)--(6,0);
\draw (5,-.1)--(6,-.1);
\draw (5,.1)--(6,.1);
\shade[ball color =white](5,0)circle(6pt);
\shade[ball color =white](6,0)circle(6pt);
\draw(5.4,0)--(5.6,.2);
\draw(5.4,0)--(5.6,-.2);
\draw(5,.5)node[black]{$\beta_2$};
\draw(6,.5)node[black]{$\beta_1$};
 \end{tikzpicture}.
 
The fundamental weights of $G_2$ are $\bar\lambda_1 = 2\beta_1+3\beta_2 = \ep_1+\ep_2$ and 
$\bar\lambda_2 = \beta_1+2\beta_2 = \tfrac{1}{3}(2\ep_1+\ep_2+\ep_3)$. The roots and weights of $G_2$ are elements of the $2$-dimensional dual Cartan subalgebra $\mathfrak k^*$, which is inside the $4$-dimensional dual Cartan subalgebra $\mathfrak h^*$ of $D_4$. It will be useful later to know the formula for the projection from $\mathfrak h^*$ to $\mathfrak k^*$, which can be written as either 
\be\label{D4toG2projectionalphas} Proj_{\mathfrak k^*}\Big( \sum_{i=1}^4 a_i \alpha_i\Big) = a_2\beta_1 + (a_1+a_3+a_4)\beta_2 \ee
or as
\be\label{D4toG2projectionepsilons} Proj_{\mathfrak k^*}\Big( \sum_{i=1}^4 e_i \ep_i\Big) 
= (e_2-e_3)\bar\lambda_1 + (e_1-e_2+2e_3)\bar\lambda_2 . \ee
These tell us that $\alpha_1$, $\alpha_3$ and $\alpha_4$ project to $\beta_2$ but $\alpha_2$ projects to $\beta_1$, and that
$\lambda_1$, $\lambda_3$ and $\lambda_4$ project to $\bar\lambda_2$ but $\lambda_2$ projects to $\bar\lambda_1$. 
The $14$-dimensional adjoint representation of $G_2$ is the irreducible $G_2$-module with highest weight $\bar\lambda_1$,
and the $7$-dimensional representation of $G_2$ is the irreducible $G_2$-module with highest weight $\bar\lambda_2$. 

\subsection*{Decompositions of $\g$ and $\cee$ as $G_2$-Modules}\label{sub:irrep}

The spinor construction described in \ref{sec:finspinconstruct} yields three isomorphic copies of $D_4$.
We may use the notation $\phi_{\mathcal C}: \g \rightarrow \e(\mathcal C)$ for the reducible representation of $\g$ on $\cee$ so that
\begin{align}\label{sigrep} \sigma\inv \circ\phi_{\cee}(x)\circ\sigma = \phi_{\cee}(\sigma\inv(x))\end{align} as in \eqref{repi}. We identify $\phi_{\cee}(\g)=\g$.
The action of $\sigma$ on $\g$ decomposes it into three eigenspaces
\be\label{eigensigma}\g= \g_0\oplus \g_1\oplus\g_2\ee
where $\g_i= \left\{ x\in \g  \ph | \ph \sigma(x)=\xi^i x \right\} $
for $i\in\{1,2,3\}$, and  where $\xi=e^{\frac{2\pi i}{3}}=\frac{-1+\sqrt{3}i}{2}$ and $\xi^2=e^{\frac{-2\pi i}{3}}=\frac{-1-i\sqrt{3}}{2}$. Note that dim$(\g _0) = 14$ and  dim$(\g _1) =$ dim$(\g _2) = 7$. In fact, $\g _0$ is a simple Lie algebra of type $G_2$ and $\g _1$ and $\g _2$ are irreducible $\g_0$-modules.

Recall that for $i\in\{1,2,3\}$ each  $V^{(i)}$  is an 8-dimensional irreducible $\g $ module. However as a  $\g_0$ module, $V^{(1)}$ is reducible into the direct sum $V^{(1)}=W^{(1)}\oplus e^{(1)}$, where $e^{(1)}$ is spanned by $e_1$, defined in \eqref{e0}, and $W^{(1)}$ has ordered basis 
\be \label{Bbasis}\mathcal B^{(1)}= \{a_1,a_2,a_3, a_4-a_4^*, a_3^*, a_2^*a_1^*\}.\ee
Applying $\sigma$ and  $\sigma^2$ to this decomposition of $V^{(1)}$ gives $V^{(2)}=W^{(2)}\oplus  e^{(2)}$ and also $V^{(3)}=W^{(3)}\oplus  e^{(3)}$. That is, $\sigma$ cyclically permutes the   $W^{(i)}$ and the  $e^{(i)}$ for  $i\in\{1,2,3\}$, which allows us to write $\cee$ as the direct sum of the irreducible $\g_0$ modules
\be\label{Wmodules} \cee= \bigoplus_{i=1}^3 \left( W^{(i)}\oplus  e^{(i)} \right).\ee We agree to use the ordered bases $B^{(2)}$ and $B^{(3)}$ of $W^{(2)}$ and $W^{(3)}$ respectively, which are induced from the action of $\sigma$ on $\mathcal B^{(1)}$.

Note that $\mathcal C$ also has a $\sigma$ eigenspace decomposition as irreducible $\g _0$ modules \begin{align} \cee = \bigoplus_{i=1}^3 \left( W^{\xi(i)}\oplus  e^{\xi(i)} \right)\label{Umodules} \end{align} where the 7-dimensional  $W^{\xi(1)}$, $W^{\xi(2)}$, and $W^{\xi(3)}$ have ordered bases
\begin{align}
B^{\xi(1)} & = \left\{b+\sigma b+\sigma^2b \st b\in \mathcal B \right\} \nonumber \\
B^{\xi(2)} & = \left\{b+\xi^2\sigma b+\xi\sigma^2b \st b\in \mathcal B \right\} \nonumber \\
B^{\xi(3)} & = \left\{b+\xi\sigma b+\xi^2\sigma^2b \st b\in \mathcal B \right\} \label{Wxibases} \end{align}
respectively, and where
\begin{align}
e^{\xi(1)} & = span\left\{ e_0 + e_1 + e_2 \right\} \nonumber \\
e^{\xi(2)} & = span\left\{ e_0 + \xi^2e_1 + \xi e_2 \right\} \nonumber \\
e^{\xi(3)} & = span\left\{ e_0 + \xi e_1 + \xi^2e_2 \right\}. \label{exibases} \end{align}

\subsection*{The $\circ$ Products}\label{sub:circle}

The fixed points of $\g$ under $\sigma$ formed a Lie algebra of type $G_2$ which we have denoted as $\g_0$. Here we use the $\circ$ product to provide the details on finding the root vectors of $\g _0$. For $\cee= A\oplus \CM^0\oplus \CM^1 $, we can use \eqref{circprod} and proposition \ref{prop:circprod} in the situation when $u_1\in \CM^0,u_2\in \CM^1$ (or vice versa), so that $u_1\circ u_2$ is the unique element of $V$ such that for all $a\in A$  $(u_1\circ u_2, a)=(au_1, u_2)$. When both $u_1$ and $u_2$ are in the same direct summand, $u_1\circ u_2=0$, and when $u\in \CM^0, a\in A$, $u\circ a= au\in \CM^1$, with the analogous result when $u\in \CM^1$. We can explicitly compute the product for pairs of elements in $\CM$. For $i,j,k,l\in \{1,2,3,4\}$,  we will use the following notation whenever $1\leq i<j<k\leq4$:
\begin{align}\label{circnote}
v&:=\bv \nonumber \\
v_i&:=a_i^*\bv \nonumber\\
v_{ij}&:=a_i^*a_j^*\bv \nonumber \\
v_{ijk}&:=a_i^*a_j^*a_k^*\bv\nonumber \\
w&:=a_1^*a_2^*a_3^*a_4^*\bv.
\end{align}

\begin{prop}\label{prop:pairb}\cite{FFR}
The pairing of any two vectors listed in  \eqref{circnote} is zero unless the total number of creation operators applied to $\bv$ sums to 4. 
\end{prop}

\begin{thm}\label{lem:circprod}
For  $i,j,k,l\in \{1,2,3,4\}$, the $\circ$ product has the following behavior on the Chevalley algebra $\cee$, in addition to the behavior expressed in Proposition \ref{prop:circprod}. 
\begin{enumerate}[\rm (a)]
\item $v\circ v_i=0$
\item $v\circ v_{ijk}=\mbox{sgn}_{(lijk)}a_l$
\item $v_i\circ v_{ij}= 0$ when any subscript, either $i$ or $j$ is repeated
\item $v_i\circ v_{jk} = \mbox{sgn}_{(iljk)}a_l$ for all distinct subscripts, $i,j,k$
\item $v_i\circ w = a_i^*$
\item $v_{ij}\circ v_{ijk}=0$ whenever two of the subscripts repeated, regardless of their order
\item $v_{ij}\circ v_{jkl}=\mbox{sgn}_{(jikl)}a_j^*$ when exactly one subscript, $j$ is repeated, regardless of the order of the subscripts
\item $v_{ijk}\circ w=0$.
\end{enumerate}
\end{thm}

\begin{prf}\label{prf:circrod}
We will illustrate these results by showing two carefully selected parts. \begin{enumerate}[(c)]
\item By Proposition \ref{prop:circprod}  $v_i\circ v_{ij} \in A$. Let $u$ be in  $A$. Using \eqref{circprod} and \eqref{cliffpairing}, the non-degenerate pairing $(v_i\circ v_{ij}, u)$  can only be identically for exactly zero for  $u=0$. But $(v_i\circ v_{ij}, u)= (v_{ij}\circ v_{i}, u)=(v_{ij}, v_{i}\circ u)= (v_{ij}, uv_{i})$, and $ uv_i \in$ span$\{ v, v_{ij}, v_{ik} \}$ for $1\leq k \leq 4, k \neq i $.  Now we apply Proposition \ref{prop:pairb} and conclude that $(v_{ij}, uv_{i})=0$ for all $u$ implying that  $v_i\circ v_{ij}=0.$
\end{enumerate}
\begin{enumerate}[(d)]
\item  Let $u$ be one of the 8 basis elements in $\mathcal A$. The pairing  $(v_i\circ v_{jk}, u)=(v_{i},uv_{jk})$  can only be non-zero for exactly one of the $u\in\mathcal A$. Moreover, $uv_{jk} \in  \{ v_{ljk}, v_{ijk}, v_{j}, v_k  \}$ and by Proposition \ref{prop:pairb}  the pairing vanishes unless $u=a_l^*$. Setting $(v_i\circ v_{jk}, u)=1$ and  by the uniqueness of $u$ in   \eqref{circprod} we conclude that $v_i\circ v_{jk}=sgn_{(iljk)} a_l.$ We determine the sign of the permutation in the standard fashion; the parity of the number of transpositions needed to put $(iljk)$ in ascending order. 
 
The omitted parts are proved similarly to either of these that we have provided. \qed
\end{enumerate}\end{prf}

\begin{prop}\label{prop:sigmaaction}\cite{FFR}
The action of $\sigma$ and $\sigma^2$ on the basis  $\mathcal A$ can be summarized in the following manner, and hence gives a description of the action of $\sigma$ on $\cee$.
$$\begin{array}{|lclcl|}
\hline
V^{(1)}& \rightarrow & V^{(2)} & \rightarrow & V^{(3)}\\
\hline
a_1 & & \phm v &  & -v_4\\
a_2 & & \phm v_{34} &  & \phm v_3\\
a_3 & & -v_{24} & & -v_2\\
a_4 & & \phm v_{14} &  & -v_{234}\\
\hline
\end{array}
\begin{array}{|lclcl|}
\hline
V^{(1)}& \rightarrow & V^{(2)} & \rightarrow & V^{(3)}\\
\hline
a_1^* & & \phm w &  & \phm v_{123}\\
a_2^* & & -v_{12} &  & \phm v_{124}\\
a_3^* & & -v_{13} &  & \phm v_{134}\\
a_4^* & & -v_{23} & & -v_1\\
\hline
\end{array}
$$
\end{prop}

\begin{rmk}\label{rmk:d4roots}
Using the notation of \eqref{finmods} we will choose the following root vectors of the $\g^{(i)}$ for $0\leq i \leq 2$, isomorphic finite dimensional Lie algebras of type $D_4$. The basis of the CSA's will be set according to the following.
$$\begin{array}{|rrr|}
\hline
\g ^{(0)}& \g ^{(1)}& \g ^{(2)}\\
\hline
\color{red}\fno a_1a_1^*\fno  &\color{red} \fno vw\fno  &\color{red} -\fno v_{4}v_{123}\fno \\
\color{red}\fno a_2a_2^*\fno &\color{red} -\fno v_{34}v_{12}\fno  &\color{red} \fno v_{3}v_{124}\fno \\
\color{red}\fno a_3a_3^*\fno &\color{red} \fno v_{24}v_{13}\fno  &\color{red} -\fno v_{2}v_{134}\fno \\
\color{red}\fno a_4a_4^*\fno &\color{red} -\fno v_{14}v_{23}\fno  &\color{red} \fno v_{234}v_{1}\fno \\
\hline
\end{array}$$
Positive root vectors are on the left side, and negative root vectors are on the right side of the following chart.  
$$\begin{array}{|rrr|}
\hline
\g ^{(0)}& \g ^{(1)}& \g ^{(2)}\\
\hline
\color{red}\fno a_1a_2\fno &\color{red} \fno vv_{34}\fno  &\color{red} -\fno v_{4}v_{3}\fno \\
\color{red}\fno a_1a_3\fno &\color{red} -\fno vv_{24}\fno  &\color{red} \fno v_{4}v_{2}\fno \\
\color{purple}\fno a_2a_3^*\fno &\color{purple} -\fno v_{34}v_{13}\fno &\color{purple} \fno v_{3}v_{134}\fno \\
\color{black}\fno a_1a_2^*\fno &\color{red} -\fno vv_{12}\fno  &\color{red} -\fno v_{4}v_{124}\fno  \\
\color{black}-\fno a_3a_4\fno &\color{red} \fno v_{24}v_{14}\fno  &\color{red} -\fno v_{2}v_{234}\fno \\
\color{red}\fno a_3a_4^*\fno &\color{red} \fno v_{24}v_{23}\fno  &\color{red} \fno v_{2}v_{1}\fno \\
\color{orange}\fno a_1a_3^*\fno &\color{orange} -\fno vv_{13}\fno  &\color{orange} -\fno v_{4}v_{134}\fno  \\
\color{orange}\fno a_2a_4\fno & \color{orange}\fno v_{34}v_{14}\fno  &\color{orange} -\fno v_{3}v_{234}\fno \\
\color{orange}-\fno a_2a_4^*\fno &\color{orange} \fno v_{34}v_{23}\fno  & \color{orange} \fno v_{3}v_{1}\fno \\
\color{brown}\fno a_2a_3\fno &\color{brown} -\fno v_{34}v_{24}\fno  &\color{brown} -\fno v_{3}v_{2}\fno \\
\color{brown}-\fno a_1a_4^*\fno &\color{brown}\fno vv_{23}\fno  &\color{brown} -\fno v_{4}v_{1}\fno \\
\color{brown}\fno a_1a_4\fno &\color{brown} \fno vv_{14}\fno  &\color{brown} \fno v_{4}v_{234}\fno \\
\hline
\end{array}
\begin{array}{|rrr|}
\hline
\g ^{(0)}& \g ^{(1)}& \g ^{(2)}\\
\hline
\color{red}\fno a_1^*a_2^*\fno &\color{red} -\fno wv_{12}\fno  &\color{red} \fno v_{123}v_{124}\fno \\
\color{red}\fno a_1^*a_3^*\fno &\color{red} -\fno wv_{13}\fno  &\color{red} \fno v_{123}v_{134}\fno \\
\color{purple}\fno a_2^*a_3\fno &\color{purple} \fno v_{12}v_{24}\fno &\color{purple} -\fno v_{124}v_{2}\fno \\
\color{red}\fno a_1^*a_2\fno &\color{red} \fno wv_{34}\fno  &\color{red} \fno v_{123}v_{3}\fno  \\
\color{red}-\fno a_3^*a_4^*\fno &\color{red} -\fno v_{13}v_{23}\fno  &\color{red} \fno v_{134}v_{1}\fno \\
\color{red}\fno a_3^*a_4\fno &\color{red} -\fno v_{13}v_{14}\fno  &\color{red} -\fno v_{134}v_{234}\fno \\
\color{orange}\fno a_1^*a_3\fno &\color{orange} -\fno wv_{24}\fno  &\color{orange} -\fno v_{123}v_{2}\fno  \\
\color{orange}\fno a_2^*a_4^*\fno & \color{orange}\fno v_{12}v_{23}\fno  &\color{orange} -\fno v_{124}v_{1}\fno \\
\color{orange}-\fno a_2^*a_4\fno &\color{orange} \fno v_{12}v_{14}\fno  & \color{orange} \fno v_{124}v_{234}\fno \\
\color{brown}\fno a_2^*a_3^*\fno &\color{brown} \fno v_{12}v_{13}\fno  &\color{brown} \fno v_{124}v_{134}\fno \\
\color{brown}-\fno a_1^*a_4\fno &\color{brown}-\fno wv_{14}\fno  &\color{brown} \fno v_{123}v_{234}\fno \\
\color{brown}\fno a_1^*a_4^*\fno &\color{brown} -\fno wv_{23}\fno  &\color{brown} -\fno v_{123}v_{1}\fno \\
\hline\end{array}$$
\end{rmk}

\begin{lem}\label{lem:identification}
For $a,b\in \mathcal A$, the action of $\fno ab\fno $, $\sigma\fno ab\fno $, and $\sigma^2\fno ab\fno $ may not agree on $\cee$, however since $\g ^{(0)}\cong \g ^{(1)}\cong \g ^{(2)}$, there is an identification of each of the operators in $\g ^{(1)}$ and $\g ^{(2)}$ with a linear combination of the operators in $\g ^{(0)}$. \begin{enumerate}[\rm(a)]
\item The identification of the operators representing the positive root vectors is given as follows.
$$\begin{array}{|rcrcr|}
\hline
\g ^{(1)}&\rightleftharpoons& \g ^{(0)} &\rightleftharpoons& \g ^{(2)}\\
\hline
\fno v_{}v_{34}\fno &&\color{red}\fno a_1a_2\fno &&-\fno v_{4}v_{3}\fno \\
-\fno v_{}v_{24}\fno &&\color{red}\fno a_1a_3\fno &&-\fno v_{4}v_{2}\fno \\
\fno v_{34}v_{13}\fno &&\color{purple}\fno a_2a_3^*\fno &&\fno v_{3}v_{134}\fno \\

-\fno v_{24}v_{23}\fno &&\color{red}\fno a_1a_2^*\fno &&\fno v_{2}v_{234}\fno \\
\fno v_{24}v_{14}\fno &&\color{red}\fno a_3a_4^*\fno &&-\fno v_{4}v_{124}\fno \\
\fno v_{}v_{12}\fno &&\color{red}\fno a_3a_4\fno &&\fno v_{2}v_{1}\fno \\
\hline
\end{array}
\begin{array}{|rcrcr|}
\hline
\g ^{(1)}&\rightleftharpoons& \g ^{(0)} &\rightleftharpoons& \g ^{(2)}\\
\hline
-\fno v_{}v_{13}\fno &&\color{orange}\fno a_2a_4\fno &&\fno v_{3}v_{1}\fno \\
-\fno v_{34}v_{14}\fno &&\color{orange}\fno a_2a_4^*\fno &&\fno v_{4}v_{134}\fno \\
\fno v_{34}v_{23}\fno &&\color{orange}\fno a_1a_3^*\fno &&-\fno v_{3}v_{234}\fno \\

\fno v_{34}v_{24}\fno &&\color{brown}\fno a_1a_4^*\fno &&-\fno v_{4}v_{234}\fno \\
\fno v_{}v_{23}\fno &&\color{brown}\fno a_1a_4\fno &&-\fno v_{3}v_{2}\fno \\
\fno v_{}v_{14}\fno &&\color{brown}\fno a_2a_3\fno &&-\fno v_{4}v_{1}\fno \\
\hline
\end{array}$$
\item The identification of the operators representing the negative root vectors is analogous.
\item The identification of the operators representing a basis of the CSA is given as follows. 

$\begin{array}{|lcl|}
\hline
\qquad \qquad \g ^{(1)}&\rightleftharpoons&\qquad \qquad \qquad  \g ^{(0)} \\
\hline
\sigma\fno a_1a_1^*\fno=\fno vw\fno &&\color{red}\frac{1}{2}(\fno a_1a_1^*\fno +\fno a_2a_2^*\fno +\fno a_3a_3^*\fno +\fno a_4a_4^*\fno )\\
\sigma\fno a_2a_2^*\fno=-\fno v_{34}v_{12}\fno  &&\color{red}\frac{1}{2}(\fno a_1a_1^*\fno +\fno a_2a_2^*\fno -\fno a_3a_3^*\fno -\fno a_4a_4^*\fno )\\
\sigma\fno a_3a_3^*\fno=\fno v_{24}v_{13}\fno  &&\color{red}\frac{1}{2}(\fno a_1a_1^*\fno -\fno a_2a_2^*\fno +\fno a_3a_3^*\fno -\fno a_4a_4^*\fno )\\
\sigma\fno a_4a_4^*\fno=-\fno v_{14}v_{23}\fno  &&\color{red}-\frac{1}{2}(\fno a_1a_1^*\fno -\fno a_2a_2^*\fno -\fno a_3a_3^*\fno +\fno a_4a_4^*\fno )\\
\hline
\hline
\qquad \qquad \g ^{(2)}&\rightleftharpoons&\qquad \qquad \qquad  \g ^{(0)}\\
\sigma^2\fno a_1a_1^*\fno=-\fno v_{4}v_{123}\fno &&\color{red}\frac{1}{2}(\fno a_1a_1^*\fno +\fno a_2a_2^*\fno +\fno a_3a_3^*\fno -\fno a_4a_4^*\fno )\\
\sigma^2\fno a_2a_2^*\fno=\fno v_{3}v_{124}\fno &&\color{red}\frac{1}{2}(\fno a_1a_1^*\fno +\fno a_2a_2^*\fno -\fno a_3a_3^*\fno +\fno a_4a_4^*\fno )\\
\sigma^2\fno a_3a_3^*\fno=-\fno v_{2}v_{134}\fno &&\color{red}\frac{1}{2}(\fno a_1a_1^*\fno -\fno a_2a_2^*\fno +\fno a_3a_3^*\fno +\fno a_4a_4^*\fno )\\
\sigma^2\fno a_4a_4^*\fno=\fno v_{234}v_{1}\fno &&\color{red}\frac{1}{2}(\fno a_1a_1^*\fno -\fno a_2a_2^*\fno -\fno a_3a_3^*\fno -\fno a_4a_4^*\fno )\\
\hline
\end{array}$
\end{enumerate}
\end{lem}
\begin{prf}\label{prf:identification}
This result is purely computational and can be proved by simply applying  each operator in $\g^{(i)}$ to each basis vector of  $V^{(j)}$ for $0\leq i \leq 2$ and $1\leq j \leq 3$.   
\qed\end{prf}

\begin{rmk}\label{rmk:matrixA}
Notice that Lemma \ref{lem:identification}(c) can be expressed as follows using the orthogonal matrix 
$$\mathscr A =[a_{ij}]= \one 
\left[ \begin{array}{rrrr} 1&1&1&-1\\ 1&1&-1&1\\1&-1&1&1\\1&-1&-1&-1  \end{array}
\right]$$ 
which represents $\sigma$ on the CSA with respect to the basis $\{h_1,h_2,h_3,h_4\}$ where
$h_j = \fno a_ja_j^*\fno$. We have $\sigma(h_j) = \sum_{i=1}^4 a_{ij} h_i$ is the linear combination given by column $j$ of $\mathscr A$, and since $\sigma^2 = \sigma^{-1}$ is represented by $\mathscr A^{-1} = \mathscr A^t$ so  $\sigma^2(h_j) = \sum_{i=1}^4 a_{ji} h_i$ is the linear combination given by row $j$ of $\mathscr A$.
\end{rmk}

According to Remark \ref{rmk:d4roots}  and \eqref{Gsystem} we now provide the root vectors for $\g_0$, the finite dimensional Lie algebra of type $G_2$ which occurs as the fixed points of $\g$ under $\sigma$. First we fix the root vectors of  our chosen CSA.
$$
\begin{array}{|lcl|}
\hline
H_1&= &\fno a_1a_1^* \fno + \fno a_3a_3^*\fno \\
H_2&= & \fno a_2a_2^* \fno - \fno a_3a_3^*\fno\\
\hline
\end{array} $$

The remaining root vectors are organized in the four blocks of the following table. The positive root vectors are on the left hand side while the negative roots are on the right. The blocks on the top contain the long root vectors, while the lower blocks contain the short root vectors.

$$
\begin{array}{|lcl|}
\hline
X_{\beta_1} & = & \fno a_2a_3^*\fno - \fno v_{34}v_{13} \fno + \fno v_{3}v_{134} \fno\\
X_{\beta_1+3\beta_2} & = & \fno a_1a_3\fno - \fno vv_{24} \fno +  \fno v_{4}v_{2} \fno\\
X_{2\beta_1+3\beta_2} & = & \fno a_1a_2\fno + \fno vv_{34} \fno - \fno v_{4}v_{3} \fno\\
\hline
X_{\beta_2} & = & \fno a_1a_2^*\fno - \fno vv_{12}\fno - \fno v_4v_{124}\fno\\
X_{\beta_1+\beta_2} & = & \fno a_1a_3^*\fno - \fno vv_{13}\fno - \fno v_4v_{134}\fno\\
X_{\beta_1+2\beta_2} & = & \fno a_2a_3\fno -   \fno v_{34}v_{24}\fno - \fno v_3v_2\fno\\
\hline
\end{array} 
\begin{array}{|lcl|}
\hline
X_{-\beta_1} & = & \fno a_3a_2^*\fno - \fno v_{12}v_{24} \fno + \fno v_{124}v_{2} \fno \\
X_{-\beta_1-3\beta_2} & = & \fno a_3^*a_1^*\fno + \fno wv_{13} \fno - \fno v_{123}v_{134} \fno \\
X_{-2\beta_1-3\beta_2} & = & \fno a_2^*a_1^*\fno + \fno wv_{12} \fno - \fno v_{123}v_{124} \fno \\
\hline
X_{-\beta_2} & = &\trd\left( \fno a_2a_1^*\fno - \fno wv_{34} \fno - \fno v_{123}v_{3} \fno \right)\\
X_{-\beta_1-\beta_2} & = & \trd\left(\fno a_3a_1^*\fno + \fno wv_{24} \fno + \fno v_{123}v_{2} \fno \right)\\
X_{-\beta_1-2\beta_2} & = & \trd\left(\fno a_3^*a_2^*\fno - \fno v_{12}v_{13} \fno - \fno v_{124}v_{134} \fno \right)\\
\hline
\end{array} $$

\begin{cor}\label{cor:G2roots}
The triality automorphism $\sigma$ induces a basis for $\g_0$ contained in $g^{(0)}$. 
 $$
\begin{array}{|lcl|}
\hline
H_1&= &\fno a_1a_1^* \fno + \fno a_3a_3^*\fno \\
H_2&= & \fno a_2a_2^* \fno - \fno a_3a_3^*\fno\\
\hline
X_{\beta_1} & = & \fno a_2a_3^*\fno \\
X_{\beta_1+3\beta_2} & = & \fno a_1a_3\fno \\
X_{2\beta_1+3\beta_2} & = & \fno a_1a_2\fno \\
\hline
X_{\beta_2} & = & \fno a_1a_2^*\fno + \fno a_3a_4^*\fno - \fno a_3a_4\fno\\
X_{\beta_1+\beta_2} & = & \fno a_1a_3^*\fno + \fno a_2a_4\fno - \fno a_2a_4^*\fno\\
X_{\beta_1+2\beta_2} & = & \fno a_2a_3\fno +  \fno a_1a_4\fno - \fno a_1a_4^*\fno\\
\hline
\end{array} 
\begin{array}{|lcl|}

\\
\\
\hline
X_{-\beta_1} & = & \fno a_3a_2^*\fno \\
X_{-\beta_1-3\beta_2} & = & \fno a_3^*a_1^*\fno \\
X_{-2\beta_1-3\beta_2} & = & \fno a_2^*a_1^*\fno \\
\hline
X_{-\beta_2} & = &\trd\left( \fno a_2a_1^*\fno + \fno a_4a_3^*\fno - \fno a_4^*a_3^*\fno\right)\\
X_{-\beta_1-\beta_2} & = & \trd\left(\fno a_3a_1^*\fno + \fno a_4^*a_2^*\fno - \fno a_4a_2^*\fno\right)\\
X_{-\beta_1-2\beta_2} & = & \trd\left(\fno a_3^*a_2^*\fno + \fno a_4^*a_1^*\fno - \fno a_4a_1^*\fno\right)\\
\hline
\end{array} $$\end{cor}

\begin{prf}\label{prf:G2roots}
This is immediate after applying Lemma \ref{lem:identification}. 
\qed\end{prf}

\begin{cor}\label{cor:G2mods}
The triality automorphism $\sigma$ induces bases for the $G_2$-modules $\g_1$ and $\g_2$ contained in $g^{(0)}$.

$\begin{array}{|c|}
\hline
\g_1\\
\hline
\fno a_1a_1^*\fno -\fno a_2a_2^*\fno -\fno a_3a_3^*\fno +(\xi^2-\xi)\fno a_4a_4^*\fno \\
\fno a_1a_2^*\fno  - \xi^2\fno a_3a_4\fno + \xi\fno a_3a_4^*\fno  \\
\fno a_1a_3^*\fno+ \xi^2\fno a_2a_4\fno -\xi\fno a_2a_4^*\fno  \\
\fno a_2a_3\fno -\xi^2\fno a_1a_4^*\fno +\xi\fno a_1a_4\fno \\
\fno a_1^*a_2\fno  - \xi^2\fno a_3^*a_4^*\fno + \xi\fno a_3^*a_4\fno \\
\fno a_1^*a_3\fno+ \xi^2\fno a_2^*a_4^*\fno -\xi\fno a_2^*a_4\fno  \\
\fno a_2^*a_3^*\fno -\xi^2\fno a_1^*a_4\fno +\xi\fno a_1^*a_4^*\fno\\
\hline
\end{array}
\begin{array}{|c|}
\hline
\g_2 \\
\hline
\fno a_1a_1^*\fno -\fno a_2a_2^*\fno -\fno a_3a_3^*\fno +(\xi-\xi^2)\fno a_4a_4^*\fno \\
\fno a_1a_2^*\fno  - \xi\fno a_3a_4\fno + \xi^2\fno a_3a_4^*\fno  \\
\fno a_1a_3^*\fno+ \xi\fno a_2a_4\fno -\xi^2\fno a_2a_4^*\fno  \\
\fno a_2a_3\fno -\xi\fno a_1a_4^*\fno +\xi^2\fno a_1a_4\fno \\
\fno a_1^*a_2\fno  - \xi\fno a_3^*a_4^*\fno + \xi^2\fno a_3^*a_4\fno  \\
\fno a_1^*a_3\fno+ \xi\fno a_2^*a_4^*\fno -\xi^2\fno a_2^*a_4\fno  \\
\fno a_2^*a_3^*\fno -\xi\fno a_1^*a_4\fno +\xi^2\fno a_1^*a_4^*\fno \\
\hline
\end{array}$

\end{cor}

\section{A Spinor Description  of $B_3$ as a Subalgebra of $D_4$}
\label{sec:B3}

In this section we use the spinor description of $D_4$ and the three order-two automorphisms $\tau, \sigma\tau$, and $\sigma^2\tau$ to give a spinor description of three copies of $B_3$ inside $D_4$.

\subsection*{Invariance of the Cartan Subalgebra Revisited}\label{sub:CartanB3}

This section is very similar to section \ref{sub:CartanG2} in that given the CSA from section \ref{sec:finspinconstruct}, $\h\subset \g $ with ordered orthonormal basis $\{h_1,...,h_4 \}$, we can find $\tau(\h) \subset \tau(\g ) = \g $. In fact, we have the analogous statement:

\begin{thm}\label{thm:} As algebras of operators on $(\mathcal C,\circ)$,  $\h = \tau(\h)$, and therefore, 
the triality group generated by $\sigma$ and $\tau$ leaves $\h$ invariant.
\end{thm}

\subsection*{Roots of $B_3$}\label{sub:rootsB3}
Recall the description of the root system for the finite dimensional rank 4 Lie algebra of type $D_4$ given in section \ref{sec:finspinconstruct};
$\Phi_{D_4} = \{\pm\ep_i\pm\ep_j\st 1\leq i<j\leq 4\}$  with simple roots
$\Delta_{D_4} = \{\alpha_i \st 1\leq i \leq 4\},$
and correspondence of root vectors with the roots given in \eqref{roots}. We now restrict $\tau$ to $\h$ to give the automorphism of the previous theorem $\tau: \h\rightarrow \h$, and induced dual map $\tau^*:\h^*\rightarrow\h^*$. \begin{align}
\tau^*(\alpha_1)=\alpha_1 ,\qquad
\tau^*(\alpha_2)=\alpha_2 ,\qquad
\tau^*(\alpha_3)=\alpha_4 ,\qquad
\tau^*(\alpha_4)=\alpha_3. \label{taustar}
\end{align}

The action of $\tau^*$ on the simple roots of $\g $ induces the action on the root spaces of $\g $ so that for $\mu \in \Phi_{D_4}$ \begin{align} \tau: \g _{\mu}\to \g _{\tau^*(\mu)}.\label{rootpermB3}\end{align}
One may also compute the action of $\tau^*$ on the individual $\ep_i$ via the appropriate use of \eqref{taustar}.

The fixed points of $\g $ under $\tau$ form a rank 3, $21$-dimensional simple Lie algebra of type $B_3$, with CSA denoted by $\mathfrak l$. We can restrict both $\tau$ and $\tau^*$,  to $\mathfrak l$ and its dual space $\mathfrak l^*$, to find the type $B_3$ root system $\Phi_{B_3}$ with two long simple roots and one short root given respectively as \be\label{Bsimple} \Delta_{B_3} = \{\gamma_1 = \alpha_1, \ph\ph \gamma_2 = \alpha_2, \ph\ph \gamma_3= \alpha_3-\alpha_4\}, \ee
which has the Dynkin diagram \qquad 
 \begin{tikzpicture}
\draw(4,0)--(5,0);
\draw (5,-.05)--(6,-.05);
\draw (5,.05)--(6,.05);
\draw(5.4,.2)--(5.6,0);
\draw(5.4,-.2)--(5.6,0);
\shade[ball color =white](4,0)circle(6pt);
\shade[ball color =white](5,0)circle(6pt);
\shade[ball color =white](6,0)circle(6pt);
\draw(4,.5)node[black]{$\gamma_1$};
\draw(5,.5)node[black]{$\gamma_2$};
\draw(6,.5)node[black]{$\gamma_3$};
\end{tikzpicture}.

Using either of the other order-two triality group elements, $\sigma\tau$ or $\sigma^2\tau$, one obtains two other copies of $B_3$ inside $D_4$, and corresponding versions of the results below. In fact, the $G_2$ subalgebra of $D_4$ found above in Section \ref{sec:G2} is the intersection of these three $B_3$ subalgebras. 

\subsection*{Decompositions of $\g$ and $\cee$ as $B_3$-Modules}\label{sub:irreptau}

There is also a reducible representation $\psi_{\mathcal C}: \g  \rightarrow \e(\mathcal C)$ similar to \eqref{repi}
\begin{align}\label{taurep} \tau\inv \circ\psi_{\cee}(x)\circ\tau = \psi_{\cee}(\tau\inv(x)).\end{align} We identify $\phi_{\cee}(\g )=\g $.

The action of $\tau$ on $\g $ decomposes it into the two eigenspaces \begin{align} \g = \mb _1\oplus\mb _{-1} \label{eigentau} \end{align} the fixed points, and the $-1$ eigenspace under $\tau$, respectively. Note that dim$(\mb _1) = 21$ and  dim$(\mb _{-1}) = 7$. It is actually the case that $\mb _1$ is a simple Lie algebra of type $B_3$ and $\mb _{-1}$ is its natural representation.

\subsection*{Root Vectors of type $B_3$}\label{sub:vectorsB3}

Here we find the root vectors of $\mb _1$ written as elements normal ordered elements in the spinor construction of $D_4$.

The action of $\tau:A\leftrightarrow A$ on the basis  $\mathcal A$ is given by the left two columns of the next table, while $\tau: \CM^0\leftrightarrow \CM^1$ via $\tau(b) = (a_4+a_4^*)b $ for all $b\in \CM^0\oplus \CM^1$ is given by the right two columns.
$$\begin{array}{|ccl|ccl|}
\hline
A & \leftrightarrow& A & A & \leftrightarrow& A \\
\hline
a_1 && -a_1 & a_1^* && -a_1^* \\
a_2 & & -a_2 & a_2^* && -a_2^* \\
a_3 & & -a_3 & a_3^* & & -a_3^* \\
a_4 & & \phm a_4^* & a_4^* & & \phm a_4 \\
\hline
\end{array}
\begin{array}{|ccl|ccl|}
\hline
\CM^0 & \leftrightarrow& \CM^1 & \CM^0 & \leftrightarrow& \CM^1 \\
\hline
v && v_4 & w & &- v_{123} \\
v_{12} && v_{124} & v_{14} & & -v_{1} \\
v_{13} & & v_{134} &v_{24} & & -v_{2}\\
v_{23} && v_{234} & v_{34} & & -v_{3}\\
\hline
\end{array}$$

As we did in the previous section for $G_2$, we fix the root vectors of  our chosen CSA for $\mb_1$, the Lie algebra of type $B_3$, and list the remaining root vectors which are organized similarly to Corollary \ref{cor:G2roots}  

\begin{cor}\label{cor:B3roots}
The triality automorphism $\tau$ induces a basis for $\mb_1$ contained in $g^{(0)}.$

$$
\begin{array}{|lcl|}
\hline
h_1 & = & \fno a_1a_1^* \fno\\
h_2 & = & \fno a_2a_2^*\fno \\
h_3 & = & \fno a_3a_3^*\fno \\
\hline
X_{\gamma_1} & = & \fno a_1a_2^*\fno \\
X_{\gamma_2} & = & \fno a_2a_3^*\fno \\
X_{\gamma_1+\gamma_2} & = & \fno a_1a_3^*\fno\\
X_{\gamma_2+2\gamma_3} & = & \fno a_2a_3\fno\\
X_{\gamma_1+\gamma_2+2\gamma_3} & = & \fno a_1a_3\fno\\
X_{\gamma_1+2\gamma_2+2\gamma_3} & = & \fno a_1a_2\fno\\
\hline
\end{array} 
\begin{array}{|lcl|}
\\
\\
\\
\hline
X_{-\gamma_1} & = & \fno a_2a_1^*\fno\\
X_{-\gamma_2} & = & \fno a_3a_2^*\fno\\
X_{-\gamma_1-\gamma_2} & = & \fno a_3a_1^*\fno\\
X_{-\gamma_2-2\gamma_3} & = & \fno a_3^*a_2^*\fno\\
X_{-\gamma_1-\gamma_2-2\gamma_3} & = & \fno a_3^*a_1^*\fno\\
X_{-\gamma_1-2\gamma_2-2\gamma_3} & = & \fno a_2^*a_1^*\fno\\
\hline
\end{array} $$$$
\begin{array}{|lcl|}
\hline
X_{\gamma_3} & = & \fno a_3a_4\fno-\fno a_3a_4^*\fno \\
X_{\gamma_2+\gamma_3} & = & \fno a_2a_4\fno-\fno a_2a_4^*\fno\\
X_{\gamma_1+\gamma_2+\gamma_3} & = & \fno a_1a_4\fno-\fno a_1a_4^*\fno\\
\hline
\end{array} 
\begin{array}{|lcl|}
\hline
X_{-\gamma_3} & = & \tfrac{1}{2}\left(\fno a_4^*a_3^*\fno-\fno a_4a_3^*\fno\right)\\
X_{-\gamma_2-\gamma_3} & = & \tfrac{1}{2}\left(\fno a_4^*a_2^*\fno-\fno a_4a_2^*\fno\right)\\
X_{-\gamma_1-\gamma_2-\gamma_3} & = & \tfrac{1}{2}\left(\fno a_4^*a_1^*\fno-\fno a_4a_1^*\fno\right)\\
\hline
\end{array} $$
\end{cor}

\vskip 5pt

\begin{cor}\label{cor:B3mods}
The triality automorphism $\tau$ induces a basis  of the $B_3$-module $\mb_{-1}$ contained in $g^{(0)}$.
$$\begin{array}{|cc|}
\hline
\qquad \qquad \qquad \qquad \mb _1&\\
\hline
 \fno a_4a_4^*\fno & \\
 \fno a_1a_4\fno + \fno a_1a_4^*\fno &  \fno a_1^*a_4\fno + \fno a_1^*a_4^*\fno\\
 \fno a_2a_4\fno + \fno a_2a_4^*\fno &  \fno a_2^*a_4\fno + \fno a_2^*a_4^*\fno\\
\fno a_3a_4\fno + \fno a_3a_4^*\fno&  \fno a_3^*a_4\fno + \fno a_3^*a_4^*\fno\\
\hline
\end{array}$$
\end{cor}

\section{Spinor Construction of $D_4^{(1)}$ and $D_4^{(1)}$-Modules}
\label{sec:affinespinconstruct}

Here we review Construction \ref{con:AffineSpinorConstruction} for the specific case of $\ell=4$ to get infinite dimensional versions of Cliff, $\CM$, and $\mathfrak{so}(8)$. 

Recall that $Z = \Z+\epsilon\in \left\{\Z, \hz\right\}$, $A = A^+\oplus A^-$ where $A^\pm\cong\C^{4}$ defined as in Section \ref{sec:finspinconstruct} with canonical bases  
$\{a_1,\cdots,a_4\}$ and $\{a^*_1,\cdots,a^*_4\}$ for $A^+$ and $A^-$, respectively, such that 
$(a_i,a_j) = 0 = (a^*_i,a^*_j)$ and $(a_i,a^*_j) = \delta_{i,j}$. 
Define $A(Z)$, with elements written $a(n)$, where $a\in A$ and $n\in Z$, and symmetric bilinear form  $(a(m),b(n))=(a,b)\dm$. Write $\clz = \Cl_4(Z)$ and note that its Clifford relations, given below, are a natural extension of the relations \eqref{maincliff} and \eqref{cliffsquare}  in the finite dimensional case. In particular, we have
\bea\label{cliffinf}
&a_i(m)a_j^*(n)+a_j^*(n)a_i(m)=\di\dm1 \notag \\
&a_i(m)a_j(n)+a_j(n)a_i(m)=0\notag\\
&a_i^*(m)a_j^*(n)+a_j^*(n)a_i^*(m)=0\textrm { for all } 1\leq i,j \leq4\notag\\
\mbox{ and } \quad & a_i(m)^2=a_j^*(n)^2=0.
\eea

The polarization $A(Z) = A(Z)^+\oplus A(Z)^-$ is used to define the left ideal
$\I(Z)$ in $\clz$ generated by $A(Z)^+.$ The irreducible left $\clz$-module $\cmz=\clz/\I(Z)$ has vacuum vector
$\bvz=1+\I(Z)$. 
This is similar to the finite dimensional case of the previous section in that $\CM(Z)$ is the cyclic left $\clz$-module generated by $\bvz$, where $\bvz$ is a non-zero element of the quotient such that $A(Z)^+\cdot \bvz=0$.
We have
\be \CM(Z)=  \mbox{ span } \{b_1(-n_1)\cdots b_r(-n_r)\bvz\st b_i(-n_i)\in A(Z)^-\}, \ee which is $\frac{1}{2}\mathbb Z$ graded. As in the finite dimensional case, $\cmz=\cmz^0\oplus \cmz^1$ is a sum of its even and odd subspaces, depending on the parity of $r$.
Recall the fermionic normal ordering $\fno a(m)b(n)\fno$ in \eqref{afno}
which is exactly as defined for the finite dimensional case in \eqref{fno} when $n=m=0$.

Using $\g=\mathfrak{so}(8)$ from Section \ref{sec:finspinconstruct} in the Definition \ref{def:KacMoody}, we get
the untwisted affine Kac-Moody algebra $\ghat$ of type $D_4^{(1)}$. 
 
These infinite dimensional Clifford modules play a central role in the spinor construction of level-$1$ representations of $\ghat$. Recall the notation,
\be \label{d4modules}
\hV^0 = \CM_4^0(\hz), \quad \hV^1 = \CM_4^1(\hz),\quad \hV^2 = \CM_4^0(\Z), \quad \hV^3 = \CM_4^1(\Z).
\ee
with highest weights 
\be \label{d4hws} 
\Lambda_0, \quad \Lambda_1 - \one\delta,\quad \Lambda_4  - \one\delta, \quad \Lambda_3  - \one\delta, \ee
and with highest weight vectors 
\be \label{d4hwvs}
\bv = \bv(\hz), \quad a_1(-\one)\bv, \quad \bv' = \bv(\Z), \quad a^*_4(0)\bv'.
\ee

The triality automorphism $\sigma$ can be lifted to $\shat: \ghat \to \ghat $ by $\shat(x(m)) = (\sigma(x))(m)$.  
We have seen that $\sigma$ is an automorphism of the $\g$-module $\cee = V^{(1)}\oplus V^{(2)} \oplus V^{(3)}$, cyclically permuting the summands. In \cite{FFR} it is shown that $\sigma$ lifts to $\shat: \hV \to \hV$ such that  $\shat(\hV^0) =\hV^0$ and $\shat$ cyclically permutes $\hV^1, \hV^2$, and $\hV^3.$ In particular, one knows that $\shat(\bv)=\bv$ and that $a_1(-\one)\bv$, $\bv'$, and $a_4^*(0)\bv'$ are cyclically permuted by $\shat$.

\newpage

The $\one\Z$-grading of $\hV$ is given by the eigenspaces $\hV_n$ of the Virasoro operator $L_0^Z$, and we refer to the eigenvalue $n$ as the \emph{depth}\index{depth} of the non-zero vectors in $\hV_n$. We have the following results on $\hV_0^0$, $\hV_1^0$, $\hV_{1/2}^i$ for $1\leq i\leq 3$. In the first table below we give a basis for $\hV_0^0$, $\hV_1^0$ and $\hV_{1/2}^1$.
$$
\begin{array}{|c|rr|c|}
\hline
\hV_n^i & Basis & & Dimension\\
\hline
\hV_0^0 &\bv&& 1\\
\hline
\hV_{1/2}^1 &a_i(-\one)\bv&1\leq i \leq 4& 4\\
            &a_i^*(-\one)\bv&1\leq i \leq 4& 4\\
\hline
\hV_1^0 &a_i(-\one)a_j^*(-\one)\bv& 1\leq i,j \leq 4& 16\\
&a_i(-\one)a_j(-\frac{1}{2})\bv&1\leq i<j \leq 4 & 6\\
&a_i^*(-\one)a_j^*(-\one)\bv&1\leq i<j \leq 4 & 6\\
\hline
\end{array}$$

The next table contains the same information for $\hV_{1/2}^2$ and $\hV_{1/2}^3$.
$$
\begin{array}{|c|rr|c|}
\hline
\hV_n^i & Basis & & Dimension\\
\hline
\hV_{1/2}^2 &\bv'& & 1\\
            &a_i^*(0)a_j^*(0)\bv'&1\leq i<j \leq 4& 6\\
    &a_1^*(0)a_2^*(0)a_3^*(0)a_4^*(0)\bv'& & 1\\
\hline
\hV_{1/2}^3 &a_i^*(0)\bv'&1\leq i \leq 4 & 4\\
    &a_i^*(0)a_j^*(0)a_k^*(0)\bv'&1\leq i<j<k \leq 4& 4\\
\hline
\end{array}$$

There are obvious $\g$-module isomorphisms $\hV^i_{1/2} \cong V^{(i)}$ for $1\leq i\leq 3$, using \eqref{naturalweights} and \eqref{spinorweights}, and $\hV^0_{1} \cong \g$ using \eqref{roots}. 
This means that there are subspaces of $\hV^0_{1}$ corresponding to the subalgebra $\mb_1$ of type $B_3$
and $\g_0$ of type $G_2$. (See \eqref{eigentau} and \eqref{eigensigma}.) Vertex operators 
$Y_k(v)$ for $v\in\hV^0_1$ such that $\shat(v)=v$ represent the affine subalgebra $\ghat_0$ of $\ghat$ of type $G_2^{(1)}$.

\section{Sugawara Constructions}\index{Sugawara construction}
\label{sec:Sugawara}

Up until this point, we have been satisfied with the definition of  fermionic normal ordering  on pairs of Clifford elements. We now note that the definition can be expanded to cover the normal orderings of any finite product of Clifford generators. For $Z=\hz$ the definition will be straight forward, however when $Z=\Z$, where there is the possibility of zero modes being present, we need to be more careful when extending the definition. 
  
\begin{dfn} \label{dfn:generalfno} \index{normal order!fermionic!expanded} For $a_{i_1}, ..., a_{i_r} \in \mathcal A$ and  $n_{1}, ..., n_{r}\in \hz $ we have the \emph{expanded fermionic normal ordering}
\be \label{generalfnohz}\fno a_{i_1}(n_1)\cdots a_{i_r}(n_r) \fno = sgn(\varphi)a_{i_{\varphi1}}(n_{\varphi1})\cdots a_{i_{\varphi r}}(n_{\varphi r}) \ee
where $\varphi$ is any permutation of the symbols $\{1,...,r\}$ such that $n_{\varphi1}\leq \cdots \leq n_{\varphi r}$. 
\newline
For  $n_{1}, ..., n_{r}\in \Z $, we want the identity $\fno a_{i_1}(n_1)\cdots a_{i_r}(n_r) \fno = sgn(\varphi)\fno a_{i_{\varphi1}}(n_{\varphi1})\cdots a_{i_{\varphi r}}(n_{\varphi r})\fno$ for any permutation $\varphi$. It is enough to make the following definition for  $n_{1}\leq \cdots \leq n_{r}$.  
Whenever each $0\neq n_i \in\Z$, we have 
\be \label{generalfnoZ}\fno a_{i_1}(n_1)\cdots a_{i_r}(n_r) \fno = a_{i_{1}}(n_{1})\cdots a_{i_{r}}(n_{r}). \ee
 However if $n_{j-1} < 0 = n_j =n_{j+1}=...=n_{j+s-1} < n_{j+s}$ for some $s\geq 1$, then define 
\be \label{generalfnoZ0}\fno a_{i_1}(n_1)\cdots a_{i_r}(n_r) \fno = a_{i_{1}}(n_{1})\cdots a_{i_{j-1}}(n_{j-1})\Big(a_{i_j}\fno(0) \cdots a_{i_j+s-1}(0)\fno\Big) a_{i_{j+s}}(n_{j+s})\cdots a_{i_{r}}(n_{r}) \ee
where for any $a_{i_1}, ...,  a_{i_k}\in \mathcal A$,
\be\label{generalfnoZno} \fno a_{i_1}(0)\cdots  a_{i_k}(0)\fno = \frac{1}{k!}\sum_{\varphi\in S_k}sgn(\varphi)a_{i_{\varphi1}}(0)\cdots a_{i_{\varphi k}}(0). \ee 
 \end{dfn}

\begin{rmk}\label{rmk:zero modes} We refer the reader to \cite{FFR} for details on the above. For this investigation we will only need to use normal ordered products of two Clifford generators, referred to henceforth as  \index{normal order!fermionic!quadratic}\emph{quadratics}, and normal ordered products of four Clifford generators, referred to as  \index{normal order!fermionic!quartic}\emph{quartics}. When proving our main results for the Ramond modules we will introduce two lemmas which are specific instances of the definition described in \eqref{generalfnoZ} - \eqref{generalfnoZno}. 
\end{rmk} 
 
We will use the following abbreviations for vectors in $\hV_2^0$:
\be\label{depth2notations}
{\bf 1} = \bv,\quad i^\circledast j^\circledast  = \fno a_i^\circledast (-\tre) a_j^\circledast( -\one)\fno \bv, \quad 
i\ca j\ca m\ca n\ca = \fno a_i^\circledast (-\one) a_j^\circledast (-\one)a_m^\circledast (-\one) a_n^\circledast (-\one)\fno\bv\ee
for $1\leq i,j,m,n\leq 4$, where any occurence of the $\circledast $ can be either blank or $*$. 
Note that, because of the Clifford 
anti-commutation relations \eqref{cliffinf}, we have, for example, $ij = -ji$ and $ijmn = mijn = mnij$. 
In this notation we can write the conformal vector \eqref{omegadell} with $\ell=4$ as 
\be\label{omegad4}\D= \frac{1}{2}\sum_{i=1}^4 \Big(ii^*+i^*i\Big).\ee

This section provides the necessary details for several constructions of Sugawara operators representing the Virasoro algebra in the spinor construction of $\hV$. We have already seen how the vertex operators $Y_m(\D)$  represent $\Vir$ with central charge $c=4$ on $\hV$. 
Here we use the full Sugawara construction for the subalgebras of types $B_3$ and $G_2$ to obtain vertex operators
$Y_m(\B)$ and $Y_m(\G)$ that also represent $\Vir$ on $\hV$, but with central charges $7/2$ and $14/5$, respectively. These will be used to obtain coset Virasoro operators $Y_m(\DB)$ and $Y_m(\BG)$ 
which also represent $\Vir$ on $\hV$, but with central charges, $1/2$ and $7/10$, respectively. By Theorem 
\ref{thm:cosetvir}, these two coset Virasoro operators commute with each other, and with the operators representing the affine subalgebra of $\ghat$ of type $G_2^{(1)}$, $Y_k(v)$ for $v\in\hV^0_1$ such that $\shat(v)=v$. Their purpose in this work is to provide the decomposition of each irreducible $\ghat$-module $\hV^i$ into the direct sum of tensor products $\Vir(1/2,h_1)\otimes \Vir(7/10,h_2)\otimes W(\Omega_j)$ for $h_1\in\{0,1/2,1/16\}$,  $h_2\in\{0,1/10,3/5, 3/2, 7/16, 3/80 \}$, where $W(\Omega_j)$, $j=0,2$, are the two level-1 irreducible 
$G_2^{(1)}$-modules. 

There are certain linear combinations of operators which we will use so frequently that we introduce  the following notations to abbreviate these long expressions.

\begin{align}\label{abreve}
\fb{$44^*$} &= 44^*+4^*4+44+4^*4^* \notag\\
\fb{$11^*22^*$} &= 11^*22^*+11^*33^*-22^*33^*\notag \\
\fb{$22^*44^*$} &= 22^*44^*+33^*44^*-11^*44^* \notag\\
\fb{$1^*234$} &= 1^*234 + 12^*3^*4^* \notag\\
\fb{$1^*234^*$} &= 12^*3^*4 + 1^*234^*.
\end{align}

\subsection{The Conformal Vector $\B$}\label{sub:sugB3}

Recall  $\mb_1$ is the Lie algebra with basis described in Corollary \ref{cor:B3roots}.  Let $\{\epsilon_1,\epsilon_2,\epsilon_3\}$ be an orthonormal basis of $\mathbb R^3$ and let $\Phi_{B_3}$ be a root system of type $B_3$ with simple roots $\Delta_{B_3}= \{\gamma_1=\epsilon_1-\epsilon_2, \quad \gamma_2=\epsilon_2-\epsilon_3, \quad \gamma_3=\epsilon_3 \}$ as in \eqref{Bsimple} . Denote the set of long roots as $\Phi(L)$ and the short roots as $\Phi(S)$. Denote the CSA as $C$, the set of long root vectors as $L$ and the set of short root vectors as $S$.

We follow the prescription for the construction of Sugawara operators (\ref{con:Sugawara}) associated with $\mb_1$ on $\hV$ as components of the vertex operator $Y(\B,w)$. The basis of $\mb_1$ given in Corollary \ref{cor:B3roots} is $\{h_1,h_2,h_3\}\cup\{X_\gamma\st \gamma\in\Phi_{B_3}\} = \{X_i\st 1\leq i \leq dim(\mb_1)=21\}$. 
For each root $\gamma$, $X_{\gamma}$ and $X_{-\gamma}$ are dual with respect to the form \eqref{spinkill}, and we write $X^{\gamma} = X_{-\gamma}$ for the dual of $X_\gamma$. We have seen that $(h_i,h_j)=\di$, so these basis vectors of the CSA are self dual, and we write the dual $h^i=h_i$. We unify these notations by writing 
$X^i$ for the dual of $X_i$, $1\leq i\leq 21$. The conformal vector is then 
\be\label{sugB3} \B = L^{B_3}_{-2} {\bf 1} =\ds \mathscr S_{\B}\sum_{i=1}^{21}\sum_{k\in\Z}\no X_i(-k)X^i(-2+k)\no \textbf{1}.\ee 
Since the dual Coxeter number of $\mb_1$ of type $B_3$ is $5$ and the level of each $\hV^i$ is $1$, the scalar factor is $\mathscr S_{\B}= \frac{1}{12}$ and from \eqref{charge} the central charge is $c_{B_3} = \frac{dim(\mb_1)\cdot 1}{5+1} = \frac{21}{6} = \frac{7}{2}$.

Whenever $-k>0$ or $-2+k>0$ we have $\no X_i(\mu)X^i(\nu)\no\textbf{1}=0$. Hence we can reduce \eqref{sugB3} to the finite sum
\be\label{suga}\ds \B =\ds \mathscr S_{\B}\left( 2\sum_{i=1}^{21}X_i(-2)X^i(0) \textbf{1}+\sum_{i=1}^{21}X_i(-1)X^i(-1) \textbf{1}\right).\ee

To simplify the first summation above, consider 
\be\label{reduction}X^i(0) \textbf{1}=\fno ab\fno_0\textbf{1}= \ds\sum_{k\in\hz} \fno a(k)b(-k)\fno \textbf{1} = 0\ee
using \eqref{normk} and the fact that Clifford operators with positive mode numbers annihilate \textbf{1}. 
Therefore, we now have the simplified formula 
\be\label{shortsuga}\ds \B =\ds \mathscr S_{\B}\left(\sum_{i=1}^{21}X_i(-1)X^i(-1) \textbf{1}\right).\ee

We will break up this summation into three parts, corresponding to the CSA, $C$, the long roots, $\Phi(L)$, and the short roots, $\Phi(S)$, and we will temporarily leave out the scalar $\tfrac{1}{12}$. The terms corresponding to the CSA give
\begin{align}
\omega_H=&\sum_{i=1}^3h_i(-1)h^i(-1)\textbf{1} \nonumber\\
=&\sum_{i=1}^3\fno a_ia_i^*\fno_{-1}\left(\fno a_ia_i^*\fno_{-1}\textbf{1}\right)\nonumber\\
= &\sum_{i=1}^3\fno a_ia_i^*\fno_{-1} \left(a_i(-\tfrac{1}{2})a_i^*(-\tfrac{1}{2}) \textbf{1}\right) \nonumber\\
= &\sum_{i=1}^3\left(\sum_{k\in K} \fno a_i(k)a_i^*(-1-k)\fno \right) \left(a_i(-\tfrac{1}{2})a_i^*(-\tfrac{1}{2}) \textbf{1}\right)\nonumber\\
=& \sum_{i=1}^3 (ii^* + i^*i) \notag \\
= &11^*+1^*1 + 22^*+2^*2 + 33^*+3^*3.\label{omegahB3}
\end{align}
\begin{rmk}\label{rmk:short}
Note that we have significantly shortened the sequence of steps that it takes to get to line \eqref{omegahB3}. The computations are routine, and in fact well known from the literature. 
\end{rmk}
Using similar techniques we can compute the contribution of the long and short root vectors to the conformal vector. 
\begin{align}
\omega_{L} &= \sum_{\gamma\in\Phi(L)}X_{\gamma}(-1)X^{\gamma}(-1)\textbf{1} \nonumber\\
&= 4\sum_{i=1}^3 (ii^*+i^*i) \nonumber\\
& = 4\times\left(11^*+1^*1 + 22^*+2^*2 + 33^*+3^*3\right)\label{omegalB3}. \\
 \hbox{ and } \notag \\
\omega_{S} &= \sum_{\gamma\in\Phi(S)}X_{\gamma}(-1)X^{\gamma}(-1)\textbf{1} \nonumber\\
&= \sum_{i=1}^3 (ii^*+i^*i) +3\times(44^*+4^*4 -44 -4^*4^*) \nonumber\\
& = \left(11^*+1^*1 + 22^*+2^*2 + 33^*+3^*3\right) +3\times(44^*+4^*4 -44 -4^*4^*)\label{omegasB3}.
\end{align}
We now combine the results of \eqref{omegahB3},\eqref{omegalB3}, and \eqref{omegasB3} to obtain 
\begin{align}\label{omegab3}
{\B}&= \tfrac{1}{12}\left(\omega_H + \omega_L + \omega_S\right)\nonumber\\
& =  \tfrac{1}{12}\Big(11^*+1^*1 + 22^*+2^*2 + 33^*+3^*3 \nonumber\\
& +  4\times\left(11^*+1^*1 + 22^*+2^*2 + 33^*+3^*3 \right)\nonumber\\
& + \left(11^*+1^*1 + 22^*+2^*2 + 33^*+3^*3\right) +3\times(44^*+4^*4 -44 -4^*4^*)\Big) \nonumber\\
& = \one\times (11^*+1^*1 + 22^*+2^*2 + 33^*+3^*3 ) + \tfrac{1}{4}\times(44^*+4^*4 -44 -4^*4^*) \nonumber\\
&=\D - \tfrac{1}{4}\fbf.
\end{align}

\subsection{The Conformal Vector $\G$}\label{sub:sugG2}

This section is analogous to Section \ref{sub:sugB3}, so we will be even more brief. 

Let $\{\epsilon_1,\epsilon_2,\epsilon_3\}$ be an orthonormal basis of $\mathbb R^3$ and let $\Phi_{G_2}$ be a root system of type $G_2$ with the simple roots $\Delta_{G_2} = \{\beta_1 = \ep_2-\ep_3\, \ph\ph \beta_2 =\tfrac{1}{3}( \ep_1-\ep_2+2\ep_3)\}$  as in \eqref{Gsimple}.  In this section, $\{X_i\st 1\leq i \leq 14\}$ is the basis of the subalgebra $\g_0$ of type $G_2$ given in Corollary \ref{cor:G2roots}, the dual Coxeter number is $4$, so that the scalar factor is $\mathscr S_{\G} = \tfrac{1}{10}$ and from \eqref{charge} the central charge is 
$c_{G_2} = \frac{dim(\g_0)\cdot 1}{4+1} = \frac{14}{5}$.

We construct the conformal vector \be\omega_{G_2}=\ds \mathscr S_{G_2} \sum_{i=1}^{14}\sum_{k\in\Z}\no 
X_i(-k)X^i(-2+k)\no \textbf{1}\label{sugG2}\end{equation}
which simplifies as above, giving

\begin{align}\label{omegag2}
{\G}&=\tfrac{1}{10}( \omega_C + \omega_L + \omega_S)\nonumber\\
&= \tfrac{1}{10}\Big( \two(11^*+1^*1 + 22^*+2^*2 + 33^*+3^*3) + \two(11^*22^* + 11^*33^* - 22^*33^*)\nonumber\\
&+ \tfrac{6}{3}(11^*+1^*1 + 22^*+2^*2 + 33^*+3^*3) + \tfrac{6}{3}(11^*22^* + 11^*33^* - 22^*33^*)\nonumber\\
&+  \tfrac{4}{3}(11^*+1^*1 + 22^*+2^*2 + 33^*+3^*3) - \two(11^*22^*+11^*33^*-22^*33^*) \nonumber\\
&+ \tfrac{6}{3}(1^*234 + 12^*3^*4^* - 1^*234^* - 12^*3^*4)+ \tfrac{6}{3} (44^*+4^*4-44-4^*4^*)\Big) \nonumber\\
& = \tfrac{1}{5}\Big(4\D-\fbf +\fbo +\fbs-\fbss\Big) 
\end{align}

\begin{thm}\label{thm:cosets} The vertex operators $Y_m(\DB) = L_m^{1/2}$, $m\in\Z$, representing the 
coset Virasoro algebra on $\hV$ with central charge $\frac{1}{2}$ are provided by the conformal vector
\be\label{OMDB}
\DB = \D - \B 
= \tfrac{1}{4} \framebox{$44^*$} \ee
and the vertex operators $Y_m(\BG) = L_m^{7/10}$, $m\in\Z$, representing the 
coset Virasoro algebra on $\hV$ with central charge $\frac{7}{10}$ are provided by the conformal vector  
\be\label{OMBG}
\BG = \B - \G
= \tfrac{1}{5}\D - \tfrac{1}{20}\fb{$44^*$} - \tfrac{1}{5}\fb{$11^*22^*$} - \tfrac{1}{5}\fb{$1^*234$} + \tfrac{1}{5}\fb{$1^*234^*$}.
\ee
The operators from these two coset Virasoros commute with each other, $[L_m^{1/2}, L_n^{7/10}] = 0$ for all $m,n\in\Z$, and both commute with all operators representing $G_2^{(1)}$.
\end{thm}
\begin{prf}  This is clear from Theorem \ref{thm:cosetvir}. The details have been presented above in \eqref{omegad4}, \eqref{omegab3}, \eqref{omegag2}, and the central charges are
$c_{D_4-B_3} = 4 - \frac{7}{2} = \frac{1}{2}$ and $c_{B_3-G_2} = \frac{7}{2} - \frac{14}{5} = \frac{7}{10}$, respectively. 
\qed \end{prf}

\section{Operators Used in the Proof of the Main Theorems}\label{sub:operators}
\subsection{Affine Operators Representing $G_2^{(1)}$}\label{ssub:G2}
There are three operators in $G_2^{(1)}$ corresponding to the three positive simple roots:
\begin{align}
X_{\beta_1}(0) &= \fno a_2(z)a_3^*(z)\fno_0  = \sum_{r\in Z} \fno 2(r)3^*(-r)\fno \label{fnobeta1}\\
X_{\beta_2}(0) &= \fno a_1(z)a_2^*(z)\fno_0  +  \fno a_3(z)a_4^*(z)\fno_0 -  \fno a_3(z)a_4(z)\fno_0 \nonumber\\
                           &= \sum_{r\in Z} \Big( \fno 1(r)2^*(-r)\fno  +  \fno 3(r)4^*(-r)\fno -  \fno 3(r)4(-r)\fno \Big)\label{fnobeta2}\\
X_{-\theta}(1) &= \fno a_2^*(z)a_1^*(z)\fno_1 =  \sum_{r\in Z} \fno 2^*(r)1^*(1-r)\fno. \label{fnotheta}
\end{align}
In these expressions, fermionic normal ordering switches the order of the operators (with a minus sign) if $r > 0$, and leaves the order alone if $r \leq 0$. But since the Clifford operators $i\ca(r)$ and $j\ca(s)$ anti-commute for $i\neq j$, we have $\fno i\ca(r)j\ca(s)\fno = i\ca(r)j\ca(s) = - j\ca(s)i\ca(r)$.

\begin{thm}\label{thm:g2ops} We may write the $G_2^{(1)}$ operators \eqref{fnobeta1},  \eqref{fnobeta2}, and \eqref{fnotheta}  as
\begin{align} 
X_{\beta_1}(0) &= \sum_{r\in Z}  2(r)3^*(-r) \label{usebeta1} \\
X_{\beta_2}(0) &= \sum_{r\in Z} \Big( 1(r)2^*(-r)  +  3(r)4^*(-r) -  3(r)4(-r) \Big) \label{usebeta2} \\
X_{-\theta}(1) &=  \sum_{r\in Z}  2^*(r)1^*(1-r). \label{usetheta}
\end{align}
\end{thm}

\begin{rmk}\label{rmk:G2ops}
In the following chapter we will only need to apply these operators to vectors in $\hV_n$ for depth $0\leq n\leq 2$.
When these operators are applied to those vectors, there are only finitely many values of $r$ which may yield nonzero results since the annihilation operators may be anti-commuted to act first. The main Clifford relations reduce these summations to $r\in \{\pm\one,\pm\tre\}$ in the Neveu-Schwarz case $Z = \hz$, and 
$r =0$ in the Ramond case $Z = \Z$.
\end{rmk}

\subsection{Virasoro Representions on the Neveu-Schwarz Module}\label{ssub:NS}

\begin{thm}\label{lem:LkNSops} The following vertex operators provide two coset representations of the Virasoro algebra on the Neveu-Schwarz module, $CM(\hz)$, with central charges $\one$ and $\sev$, respectively. 
For all $k\in\Z$, we have
\begin{align}
L_k^{1/2} = - \tfrac{1}{4}\sum_{r\in \hz}\left(r+\one\right)\Big(\fno 4(r)4(k-r) \fno + \fno 4^*(r)4^*(k-r) \fno +\fno 4(r)4^*(k-r) \fno + \fno 4^*(r)4(k-r) \fno \Big)
 \label{LkDBNS} \end{align}
and for $r_1,r_2,r_3,r_4\in\hz$, we have
\begin{align}
&L_k^{7/10} = -\tfrac{1}{10}\sum_{i=1}^4\sum_{r\in\hz}(r+\tfrac{1}{2})\Big(\fno i^*(r)i(k-r)\fno+\fno i(r)i^*(k-r) \fno\Big)\nonumber\\
&+\tfrac{1}{20}\sum_{r\in\hz}(r+\tfrac{1}{2})\Big(\fno 4(r)4(k-r)\fno + \fno 4^*(r)4^*(k-r)\fno+\fno 4^*(r)4(k-r)\fno+\fno 4(r)4^*(k-r)\fno\Big) \nonumber\\
&-\tfrac{1}{5}\sum_{r_1+r_2+r_3+r_4=k}\Big( \fno 1(r_1)1^*(r_2)2(r_3)2^*(r_4)\fno + \fno 1(r_1)1^*(r_2)3(r_3)3^*(r_4)\fno - \fno 2(r_1)2^*(r_2)3(r_3)3^*(r_4)\fno    \nonumber\\
&\qquad \qquad\qquad+ \fno 1^*(r_1)2(r_2)3(r_3)4(r_4)\fno + \fno 1(r_1)2^*(r_2)3^*(r_3)4^*(r_4)\fno \nonumber\\
&\qquad\qquad\qquad- \fno 1^*(r_1)2(r_2)3(r_3)4^*(r_4)\fno -\fno 1(r_1)2^*(r_2)3^*(r_3)4(r_4)\fno    \Big).\label{LkBGNS}
\end{align}
\end{thm}

\subsection{Virasoro Representions on the Ramond Module}\label{ssub:R}

For $v\in\CM(\hz) = \hV^0\oplus \hV^1$ the defintion of vertex operators $Y(v,z)$ acting on the Ramond module $\CM(\Z) = \hV^2\oplus \hV^3$ is modified as in \eqref{RamondVO} with $\ell=4$. 
These vertex operators provide the Ramond module $\CM(\Z)$ with the structure of a vertex operator superalgebra module, as shown in \cite{FFR}. In particular, for $v = \D$ we have 
\be\label{Deltad4} \exp(\Delta(z))\D = \D + \tfrac{1}{2}  z^{-2} {\bf 1}\ee
so 
\be\label{Ramondd4VirasoroVO} Y(\D,z) = {\overline Y}(\D,z) 
+ \tfrac{1}{2} z^{-2} {\overline Y}({\bf 1},z) 
= {\overline Y}(\D,z) + \tfrac{1}{2} z^{-2} I = L^{\Z}(z) \ee
where $I$ is the identity operator on $\CM(\Z)$. The effect of this ``correction'' is to add the scalar
$\tfrac{1}{2}$ to the unmodified operator ${\bar L}^\Z_0$, which explains the $-\one\delta$ in the weight of
$\bv'$ in \eqref{d4hws}. We will now establish formulas for the coset Virasoro operators $L^{1/2}(z)$ and 
$L^{7/10}(z)$ on the Ramond module.
 
\begin{lem}\label{lem:expDelta} We have
\bea \exp(\Delta(z))\DB &=& \DB + \tfrac{1}{16}z^{-2}{\bf 1} \qquad\hbox{and}\\
\exp(\Delta(z))\BG &=& \BG + \tfrac{7}{80}z^{-2}{\bf 1}\eea 
therefore, 
\bea
L^{1/2}(z) &=& \overline Y_{\Z}(\DB,z)+\tfrac{1}{16}z^{-2}I = \sum_{k\in\Z}L^{1/2}_kz^{-k-2}\\
L^{7/10}(z) &=& \overline Y_{\Z}(\BG,z)+\tfrac{7}{80}z^{-2}I = \sum_{k\in\Z}L^{1/2}_kz^{-k-2}.
\eea
\end{lem}
\begin{prf} In applying $\Delta(z)$ to $\DB$, note that the only possible non-zero contributions are
\bea &\for \Big(C_{0,1} a_4(\tre)a^*_4(\one) + C_{1,0} a_4(\one)a^*_4(\tre) \Big) 
\Big(a_4(-\tre)a^*_4(-\one){\bf 1} + a^*_4(-\tre)a_4(-\one){\bf 1} \Big) z^{-2}\notag \\
&= \for (C_{0,1}-C_{1,0}){\bf 1} z^{-2} = \tfrac{1}{16} {\bf 1} z^{-2}\eea
since $C_{0,1}= -C_{1,0} = \tfrac{1}{8}$. Since $\Delta(z){\bf 1} = 0$, we get the first equation.
In applying $\Delta(z)$ to $\BG$, the $(m,n)=(0,0)$ terms in $\Delta(z)$ that might give a non-zero contribution
when applied to the quartic terms in $\BG$ have coefficients $C_{0,0} = 0$. So we have
\be \Delta(z)\BG = \Delta(z) \Big( \tfrac{1}{5}\D - \tfrac{1}{20}\fb{$44^*$}\ \Big) =
\tfrac{1}{5} \Big(\tfrac{1}{2} z^{-2} {\bf 1}\Big) - \tfrac{1}{20}\Big(\tfrac{1}{4} z^{-2} {\bf 1}\Big) 
= \tfrac{7}{80} z^{-2} {\bf 1}
\ee
giving the second equation.
\qed\end{prf}

\begin{thm}\label{lem:LkRops} The following vertex operators provide two coset representations of the Virasoro algebra on the Ramond module, $\CM(\Z)$, with central charges $\one$ and $\sev$ respectively. 
For all $k\in\Z$, we have

$L_k^{1/2} = \tfrac{1}{16}\delta_{k,0}I $ 
\begin{align}
  - \tfrac{1}{4}\sum_{r\in \Z}\left(r+\one\right)\Big(\fno 4(r)4(k-r) \fno + \fno 4^*(r)4^*(k-r) \fno +\fno 4(r)4^*(k-r) \fno + \fno 4^*(r)4(k-r) \fno \Big) \label{LkDBR} 
\end{align}
and for $r_1,r_2,r_3,r_4\in\Z$, we have
\begin{align}
&L^{7/10}_k = \tfrac{7}{80}\delta_{k,0}I -\tfrac{1}{10}\sum_{i=1}^4\sum_{r\in\Z}(r+\tfrac{1}{2})\Big(\fno i^*(r)i(k-r)\fno+\fno i(r)i^*(k-r) \fno\Big)\nonumber\\
&+\tfrac{1}{20}\sum_{r\in\Z}(r+\tfrac{1}{2})\Big(\fno 4(r)4(k-r)\fno + \fno 4^*(r)4^*(k-r)\fno+\fno 4^*(r)4(k-r)\fno+\fno 4(r)4^*(k-r)\fno\Big) \nonumber\\
&-\tfrac{1}{5}\sum_{r_1+r_2+r_3+r_4=k}\Big( \fno 1(r_1)1^*(r_2)2(r_3)2^*(r_4)\fno + \fno 1(r_1)1^*(r_2)3(r_3)3^*(r_4)\fno - \fno 2(r_1)2^*(r_2)3(r_3)3^*(r_4)\fno    \nonumber\\
&\qquad\qquad\qquad+ \fno 1^*(r_1)2(r_2)3(r_3)4(r_4)\fno + \fno 1(r_1)2^*(r_2)3^*(r_3)4^*(r_4)\fno \nonumber\\
&\qquad\qquad\qquad- \fno 1^*(r_1)2(r_2)3(r_3)4^*(r_4)\fno -\fno 1(r_1)2^*(r_2)3^*(r_3)4(r_4)\fno    \Big).\label{LkBGR}
\end{align}
\end{thm}

\begin{rmk}\label{rmk:vops}
As in Remark \ref{rmk:G2ops}, $L_k^{1/2}$ acting on vectors of depth at most $2$ reduces to a sum over the same small set of $r$ since these operators are also quadratic. We handle the quadratic operators in $L_k^{7/10}$  similarly, however, the \emph{quartic} operators require a different and less obvious treatment. In the following chapter we verify that certain vectors $v$ of depth at most $2$ are highest weight vectors with respect to $G_2^{(1)}$ and the two coset representations of the Virasoro algebra. To do so we will only need to check the following:
 \begin{enumerate}[(a)]
 \item $L_k^{1/2}v = 0 = L_k^{7/10}v$ for $k\in\{1,2\}$.
 \item The operators in Theorem \ref{thm:g2ops} also annihilate $v$.
 \item  $L_0^{1/2}v = h^{1/2}_v v$, $L_0^{7/10}v = h^{7/10}_v v$ for certain rational numbers $h^{1/2}_v$ and $h^{7/10}_v$.
 \item  $v$ has $G_2^{(1)}$-weight $\om_j$ for $j\in\{0,2\}$.
\end{enumerate}
\end{rmk}

\begin{ex}\label{ex:depth3} We will now show how to treat the quartic terms in $L_k^{7/10}$ acting on a basis vector $v\in \hV^1_{3/2}$. Such a vector has the form $a_i(-\one)a_i^*(-\one)a_j\ca(-\one){\bf1}$,  $a_i^\circledast(-\one)a_j^\circledast(-\one)a_k^\circledast(-\one){\bf1}$, or $a_i^\circledast(-\tre){\bf1}$ as shown in Appendix \ref{sec:bases}. Any Clifford operator $b\ca(r_n)$ in a quartic term with mode number $r_n > \tre$ annihilates $v$ since it can anti-commute past the creation operators in $v$ and annihilate ${\bf1}$. Hence, we need only consider $\tre \geq r_n\in \hz$ for $1\leq n\leq 4$.

When $k=0= r_1+r_2+r_3+r_4$, and $v\in \hV^1_{3/2}$, the only quartic operators which may have a non-zero contribution when acting on $v$ have multiset $\{r_1,r_2,r_3,r_4\}$ equal to one of the following multisets:    \begin{align}\label{quartic} \Big\{-\one,-\one,+\one,+\one\Big\}, \quad \Big\{-\tre,+\one, +\one,+\one\Big\}, \quad \Big\{-\one,-\one,-\one,+\tre\Big\} \end{align} henceforth denoted  \begin{align} --++, \quad -+++, \quad ---+ .\end{align}
When the multiset is $--++$, there are six distinct quartic operators which contribute to the sum $\ds\sum_{r_1+r_2+r_3+r_4=0} \fno a_i\ca(r_1)a_j\ca(r_2)a_k\ca(r_3)a_l\ca(r_4)\fno$. After normal ordering has been applied, permuting the Clifford generators in each of those six quartic operators so that the mode numbers are in non-decreasing order, they are
\be\label{pattern1} 
i\ca j\ca k\ca l\ca \qquad  i\ca k\ca l\ca j\ca   \qquad i\ca l\ca j\ca k\ca  \qquad  j\ca k\ca i\ca l\ca   \qquad j\ca l\ca k\ca i\ca  \qquad k\ca l\ca i\ca j\ca .
\ee
We have shortened the notation for brevity, for example, by writing $i\ca j\ca k\ca l\ca$ for \newline $a_i\ca(-\one)a_j\ca(-\one)a_k\ca(\one)a_l\ca(\one)$. 
When the multiset is either  $-+++$ or  $---+$, there are only four distinct quartic operators in each case,
and after normal ordering has been applied, we get
\be\label{pattern2}
i\ca j\ca k\ca l\ca  \qquad  j\ca i\ca l\ca k\ca  \qquad k\ca l\ca i\ca j\ca  \qquad  l\ca k\ca j\ca i\ca  .
\ee

When $k=1=r_1+r_2+r_3+r_4$ we can see that the only quartic operators which may have a non-zero contribution have multiset of the form $\Big\{-\one,\one,\one,\one\Big\}$, also abbreviated by $-+++$, as there is no chance for ambiguity here. If any $r_n=\tre$, there must be at least one more positive mode in order for the sum to equal 1. It is clear from the basis of $\hV^1_{3/2}$ that such quartic operators of this form annihilate $v$.  If the largest mode is $\one$, then the only chance for the sum to equal 1 is the multiset presented above. 
\end{ex}

\chapter{Coset Virasoro Representations}\label{cha:cosetvir}

\section{Neveu-Schwarz Representations }\label{sec:NS}

In this chapter we investigate the decomposition of the Neveu-Schwarz representation,  $\hat V^0\oplus \hat V^1$ of the affine Kac-Moody Lie algebra $D_4^{(1)}$ with respect to its subalgebra $G_2^{(1)}$. Using the coset Virasoro construction \cite{GKO} from the previous section, we provide highest weight vectors for $G_2^{(1)}$ in $\hat V^0\oplus \hat V^1$  that will be used to decompose $\hat V^0\oplus \hat V^1$ as the tensor products of the irreducible Virasoro modules $L(1/2,h^{1/2})$ and $L(7/10,h^{7/10})$, where $h^{1/2}\in\{0,1/2\}$ and $h^{7/10}\in\{0,1/10,3/5,3/2, \}$. Specifically we provide the details that show that
\begin{align*} \hat V^0&\supseteq L(1/2,0)\otimes L(7/10,0)\otimes W(\om_0) \\
&\oplus L(1/2,0)\otimes L(7/10,3/5)\otimes W(\om_2) \\
&\oplus L(1/2,1/2)\otimes L(7/10,1/10)\otimes W(\om_2)\\
&\oplus L(1/2,1/2)\otimes L(7/10,3/2)\otimes W(\om_0)
\end{align*}
and that
\begin{align*} \hat V^1&\supseteq L(1/2,1/2)\otimes L(7/10,0)\otimes W(\om_0) \\
&\oplus L(1/2,1/2)\otimes L(7/10,3/5)\otimes W(\om_2)\\
&\oplus L(1/2,0)\otimes L(7/10,1/10)\otimes W(\om_2)\\
&\oplus L(1/2,0)\otimes L(7/10,3/2)\otimes W(\om_0) 
\end{align*}
and, in fact, the containments above are equalities, and combining them we have
\begin{align*}
 \hat V^0\oplus \hat V^1 &= \Big(L(1/2,0)\oplus L(1/2,1/2)\Big) \otimes \Big(L(7/10,0)\oplus L(7/10,3/2)\Big) \otimes W(\om_0) \\
 &\oplus \Big(L(1/2,0)\oplus L(1/2,1/2)\Big) \otimes \Big(L(7/10,1/10)\oplus L(7/10,3/5)\Big) \otimes W(\om_2).
\end{align*}

Note that $L(1/2,0)\oplus L(1/2,1/2) = \CM_{e}(\hz)$ has the structure of a vertex operator superalgebra as well as being a one-fermion
Neveu-Schwarz Clifford module, where the Clifford generators are of the form $e(m) = a_4(m) + a^*_4(m)$. We believe that
$L(7/10,0)\oplus L(7/10,3/2)$ also has the structure of a vertex operator superalgebra, but the primary vector generating its super part
is of weight $3/2$, so it is not a Clifford module. The sum $L(7/10,1/10)\oplus L(7/10,3/5)$ also appears in the decomposition below, and
we believe it is a module for the superalgebra $L(7/10,0)\oplus L(7/10,3/2)$. These superalgebras may already have appeared in the
literature and may have standard names. 

As references for the construction of vertex operator superalgebras and their modules we mention \cite{FFR} and \cite{Kac}.
Later we will prove this equality by using the principal characters and graded dimensions of each of these representations.

We introduce the notation
\be\label{hwvnotation} \pu(\hV^i_n, h^{1/2},h^{7/10},\om_j)\ee
for a vector in $\hV^i_n$ which is a highest weight vector with respect to both coset Virasoro operators and the $G_2^{(1)}$ subalgebra, 
which has $L_0^{1/2}$ eigenvalue $h^{1/2}$, $L_0^{7/10}$ eigenvalue 
$h^{7/10}$, and $G_2^{(1)}$ weight $\om_j$. We use the following abbreviations as in \eqref{depth2notations}. 
For Clifford operators we write
$i\ca(r):=a\ca_i(r)$ for $1\leq i\leq 4$, $r\in Z$, and for vectors we write
\be\label{abbrev4}i\ca j\ca k\ca l\ca := a\ca_i(-\one)a\ca_j(-\one)a\ca_k(-\one)a\ca_l(-\one){\bf1} \in\hV^0_2\ee 
and 
\be\label{abbrev3}i\ca j\ca k\ca := a\ca_i(-\one)a\ca_j(-\one)a\ca_k(-\one){\bf1} \in\hV^0_{3/2}.\ee
The $D_4$ weight of such vectors is easily found since each occurrence 
of a Clifford generator $a_i(r)$ adds weight $\ep_i$, while $a^*_i(r)$ subtracts weight $\ep_i$, and the $D_4$ weight of $\bf1$ is $0$, corresponding to the fact that $\hV^0_0$ is the $1$-dimensional trivial $D_4$ module.
In the Ramond modules, the $D_4$ weight of ${\bf1}'$ is $\lambda_4$, corresponding to the fact that $\hV^2_{1/2}$ is the $8$-dimensional 
even spinor $D_4$-module. We also use the fact that the top graded piece of the $G_2^{(1)}$-module $W(\om_0)$ is the $1$-dimensional trivial $G_2$ module, and the top graded piece of $W(\om_2)$ is the $7$-dimensional $G_2$ module. This means that if a highest weight vector 
\eqref{hwvnotation} has $\om_j$ then the $D_4$ weight of the vector must project under \eqref{D4toG2projectionepsilons} to the $G_2$ weight
$0$ if $j=0$ or to $\bar\lambda_2$ if $j=2$. Knowing a basis for $\hV^i_n$, $0\leq n\leq 2$, we can list those basis vectors whose weights have the correct projection, and express the desired highest weight vector as a linear combination of them. We can find the conditions on the coefficients of the combination that correspond to it being annihilated by $L_1^{1/2}$, $L_2^{1/2}$, $L_1^{7/10}$, $L_2^{7/10}$, the positive simple root vectors of $G_2^{(1)}$, and so that $L_0^{1/2}$ and $L_0^{7/10}$ act with the eigenvalues $h^{1/2}$ and $h^{7/10}$, respectively.
In the following Lemma we begin this process by giving the appropriate lists of basis vectors for each highest weight vector that was found and
contributed to the answer. 

\begin{lem}\label{lem:hwvs}  In the Neveu-Schwarz modules, for highest weight vectors $ \pu(\hV^i_n, h^{1/2},h^{7/10},\om_j)$ in $\hV^0_n$ or in 
$\hV^1_n$ for $0\leq n\leq 2$, we can express each vector as some linear combination of specific basis vectors as follows:

\begin{enumerate}[\rm (a)]
\item $\hV^0_0$ is $1$-dimensional with basis $\{{\bf1}\}$, which has $G_2$ weight $0$.
\item In $\hV^1_{1/2}$ any vector of $G_2$ weight $0$ is a linear combination of $\left\{ 4(-\one){\bf1}, 4^*(-\one){\bf1}\right\}$, and any vector of $G_2$ weight $\bar\lambda_2$ is a linear combination of $\left\{1(-\one){\bf1}\right\}$. 
\item In $\hV^0_1$ there are no highest weight vectors with $G_2$ weight 0, and any vector of $G_2$ weight $\bar\lambda_2$ is a linear combination of $\left\{ 1(-\one)4(-\one){\bf1}, 1(-\one)4^*(-\one){\bf1}, 2(-\one)3(-\one){\bf1}\right\}$. 
\item In $\hV^1_{3/2}$ any vector of $G_2$ weight $0$ is a linear combination of \newline
$\left\{4(-\tre){\bf1}, 4^*(-\tre){\bf1}, 11^*4, 11^*4^*, 22^*4, 22^*4^*, 33^*4,33^*4^*, 1^*23, 12^*3^*\right\}$,  
and any vector of $G_2$ weight $\bar\lambda_2$ is a linear combination of $\left\{1(-\tre){\bf1}, 122^*,133^*,144^*, 234, 234^*\right\}$.  
\item In $\hV^0_2$ any vector of $G_2$ weight $0$ is a linear combination of \newline
$\left\{11^*44^*, 22^*44^*, 33^*44^*, 1^*234, 1^*234^*, 12^*3^*4, 12^*3^*4^*, 11^*22^*, 11^*33^*, 22^*33^*\right\}$, \newline and there are no highest weight vectors with $G_2$ weight $\bar\lambda_2$.
\end{enumerate}
\end{lem}

\begin{thm}\label{thm:hwvs} The following eight vectors are annihilated by the operators \eqref{usebeta1} - \eqref{usetheta} which represent the positive simple root vectors of $G_2^{(1)}$, by the $c=\one$ coset Virasoro positive operators $L_k^{1/2}$ for $k=1,2$, and by the $c = \frac{7}{10}$ coset Virasoro positive operators $L_k^{7/10}$ for $k=1,2$.

\begin{enumerate}[\rm(a)]
\item $\pu(\hat V^0_0 ,0 ,0 ,\om_0 ) = {\bf1}$
\item $ \pu(\hat V^0_1 ,\frac{1}{2} , \frac{1}{10},\om_2 ) = 1(-\one)4(-\one){\bf1} +1(-\one)4^*(-\one){\bf1}$
\item  $\pu(\hat V^0_1 ,0 ,\frac{3}{5} ,\om_2 ) = 2\Big(2(-\one)3(-\one){\bf1}\Big) - 1(-\one)4(-\one){\bf1} + 1(-\one)4^*(-\one){\bf1}$
\item $ \pu(\hat V^0_2 ,\frac{1}{2} , \frac{3}{2},\om_0 )  = 11^*44^* - 22^*44^* - 33^*44^* + 1^*234 + 1^*234^* + 12^*3^*4 + 12^*3^*4^* $
\item $\pu(\hat V^1_{1/2} ,\frac{1}{2} ,0 ,\om_0 )  = 4(-\one){\bf1}+4^*(-\one){\bf1}$
\item $\pu(\hat V^1_{1/2} , 0,\frac{1}{10} ,\om_2 )  =  1(-\one){\bf1}$
\item $\pu(\hat V^1_{3/2} ,\frac{1}{2} ,\frac{3}{5} ,\om_2 )  =  234+234^*-144^* $
\item $\pu(\hat V^1_{3/2} ,0 , \frac{3}{2},\om_0 ) = 11^*4-11^*4^* - 22^*4+22^*4^* - 33^*4+33^*4^* + 2\Big(1^*23+12^*3^*\Big) $
\end{enumerate}

\end{thm}

\begin{prf} We will only show the proofs for the easiest part, (a), and for part (g), which exhibits all of the techniques needed to complete the other parts, using the appropriate tables from Appendix \ref{sec:NSunused}.

(a) This is the easiest vector to check. Since ${\bf1}$ is at depth $0$, it is annihilated by the operators \eqref{usebeta1} - \eqref{usetheta} which represent the positive simple root vectors  in $G_2^{(1)}$ because each normal ordered component of each operator contains a Clifford element which annihilates ${\bf1}$.  In particular,  $X_{\beta_1}(0) {\bf1} = 0 =X_{\beta_2}(0) {\bf1} = X_{-\theta}(1) {\bf1} $.  The same is true for the Virasoro operators with positive subscripts found in \eqref{LkDBNS} and \eqref{LkBGNS}. In particular,   $L_1^{1/2}{\bf1} = 0 = L_1^{7/10}{\bf1}$ and $L_2^{1/2}{\bf1} = 0 = L_2^{7/10}{\bf1} $. 

 (g) By Lemma \ref{lem:hwvs}(d), any highest weight vector $\pu(\hat V^1_{3/2} ,\frac{1}{2} ,\frac{3}{5} ,\om_2 )$ must be a linear combination of the form $v=\sum_{i=1}^6 b_iv_i$, where 
$$ v_1 = 1(-\tre){\bf1},\ \ v_2 = 122^*,\ \  v_3 = 133^*,\ \  v_4 =144^*,\ \  v_5 = 234,\ \  v_6 =234^*.$$
From its depth, any such combination is clearly annihilated by the Virasoro operators $L_k^{1/2}$ and $L_k^{7/10}$ for $k\geq 2$. In order for $v$ to be a HWV, it must be annihilated by \eqref{usebeta1} - \eqref{LkBGNS}. First we  apply the  $G_2^{(1)}$  operators \eqref{usebeta1} - \eqref{usetheta}, keeping in mind Remark \ref{rmk:G2ops} which in this situation implies that $r\in \{\pm\one, \pm\tre\}$.
Tables \ref{tab:beta1} - \ref{tab:theta} show that 
\begin{align}\label{usedG2}
X_{\beta_1}(0) v &= b_3(123^*) -b_2(123^*)  \nonumber\\
X_{\beta_2}(0) v&= b_5(134) + b_6(134^*) + b_3(-134^*) + b_4(134^*) - b_3(-134) - b_4(-134)\nonumber \\
X_{-\theta}(1) v& = b_1\Big(2^*(-\one){\bf1}\Big)+b_2\Big(2^*(-\one){\bf1}\Big)
\end{align}
so these are all $0$ when the coefficients $b_i$ satisfy the linear equations, 
\begin{align}\label{maple1}
0&= -b_2 +b_3  \nonumber\\
0&= b_3+b_4+b_5  \nonumber\\
0&= -b_3 + b_4 + b_6  \nonumber\\
0& = b_1+b_2 .
\end{align}

Next we apply the Virasoro operators  \eqref{LkDBNS} and \eqref{LkBGNS} with $k=1$.  
For the action of the quadratic operators on $v$, contrary to the case above, we do not need all  $r\in \{\pm\one, \pm\tre\}$; we need only consider $r=\tre$.  The justification is as follows; if $r>\tre$ the Clifford generators with positive modes will annihilate $v$. If $r=\one$ then $\fno i(\one)i^*(\one)\fno + \fno i^*(\one)i(\one)\fno $ as well as   $\fno 4(\one)4(\one)\fno $ and $\fno 4^*(\one)4^*(\one)\fno $ are all the zero-operator by the anti-symmetry of normal ordering.  If $r=-\one$ then the coefficient $(r+\one)=0$.  (See \eqref{normderk}). Finally if $ r<-\one$ then $1-r>\tre$ and the positive modes will annihilate $v$. We do not need a table for this calculation as the explanation above tells us that $L_1^{1/2}\cdot v=0$, and also that the contribution of the quadratic operators of $L_1^{7/10}$ on $v$ come only from the terms $i=1$ on $b_1v_1$,
\bea\label{L1quadterms} &-&\tfrac{1}{10}(\tre+\one)\Big(\fno 1^*(\tre)1(-\one)\fno+\fno 1(\tre)1^*(-\one)\fno\Big)  b_1 1(-\tre){\bf1} \notag \\
&=& -\tfrac{1}{10}(2)\Big(- 1(-\one)1^*(\tre) - 1^*(-\one) 1(\tre)\Big) b_1 1(-\tre){\bf1} \notag\\
&=& \tfrac{1}{5} b_1 1(-\one)1^*(\tre) 1(-\tre){\bf1} \notag\\
&=& \tfrac{1}{5} b_1 1(-\one){\bf1}.\eea
Table \ref{tab:L1} summarizes the actions of the quartic operators on $v$ based on the observations in Remark \ref{rmk:vops} and Example \ref{ex:depth3}. Suppressed in each entry is the sum over all possible multisets  \eqref{pattern2} applied to each $v_i$.  The table and the discussion above show that 
\begin{align}\label{usedL1}
L_1^{1/2} v&=0\nonumber\\
L_1^{7/10} v&= \Big(\tfrac{1}{5}b_1 -\tfrac{1}{5}b_2 - \tfrac{1}{5}b_3 + \tfrac{1}{5}b_5 - \tfrac{1}{5}b_6\Big) \Big(1(-\one){\bf1}\Big) . 
\end{align}
In order for $L_1^{7/10}\cdot v = 0$, the coefficients must satisfy
\begin{align}\label{maple2}
0&= b_1 -  b_2 -  b_3 + b_5 - b_6 . 
\end{align}
The two systems \eqref{maple1} and \eqref{maple2} together imply that 
\begin{align}\label{maple3} 
b_1& =  b_2 = b_3 =0\notag\\
b_4 &=  -b_6 \notag\\
b_5 &=  b_6 \notag\\
b_6 &= s \qquad \mbox{is free.}  
\end{align}
Thus the vector proposed in the theorem is the unique (up to scalar multiples) linear combination of basis vectors which is highest with respect to $G_2^{(1)}$ as well as the $c=\one$ coset Virasoro and the $c = \frac{7}{10}$ coset Virasoro.   
 \qed
\end{prf}

\begin{thm}\label{thm:eigen}
The eight vectors from Theorem \ref{thm:hwvs} have the appropriate $L^c_0$ eigenvalues for the $c=\frac{1}{2}$ and $c=\frac{7}{10}$ coset Virasoro representations, as found in their labels. 
\end{thm}

\begin{prf} Again we only show parts (a) and (g). 

(a) It is easy to see that both   $L_0^{1/2}{\bf1} = 0 $ and $L_0^{7/10}{\bf1} = 0 $ since each normal ordered component of each operator contains a Clifford element with a positive mode number, which annihilates ${\bf1}$.  

 (g) Let $v =  \pu(\hat V^1_{3/2} ,\frac{1}{2} ,\frac{3}{5} ,\om_2 )  =  234+234^*-144^* $. We claim $L_0^{1/2} v=\one v$ and $L_1^{7/10} v=\tfrac{3}{5} v$. 
When $L_0^{1/2}$ or $L_0^{7/10}$ acts on $v$, the possible non-zero contributions of the action of the $L^c_0$ operator come  from linear combinations of quartic operators whose actions are computed in Tables \ref{tab:--++} - \ref{tab:---+} and summarized in Table \ref{tab:sum2}, as well as
from linear combinations of quadratic operators whose actions are computed in Tables \ref{tab:quadlam2a} and \ref{tab:quadlam2b} and summarized in Table \ref{tab:quadlam2c}; all according to Remark \ref{rmk:vops} and Example \ref{ex:depth3}. We then have
\begin{align}\label{LeigenDB} L_0^{1/2} v &=  L_0^{1/2}( 234+234^*-144^*)\nonumber\\
&=\tfrac{1}{4}(2)\Big(234^*+234+234+234^*-144^*-144^*\Big)\nonumber\\
&=\tfrac{1}{2}\Big(234+234^*-144^*\Big)\nonumber\\
&=\one v \end{align}
as claimed. Now we compute
\begin{align}\label{LeigenBG} L_0^{7/10} v &=  L_0^{7/10}( 234+234^*-144^*)\nonumber\\
&= -\tfrac{1}{10}\Big(v_4+v_5+v_6+v_5+v_6+v_4+v_6+v_4+v_5 \Big) \nonumber\\
&\phm\ph +\tfrac{1}{20}\Big(v_4+v_6+v_4+v_5+v_5+v_6 \Big) \nonumber\\
&\phm\ph -\tfrac{1}{5}\Big( -v_5-v_6-v_5+v_2+v_3+v_4-v_1-v_6-v_2-v_3+v_4+v_1\Big) \nonumber\\
&= \tfrac{12}{20}(-v_4+v_5+v_6) \nonumber\\
&=\tfrac{3}{5}v
\end{align}
verifying the second part of the claim.
 \qed
\end{prf}

\section{Ramond Representations }\label{sec:R}

Now we investigate the decomposition of the Ramond representation $\hat V^2\oplus \hat V^3$ using the tensor products of the irreducible Virasoro modules $L(1/2,h_3)$ and $L(7/10,h_4)$, where $h_3\in\{0,1/2,1/16\}$ and $h_4\in\{7/16, 3/80 \}$. We show that
\begin{align*} \hat V^2&\supseteq L(1/2,1/16)\otimes L(7/10,3/80)\otimes W(\om_2) \\
&\oplus L(1/2,1/16)\otimes L(7/10,7/16)\otimes W(\om_0) 
\end{align*}
and with the same decomposition, 
\begin{align*} \hat V^3&\supseteq L(1/2,1/16)\otimes L(7/10,3/80)\otimes W(\om_2) \\
&\oplus L(1/2,1/16)\otimes L(7/10,7/16)\otimes W(\om_0). 
\end{align*}
In this section we continue to use the notation \eqref{hwvnotation}, and we have analogous theorems and lemmas.

\begin{lem}\label{lem:hwvsR}  In the Ramond modules, for highest weight vectors $ \pu(\hV^i_n, h^{1/2},h^{7/10},\om_j)$ in $\hV^2_{1/2}$ or in 
$\hV^3_{1/2}$, we can express each vector as some linear combination of specific basis vectors as follows:  
\begin{enumerate}[\rm (a)]
\item In $\hV^2_{1/2}$ any vector of $G_2$ weight $0$ is a linear combination of $\left\{ 1^*(0)4^*(0){\bf1'},  2^*(0)3^*(0){\bf1'} \right\}$, 
and any vector of $G_2$ weight $\bar\lambda_2$ is a linear combination of $\left\{{\bf1'}\right\}$. 
\item In $\hV^3_{1/2}$ any vector of $G_2$ weight $0$ is a linear combination of $\left\{ 1^*(0){\bf1'} , 2^*(0)3^*(0)4^*(0) {\bf1'}\right\}$, 
and any vector of $G_2$ weight $\bar\lambda_2$ is a linear combination of $\left\{4^*(0){\bf1'}\right\}$. 
\end{enumerate}
\end{lem}

\begin{thm}\label{thm:hwvsR} The following four vectors are annihilated by the operators \eqref{usebeta1} - \eqref{usetheta} which represent the positive simple root vectors of $G_2^{(1)}$, by the $c=\one$ coset Virasoro positive operators $L_k^{1/2}$ for $k=1,2$, and by the $c = \frac{7}{10}$ coset Virasoro positive operators $L_k^{7/10}$ for $k=1,2$.
\begin{enumerate}[\rm(a)]
\item $\pu(\hat V^2_{1/2} ,\frac{1}{16} ,\frac{3}{80},\om_2 ) = {\bf1'}$
\item $\pu(\hat V^2_{1/2} ,\frac{1}{16} ,\frac{7}{16},\om_0 ) =1^*(0)4^*(0){\bf1'} - 2^*(0)3^*(0) {\bf1'}$
\item $\pu(\hat V^3_{1/2} ,\frac{1}{16} ,\frac{3}{80},\om_2 ) = 4^*(0){\bf1'}$
\item $\pu(\hat V^3_{1/2} ,\frac{1}{16} ,\frac{7}{16},\om_0 ) =1^*(0){\bf1'} + 2^*(0)3^*(0)4^*(0) {\bf1'}$
\end{enumerate}

\end{thm}

\begin{prf} Each of the vectors from parts (a)-(d) are at depth $\one$, therefore they are annihilated by $L_k^{1/2}$ and $L_k^{7/10}$ for $k\in\{1,2\}$.  Recall Remark \ref{rmk:G2ops} which says that we need only consider $r=0$ when applying the $G_2^{(1)}$ operators. When $r=0$,  $X_{-\theta}(1) = 2^*(0)1^*(1)$, whose positive mode number will annihilate the vectors from parts (a)-(d). When applying $X_{\beta_1}(0)$ it is easy to see that it also acts trivially on each vector (a)-(d). Either  the Clifford element  $2(0)$ anti-commutes to annihilate ${\bf1'}$, or $3^*(0)$ anti-commutes to square (and hence annihilate) another factor of $3^*(0)$. Thus we need only apply $X_{\beta_2}(0)$ with care. Those calculations are summarized in Table \ref{tab:beta2R} located in Appendix \ref{sec:Ramond}. Since $X_{\beta_2}(0)$ also annihilates each vector, this completes the proof.
 \qed
\end{prf}

The following lemmas are proved by simply invoking Definition \ref{dfn:generalfno}. 
\begin{lem}\label{lem:4zeros1} For $i\neq j, \fno i(0)i^*(0)j(0)j^*(0)\fno = \frac{1}{4}I -i^*(0)j^*(0)i(0)j(0)- \one i^*(0)i(0)-\one j^*(0)j(0)$. \end{lem}
\begin{lem}\label{lem:4zeros2} For $i,j,m,n$ all distinct, $\fno i\ca(0)j\ca(0)m\ca(0)n\ca(0) \fno = i\ca(0)j\ca(0)m\ca(0)n\ca(0)$. \end{lem}

\begin{thm}\label{thm:eigenR}
The four vectors from Theorem \ref{thm:hwvsR} have the appropriate $L^c_0$ eigenvalues for the $c=\frac{1}{2}$ and $c=\frac{7}{10}$ coset Virasoro representations, as found in their labels. 
\end{thm}

\begin{prf}
When applying the quadratic operators to each of the vectors from parts (a)-(d) we notice first that we can only have non-zero contributions when $r=0$; otherwise the presence of a positive mode number will annihilate ${\bf1'}$. By the anti-symmetry of normal ordering, $\fno i^*(0)i(0)\fno + \fno i(0)i^*(0) \fno$, $\fno 4(0)4(0)\fno$ and 
$\fno 4^*(0)4^*(0)\fno$ each equal the zero-operator. Hence the action of $L_0^{1/2}$ given in \eqref{LkDBR} reduces to $\tfrac{1}{16}I$ on all four vectors, verifying that $h^{1/2} = \tfrac{1}{16}$ for each vector. 

The action of $L_0^{7/10}$ on the  vectors from parts (a)-(d) can be discerned from the observation made above about the trivial action of all of the  quadratic operators, and from Table \ref{tab:L1R}.  Inspection of the entries of the table leads us to infer that the eigenvalues will be exactly the same for  $\pu(\hat V^2_{1/2} ,\frac{1}{16} ,\frac{3}{80},\om_2 )$ and  $\pu(\hat V^3_{1/2} ,\frac{1}{16} ,\frac{3}{80},\om_2 )$, as well as for $\pu(\hat V^2_{1/2} ,\frac{1}{16} ,\frac{7}{16},\om_0 )$ and $\pu(\hat V^3_{1/2} ,\frac{1}{16} ,\frac{7}{16},\om_0 )$. Since the proofs are so similar, we only present the calculation for part (a).

By \eqref{LkBGR}, Lemmas \ref{lem:4zeros1} and \ref{lem:4zeros2} and Table \ref{tab:L1R}, we have
\begin{align*}
L^{7/10}_0{\bf1'} = \tfrac{7}{80}{\bf1'} -\tfrac{1}{5}\Big(&\fno 1(0)1^*(0)2(0)2^*(0)\fno + \fno 1(0)1^*(0)3(0)3^*(0)\fno - \fno 2(0)2^*(0)3(0)3^*(0)\fno \\
+&\fno 1^*(0)2(0)3(0)4(0)\fno + \fno 1(0)2^*(0)3^*(0)4^*(0)\fno \\
-&\fno 1^*(0)2(0)3(0)4^*(0)\fno -\fno 1(0)2^*(0)3^*(0)4(0)\fno \Big){\bf1'}.
\end{align*}
This yields, $L^{7/10}_0 {\bf1'} = \tfrac{7}{80}{\bf1'} -\tfrac{1}{5}(\for {\bf1'}+ \for {\bf1'} - \for {\bf1'}) =  \tfrac{3}{80}{\bf1'}$, which verifies $h^{7/10} =\tfrac{3}{80}.$ 
 \qed
\end{prf}

\chapter{Conclusions}

The results of Chapter \ref{cha:cosetvir} show that $\hV$ decomposes into at least twelve highest weight modules
for $\Vir^{1/2} \otimes\Vir^{7/10}\otimes G_2^{(1)}$, each of the form 
$L(1/2,h^{1/2})\otimes L(7/10,h^{7/10})\otimes W(\om_j)$, where $h^{1/2}\in\{0,1/2,1/16\}$,  $h^{7/10}\in\{0,1/10,3/5,3/2,7/16,3/80 \}$, and $j\in\{0,2\}$. 
We found eight of them inside the Neveu-Schwarz module, and four more inside of the Ramond module. 
The content of Chapter \ref{cha:punchline} is that $\hV$ contains no more than these twelve highest weight modules, as the character of $\hV$ equals the character of the sum of those twelve modules. Hence we have proved our main theorem: 
\begin{thm}\label{thm:conclusion} The direct sum of the four level-1 irreducible highest weight modules for the affine Kac-Moody Lie algebra $D_4^{(1)}$, $\hV = \hV^0 \oplus \hV^1 \oplus\hV^1 \oplus\hV^3$, decomposes with respect to its affine subalgebra $G_2^{(1)}$ into the direct
sum of twelve $\Vir^{1/2} \otimes\Vir^{7/10}\otimes G_2^{(1)}$-modules as follows:
 \begin{align*}
\hV &= \Big(L(1/2,0)\oplus L(1/2,1/2)\Big) \otimes \Big(L(7/10,0)\oplus L(7/10,3/2)\Big) \otimes W(\om_0) \\
&\oplus \Big(L(1/2,0)\oplus L(1/2,1/2)\Big) \otimes \Big(L(7/10,1/10)\oplus L(7/10,3/5)\Big) \otimes W(\om_2)\\
&\oplus \Big(  L(1/2,1/16)\otimes L(7/10,3/80)\oplus  L(1/2,1/16)\otimes L(7/10,3/80) \Big)\otimes W(\om_2)\\
&\oplus\Big( L(1/2,1/16)\otimes L(7/10,7/16)\oplus L(1/2,1/16)\otimes L(7/10,7/16)\Big) \otimes W(\om_0).
\end{align*}
Each summand is determined by a vector $\pu(\hV^i_n, h^{1/2},h^{7/10},\om_j)$ which is a highest weight vector with respect to both coset Virasoro operators and the $G_2^{(1)}$ subalgebra, which has $L_0^{1/2}$ eigenvalue $h^{1/2}$, $L_0^{7/10}$ eigenvalue $h^{7/10}$, and $G_2^{(1)}$ weight $\om_j$. The explicit form of these highest weight vectors is given in Theorems \ref{thm:hwvs} and \ref{thm:hwvsR}, showing which summands occur in each $D_4^{(1)}$-module $\hV^i$. 
\end{thm}

\begin{rmk}\label{rmk:conclusion1}
During the process of finding and proving the decomposition above, we also learned more about the structure of $\hV$ and the branching rules for decomposing $D_4$-modules into $B_3$-modules and $G_2$-modules.  For instance, at depth $ \one$ the triality automorphism $\sigma$ cyclically permutes the 8-dimensional natural, even, and odd spinor representations of the finite dimensional Lie algebra $D_4$, and $\shat$ permutes the highest weight vectors found in $\hV^i_{1/2}$, $1\leq i\leq 3$, 
\begin{align}
\shat^2 \Big(\pu(\hat V^1_{1/2} , 0,\tfrac{1}{10} ,\om_2 )\Big) = \shat \Big(\pu(\hat V^2_{1/2} , \tfrac{1}{16},\tfrac{3}{80} ,\om_2 )\Big) = \pu(\hat V^3_{1/2} , \tfrac{1}{16},\tfrac{3}{80} ,\om_2 )\notag \\
\shat^2 \Big(\pu(\hat V^1_{1/2} ,\tfrac{1}{2} ,0 ,\om_0 )\Big) = \shat \Big(\pu(\hat V^2_{1/2} , \tfrac{1}{16},\tfrac{7}{16} ,\om_0 )\Big) = \pu(\hat V^3_{1/2} , \tfrac{1}{16},\tfrac{7}{16} ,\om_0 ).
\end{align}

This might lead one to naively believe that the two highest weight vectors in $\hV^1_{3/2}$ would be sent by $\shat$ and $\shat^2$ to highest weight vectors in $\hV^2_{3/2}$ and $\hV^3_{3/2}$. This is just not true, as shown by the character theory supplied in Chapter \ref{cha:punchline}. 
The point is that the vertex operators representing the coset Virasoro algebras transform under conjugation by $\shat$ into different coset Virasoro operators which are associated with different $B_3$ subalgebras inside $D_4$. While we did not need these ``twisted" coset Virasoro operators to achieve our goal, we think they are interesting and will discuss them below. 
\end{rmk}

\begin{rmk}\label{rmk:conclusion2}
Throughout, we have mentioned that as we focus our attention from $D_4$ to $B_3$, we see that there are three isomorphic copies to choose from, the fixed points under $\tau, \sigma\tau$, or $\sigma^2\tau$. Even though our exposition mainly discusses the role of the former, $B_3 \cong D_4^{\tau}$, its affinization, and its role in the coset Virasoro constructions, one can (and arguably should) investigate the behavior of its $\sigma$-``twisted sisters",  and their coset constructions. Fortunately, we are able to provide the following information: 
\begin{align*}
\fb{$44^*$} &= 44^*+4^*4+44+4^*4^* \notag\\
\fb{$11^*22^*$} &= 11^*22^*+11^*33^*-22^*33^*\notag \\
\fb{$22^*44^*$} &= 22^*44^*+33^*44^*-11^*44^* \notag\\
\fb{$1^*234$} &= 1^*234 + 12^*3^*4^* = - \shat^2 (44 + 4^*4^*)\notag\\
\fb{$1^*234^*$} &= 12^*3^*4 + 1^*234^*= \shat (44 + 4^*4^*) .\notag
\end{align*}
which  allows us to write
\begin{align*}
\shat \fb{$44^*$} &= \one\D -\one \fb{$11^*22^*$}-\one \fb{$22^*44^*$}  + \fb{$1^*234^*$} \\
\shat \fb{$11^*22^*$} &=  \one\D - (44^*+4^*4) +\one\fb{$11^*22^*$} -\one\fb{$22^*44^*$} \\
\shat \fb{$22^*44^*$} &= -\one\D + (44^*+4^*4) +\one\fb{$11^*22^*$} -\one\fb{$22^*44^*$} \\
\shat \fb{$1^*234$} &=  -( 44 + 4^*4^*)\\
\shat \fb{$1^*234^*$} &=- \fb{$1^*234$} 
\end{align*}
and
\begin{align*}\shat^2 \fb{$44^*$} &= \one\D -\one \fb{$11^*22^*$}+\one \fb{$22^*44^*$}  - \fb{$1^*234$} \\
\shat^2 \fb{$11^*22^*$} &=  \one\D - (44^*+4^*4) +\one\fb{$11^*22^*$} +\one\fb{$22^*44^*$} \\
\shat^2 \fb{$22^*44^*$} &= \one\D - (44^*+4^*4) -\one\fb{$11^*22^*$} -\one\fb{$22^*44^*$} \\
\shat^2 \fb{$1^*234$} &= - \fb{$1^*234^*$}\\
\shat^2 \fb{$1^*234^*$} &= (44 + 4^*4^*)
\end{align*}
as well as
\begin{align*}
\DB
&= \tfrac{1}{4} \framebox{$44^*$}\\
\shat\DB
&= \tfrac{1}{8}\omega_{D_4}-\tfrac{1}{8}\fb{$11^*22^*$} - \tfrac{1}{8}\fb{$22^*44^*$} + \tfrac{1}{4}\fb{$1^*234^*$} \\
\shat^2\DB
&= \tfrac{1}{8}\omega_{D_4}-\tfrac{1}{8}\fb{$11^*22^*$} +\tfrac{1}{8}\fb{$22^*44^*$} - \tfrac{1}{4}\fb{$1^*234$}
\end{align*}
and 
\begin{align*}
\BG
&= \tfrac{1}{5}\D - \tfrac{1}{20}\fb{$44^*$} - \tfrac{1}{5}\fb{$11^*22^*$} - \tfrac{1}{5}\fb{$1^*234$} + \tfrac{1}{5}\fb{$1^*234^*$}\\
\shat\BG
& = \tfrac{3}{40}\D+ \tfrac{1}{5} \fb{$44^*$}   -\tfrac{3}{40} \fb{$11^*22^*$} +\tfrac{1}{8}\fb{$22^*44^*$}   - \tfrac{1}{5}\fb{$1^*234$}- \tfrac{1}{20}\fb{$1^*234^*$} \\
\shat^2\BG
&= \tfrac{3}{40}\D+ \tfrac{1}{5} \fb{$44^*$}   -\tfrac{3}{40} \fb{$11^*22^*$} -\tfrac{1}{8}\fb{$22^*44^*$}   + \tfrac{1}{20}\fb{$1^*234$} + \tfrac{1}{5}\fb{$1^*234^*$}. \\
\end{align*}

It is shown in \cite{FFR} that for any automorphism $g$ in the triality group
$G = \langle \sigma,\tau\rangle \cong S_3$, for any $v\in\hV$ we have 
\be
gY(v,z)g\inv = Y(gv,z),
\ee
and, in particular, for any $u, v\in \hV$, we have
\be
gY(v,z) u = Y(gv,z)gu \label{ygvzgu}
\ee
which means that the coefficients satisfy $gY_m(v)u = Y_m(gv)gu$ for all $m\in Z$. For $1\leq i \leq 3$, we have the vertex operators
\bea\label{twistedcosetVirs}
\shat^iY\left(\DB,z\right)\shat^{-i} &=&Y(\shat^i(\DB),z) = \sum_{k\in\Z} L_k^{1/2,\shat^i} z^{-k-2}  \\
\shat^iY\left(\BG,z\right)\shat^{-i} &=& Y(\shat^i(\BG),z)  = \sum_{k\in\Z} L_k^{7/10,\shat^i}  z^{-k-2}
\eea
so that 
\be
L_k^{1/2,\shat^i} = \shat^iL_k^{1/2}\shat^{-i}  \qquad\hbox{and}\qquad
L_k^{7/10,\shat^i} = \shat^iL_k^{7/10}\shat^{-i}  .
\ee
It is then clear that, as $\shat$ conjugates of operators representing the Virasoro algebra, these new operators provide $\shat$-twisted 
representations of those coset Virasoro algebras. If $v\in\hV^1_n$ is a highest weight vector with respect to the untwisted coset Virasoro
operators, that is, $L_k^c v = \delta_{k,0} h^c v$ for $k\geq 0$, $c = \one$ or $c = \sev$, then 
\be
L_k^{c,\shat^i}  (\shat^i v) = \delta_{k,0} h^c (\shat^i v)
\ee
shows that $\shat^i v$ is a highest weight vector with respect to the $\shat^i$-twisted coset Virasoro operators rather than with respect to
the untwisted coset Virasoro operators. Since $\shat(\D) = \D$ and $\shat(\G) = \D$, the vertex operators representing the Virasoro algebras
associated with $D_4^{(1)}$ and $G_2^{(1)}$ are invariant under conjugation by $\shat$, so if a vector $v$ is highest with respect to either of
those two Virasoro algebras, then $\shat^i v$ will also have that property.

We can  now supply the ``reptilian" vertex operators ... scaled for their beauty:  
\begin{align}
8Y_{k}(\shat\DB)
= Y_{k}(\omega_{D_4}) -Y_{k}\fb{$11^*22^*$} - Y_{k}\fb{$22^*44^*$} + 2Y_{k}\fb{$1^*234^*$} \nonumber \\
8Y_{k}(\shat^2\DB)
=Y_{k}(\omega_{D_4})v-Y_{k}\fb{$11^*22^*$} +Y_{k}\fb{$22^*44^*$} -2Y_{k}\fb{$1^*234$}  \nonumber\\
40Y_{k}(\shat\BG)
 =3Y_{k}(\D)+8 Y_{k}\fb{$44^*$}   -3Y_{k} \fb{$11^*22^*$} +5Y_{k}\fb{$22^*44^*$}   
 \nonumber \\-8Y_{k}\fb{$1^*234$}-2Y_{k}\fb{$1^*234^*$}  \nonumber \\
40Y_{k}(\shat^2\BG)
 =3Y_{k}(\D)+8 Y_{k}\fb{$44^*$}   -3Y_{k} \fb{$11^*22^*$} -5Y_{k}\fb{$22^*44^*$}   \nonumber \\
 +2Y_{k}\fb{$1^*234$} +8Y_{k}\fb{$1^*234^*$}\ .  \label{conclusionops} 
\end{align}

We found that changing perspective on which copy of $B_3^{(1)}$ would serve as our intermediate algebra was a very enlightening step in understanding this project. 

\end{rmk}

Our plan for future work is to investigate the branching rules of other affine Kac-Moody Lie algebra representations, including \emph{twisted} algebras such as  $D_4^{(3)}$. We also wish to apply these techniques to the study of the symplectic affine Kac-Moody algebras $C_\ell^{(1)}$,  which have analogous level $-\one$ constructions from a Weyl algebra instead of a Clifford algebra. 
\chapter{Appendix}

\section{Verification that $\pu(\hat V^1 ,\frac{1}{2} ,\frac{3}{5} , \Omega_2 )$ is in fact a HWV}

\subsection{Tables Relevant to $G_2^{(1)}$}\label{sub:g2tab}

\begin{table}[h!]
$$
\begin{array}{|c|c|c|c|c|c|c|}
\hline
\color{red} r < 0&1(-\tre){\bf1}&122^*& 133^*&144^*&234&234^*\\
\hline
2(-\tre)3^*(\tre)&0&0&0&0&0&0\\
2(-\one)3^*(\one) &0&0&123^*&0&0&0\\
\hline
\hline
\color{red} r > 0&1(-\tre){\bf1}&122^*& 133^*&144^*&234&234^*\\
\hline
3^*(-\tre)2(\tre)&0&0&0&0&0&0\\
3^*(-\one)2(\one) &0&123^*&0&0&0&0\\
\hline
\end{array}
$$
\caption{Application of operators \eqref{usebeta1}}\label{tab:beta1}
\end{table}
\vskip -10pt

\begin{table}[h!]

$$
\begin{array}{|c|c|c|c|c|c|c|}
\hline
\color{red} r < 0&1(-\tre){\bf1}&122^*& 133^*&144^*&234&234^*\\
\hline
1(-\tre)2^*(\tre)&0&0&0&0&0&0\\
1(-\one)2^*(\one) &0&0&0&0&134&134^*\\
3(-\tre)4^*(\tre)&0&0&0&0&0&0\\
3(-\one)4^*(\one) &0&0&0&134^*&0&0\\
3(-\tre)4(\tre)&0&0&0&0&0&0\\
3(-\one)4(\one) &0&0&0&-134&0&0\\
\hline
\hline
\color{red} r > 0&1(-\tre){\bf1}&122^*& 133^*&144^*&234&234^*\\
\hline
2^*(-\tre)1(\tre)&0&0&0&0&0&0\\
2^*(-\one)1(\one) &0&0&0&0&0&0\\
4^*(-\tre)3(\tre) &0&0&0&0&0&0\\
4^*(-\one)3(\one) &0&0&134^*&0&0&0\\
4(-\tre)3(\tre) &0&0&0&0&0&0\\
4(-\one)3(\one) &0&0&134&0&0&0\\
\hline
\end{array}
$$
\caption{Application of operators \eqref{usebeta2}}\label{tab:beta2}
\end{table}

\begin{table}[h!]

$$
\begin{array}{|c|c|c|c|c|c|c|}
\hline
{\color{orange} net} = {\color{red} 1} \color{red} - 1 &1(-\tre){\bf1}&122^*& 133^*&144^*&234&234^*\\
\hline
\fno 2(r)3^*(-r) \fno &0&-123^*&123^*&0&0&0\\
\hline
\fno 1(r)2^*(-r) \fno &0&0&0&0&134&134^*\\
\fno 3(r)4^*(-r) \fno &0&0&-134^*&134^*&0&0\\
\fno 3(r)4(-r) \fno &0&0&-134&-134&0&0\\
\hline
\end{array}
$$
\caption{Summary of Tables \ref{tab:beta1} and \ref{tab:beta2}}\label{tab:sum1}
\end{table}

\begin{table}[h!]

$$
\begin{array}{|c|c|c|c|c|c|c|}
\hline
 \sum = 1  &1(-\tre){\bf1}&122^*& 133^*&144^*&234&234^*\\
\hline
\fno 2^*(\tre)1^*(-\one)\fno &0&0&0&0&0&0\\
\fno 2^*(\one)1^*(\one)\fno   &0&2^*(-\one){\bf1}&0&0&0&0\\
\fno 2^*(-\one)1^*(\tre)\fno  &2^*(-\one){\bf1}&0&0&0&0&0\\
\hline
\end{array}
$$
\caption{Application of operators \eqref{usetheta}}\label{tab:theta}
\end{table}

\subsection{Table Relevant to $L_1^c$}\label{sub:L1tab}

\begin{table}[h!]
$$
\begin{array}{|c|c|c|c|c|c|c|}
\hline
-+++&1(-\tre){\bf1}&122^*& 133^*&144^*&234&234^*\\
\hline
11^*22^*&0&1(-\one){\bf1}&0&0&0&0\\
11^*33^*&0&0&1(-\one){\bf1}&0&0&0\\
22^*33^* &0& 0&0&0&0&0\\
\hline
1^*234^*&0& 0&0 &0&0&0\\
12^*3^*4& 0&0&0&0&0&-1(-\one){\bf1}\\
\hline
1^*234&0&0&0&0&0&0 \\
12^*3^*4^*&0&0&0&0&-1(-\one){\bf1}&0\\
\hline
\end{array}
$$
\caption{Application of quartic operators for $L_1^c$}\label{tab:L1}
\end{table}


\subsection{Tables Relevant to $L_0^c$}\label{sub:L0tab}

\begin{table}[h!]

$$
\begin{array}{|c|c|c|c|c|c|c|}
\hline
--++&1(-\tre){\bf1}&122^*& 133^*&144^*&234&234^*\\
\hline
11^*22^*&0&0&0&0&0&0\\
11^*33^*&0&0&0&0&0&0\\
22^*33^* &0&133^*&122^*&0&234&234^*\\
\hline
1^*234^*&0& 234^*&234^* &-234^*&0&0\\
12^*3^*4& 0&0&0&0&0&122^*+133^*-144^*\\
\hline
1^*234&0&234&234&234&0&0 \\
12^*3^*4^*&0&0&0&0&122^*+133^*+144^*&0\\
\hline
\end{array}
$$
\caption{Application of quartic operators for $r_i$ multiset $--++$}\label{tab:--++}
\end{table}

\begin{table}[h!] 
$$
\begin{array}{|c|c|c|c|c|c|c|}
\hline
-+++&1(-\tre){\bf1}&122^*& 133^*&144^*&234&234^*\\
\hline
11^*22^*&0&1(-\tre){\bf1}&0&0&0&0\\
11^*33^*&0&0&1(-\tre){\bf1}&0&0&0\\
22^*33^* &0& 0&0&0&0&0\\
\hline
1^*234^*&0& 0&0 &0&0&0\\
12^*3^*4& 0&0&0&0&0&-1(-\tre){\bf1}\\
\hline
1^*234&0&0&0&0&0&0 \\
12^*3^*4^*&0&0&0&0&-1(-\tre){\bf1}&0\\
\hline
\end{array}
$$
\vskip -5pt
\caption{Application of quartic operators for $r_i$ multiset $-+++$}\label{tab:-+++}
\end{table}

\begin{table}[h!]
$$
\begin{array}{|c|c|c|c|c|c|c|}
\hline
---+&1(-\tre){\bf1}&122^*& 133^*&144^*&234&234^*\\
\hline
11^*22^*&122^*&0&0&0&0&0\\
11^*33^*&133^*&0&0&0&0&0\\
22^*33^* &0&0&0&0&0&0\\
\hline
1^*234^*&-234^*& 0&0 &0&0&0\\
12^*3^*4& 0&0&0&0&0&0\\
\hline
1^*234&-234&0&0&0&0&0 \\
12^*3^*4^*&0&0&0&0&0&0\\
\hline
\end{array}
$$
\vskip -5pt
\caption{Application of quartic operators for $r_i$ multiset $---+$}\label{tab:---+}
\end{table}

\begin{table}[h!]
$$
\begin{array}{|c|c|c|c|c|c|c|}
\hline
\color{red}r=\one &1(-\tre){\bf1}&122^*& 133^*&144^*&234&234^*\\
\hline
1^*(-\one)1(\one)&0&0&0&0&0&0\\
1(-\one)1^*(\one)&0&122^*&133^*&144^*&0&0\\
2^*(-\one)2(\one)&0&122^*&0&0&0&0\\
2(-\one)2^*(\one)&0&122^*&0&0&234&234^*\\
3^*(-\one)3(\one)&0&0&133^*&0&0&0\\
3(-\one)3^*(\one)&0&0&133^*&0&234&234^*\\
4^*(-\one)4(\one)&0&0&0&144^*&0&234^*\\
4(-\one)4^*(\one)&0&0&0&144^*&234&0\\
\hline
4(-\one)4(\one)&0&0&0&0&0&234\\
4^*(-\one)4^*(\one)&0&0&0&0&234^*&0\\
\hline
\end{array}
$$
\vskip -5pt
\caption{Application of quadratic operators for $r=\one$}\label{tab:quadlam2a}
\end{table}

\vskip -5pt

\begin{table}[h!]
$$
\begin{array}{|c|c|c|c|c|c|c|}
\hline
\color{red}r=\tre &1(-\tre){\bf1}&122^*& 133^*&144^*&234&234^*\\
\hline
1^*(-\tre)1(\tre)&0&0&0&0&0&0\\
1(-\tre)1^*(\tre)&1(-\tre)&0&0&0&0&0\\
\hline

\hline
\color{red}r=-\tre &&&&&&\\
\hline
1^*(-\tre)1(\tre)&0&0&0&0&0&0\\
1(-\tre)1^*(\tre)&1(-\tre)&0&0&0&0&0\\
\hline
\end{array}
$$
\vskip -5pt
\caption{Application of quadratic operators for $r=\pm\tre$}\label{tab:quadlam2b}
\end{table}

\begin{table}[h!]
$$
\begin{array}{|c|c|c|c|}
\hline
& 144^* & 234 & 234^*\\
\hline
11^*22^* & 0&0 &0\\
11^*33^* & 0&0 &0\\
22^*33^* & 0&234 &234^*\\
\hline
1^*234^*& -234^* &0 &0\\
12^*3^*4& 0 & 0 &122^*+133^*-144^*-1(-\tre){\bf1} \\
\hline
1^*234& 234 & 0& 0\\  
12^*3^*4^*&  0 & 122^*+133^*+144^*-1(-\tre){\bf1} &0 \\
\hline
\end{array}
$$
\vskip -5pt
\caption{Summary of Tables \ref{tab:--++} - \ref{tab:---+} (for use in Theorem \ref{thm:eigen})}\label{tab:sum2} 
\end{table}

\begin{table}[h!]
$$
\begin{array}{|c|c|c|c|}
\hline
&144^*&234&234^*\\
\hline
\fno 1(r)1^*(-r)\fno &0&0&0\\
\fno 1^*(r)1(-r)\fno &144^*&0&0\\
\fno 2(r)2^*(-r)\fno &0&0&0\\
\fno 2^*(r)2(-r)\fno &0&234&234^*\\
\fno 3(r)3^*(-r)\fno &0&0&0\\
\fno 3^*(r)3(-r)\fno &0&234&234^*\\
\fno 4(r)4^*(-r)\fno &144^*&0&234^*\\
\fno 4^*(r)4(-r)\fno &144^*&234&0\\
\hline
\fno 4(r)4(-r)\fno &0&0&234\\
\fno 4^*(r)4^*(-r)\fno &0&234^*&0\\
\hline
\end{array}
$$
\vskip -5pt
\caption{Summary of Tables \ref{tab:quadlam2a} and \ref{tab:quadlam2b}  (for use in Theorem \ref{thm:eigen})}\label{tab:quadlam2c}
\end{table}

\vskip -5pt

\section{Bases for  $\hV^1_{3/2}$  and $\hV_2^0$ }\label{sec:bases}

$$
\begin{array}{|c|rr|c|}
\hline
\hV_n^i & Basis & & Dimension\\
\hline
\hV_{3/2}^1 &i\ca(-\tre)\bv&1\leq i \leq 4& 8\\
            &i(-\one)i^*(-\one)j\ca(-\one)\bv&1\leq i\neq j \leq 4& 24\\
            &i\ca(-\one)j\ca(-\one)k\ca(-\one)\bv& 1 \leq i<j<k  \leq 4 & 32\\

\hline
\hV_{2}^0 &i\ca(-\tre) j\ca(-\one)\bv&1\leq i,j \leq 4& 64\\
            &i(-\one)i^*(-\one)j(-\one)j^*(-\one)\bv&1\leq i<j \leq 4& 6\\
            &i(-\one)i^*(-\one)j\ca(-\one)k\ca(-\one)\bv& 1 \leq i  \leq 4, & \\
            &&i \neq j, i\neq k, &48\\
             &&1 \leq  j<k \leq 4&\\
            &1\ca(-\one)2\ca(-\one)3\ca(-\one)4\ca(-\one)\bv&&16\\
\hline
\end{array}
$$

\begin{landscape}
\section{Ramond Tables}\label{sec:Ramond}

\subsection{Table Relevant to $G_2^{(1)}$}\label{sub:g2tabR}
\begin{table}[h!]
$$
\small
\begin{array}{|c|c|c|c|c|c|c|}
\hline
\color{orange} X_{\beta_2}(0) &  {\bf1'}& 1^*(0)4^*(0){\bf1'}  & 2^*(0)3^*(0) {\bf1'} & 4^*(0){\bf1'} & 1^*(0){\bf1'}  & 2^*(0)3^*(0)4^*(0) {\bf1'}\\
\hline
 - 2^*(0)1(0) &0&-2^*(0)4^*(0) {\bf1'}&0&0&-2^*(0) {\bf1'}&0\\
 \hline
 - 4^*(0)3(0) &0&0&-2^*(0)4^*(0) {\bf1'}&0&0&0\\
 \hline
  4(0)3(0) &0&0&0&0&0&2^*(0) {\bf1'}\\
\hline
\end{array}
$$
\caption{Application of operators \eqref{usebeta2}}\label{tab:beta2R}
\end{table}

\subsection{Table Relevant to $L_0^{7/10}$}\label{sub:g2tabR}

\begin{table}[h!]
{\small
\noindent
$$
\begin{array}{|c|c|c|c|c|c|c|}
\hline
r_i=0 &{\bf1'}& 1^*(0)4^*(0){\bf1'}  & 2^*(0)3^*(0) {\bf1'} & 4^*(0){\bf1'} & 1^*(0){\bf1'}  & 2^*(0)3^*(0)4^*(0) {\bf1'}\\
\hline
1^*(0)2^*(0)1(0)2(0) & 0  &  0 & 0  & 0  & 0  & 0    \\
1^*(0)3^*(0)1(0)3(0) &  0 & 0  & 0  & 0  & 0  & 0    \\
2^*(0)3^*(0)2(0)3(0) &  0 & 0  & -2^*(0)3^*(0){\bf1'}   & 0  & 0  &  -2^*(0)3^*(0)4^*(0){\bf1'}     \\
\hline
1^*(0)1(0)  &  0 &  1^*(0)4^*(0){\bf1'} &  0 & 0  & 1^*(0){\bf1'}  & 0    \\
2^*(0)2(0)  &  0 &  0 & 2^*(0)3^*(0){\bf1'}   & 0  & 0  &   2^*(0)3^*(0)4^*(0){\bf1'}    \\
3^*(0)3(0)  & 0  &  0 & 2^*(0)3^*(0){\bf1'}  &  0 & 0  &  2^*(0)3^*(0)4^*(0){\bf1'}    \\
\hline
 1^*(0)2(0)3(0)4^*(0)& 0 & 0 &  -1^*(0)4^*(0){\bf1'}    & 0  & 0  & 0    \\
 1(0)2^*(0)3^*(0)4(0) & 0  & -2^*(0)3^*(0){\bf1'}   & 0  & 0  & 0  & 0    \\
\hline
 1^*(0)2(0)3(0)4(0) & 0  & 0  & 0  & 0  & 0  & - 1^*(0){\bf1'}    \\
 1(0)2^*(0)3^*(0)4^*(0)& 0  &  0 & 0  & 0  & -2^*(0)3^*(0)4^*(0){\bf1'}   &  0   \\
\hline
\end{array}
$$
}
\caption{Application of quartic operators for $L_0^{7/10}$}\label{tab:L1R}
\end{table}
\end{landscape}

\section{Tables for the Remaining Neveu-Schwarz  HWV's  }\label{sec:NSunused}
\subsection{Tables Relevant to $G_2^{(1)}$}\label{sub:unusedG2}

\begin{table}[h!]
$$
\begin{array}{|c|c|c|c||c|c|c|}
\hline
\color{red} X_{\beta_1}(0) &4(-\one)\bf1&4^*(-\one)\bf1 & 1(-\one)\bf1 & 1(-\one)4(-\one)\bf1  & 1(-\one)4^*(-\one)\bf1 & 2(-\one)3(-\one)\bf1 \\
\hline
2(-\one)3^*(\one) & 0& 0& 0  & 0 & 0  & 0 \\
3^*(-\one)2(\one) & 0& 0 & 0 & 0 & 0 & 0\\
\hline
\hline
\color{red} X_{\beta_2}(0)&&&&&&                \\
\hline
1(-\one)2^*(\one) &0&0& 0 & 0 & 0&1(-\one)3(-\one)\bf1\\
2^*(-\one)1(\one)&0&0& 0 & 0 & 0 & 0\\
\hline
3(-\one)4^*(\one) &3(-\one){\bf1}& 0&0&1(-\one)3(-\one)\bf1 & 0 & 0\\
4^*(-\one)3(\one)&0&0& 0 & 0 & 0 & 0\\
\hline
3(-\one)4(\one) &0&3(-\one){\bf1}& 0 & 0&1(-\one)3(-\one)\bf1 & 0\\
4(-\one)3(\one)&0&0& 0 & 0 & 0 & 0\\
\hline
\hline
\color{red} X_{-\theta}(1) &&&&&&                   \\
\hline
2^*(\one)1^*(\one) & 0& 0& 0 & 0 & 0 & 0 \\
 \hline
\end{array}
$$
\caption{Action of the $G_2^{(1)}$ operators on the basis vectors discussed in \ref{thm:hwvs}} \label{tab:rest}
\end{table}

\begin{table}[h!]
$$
\begin{array}{|c|c|c|c|c|c|c|c|c|}
\hline
\color{red} X_{\beta_1}(0) &11^*4 & 11^*4^* & 22^*4 & 22^*4^* & 33^*4 & 33^*4^* &1^*23 &12^*3^*\\
\hline
2(-\one)3^*(\one) &0&0&0&0&23^*4&23^*4^*&0&0\\
3^*(-\one)2(\one) &0&0&23^*4&23^*4^*&0&0&0&0\\
\hline
\hline
\color{red} X_{\beta_2}(0)  &&&&&&&&\\
\hline
1(-\one)2^*(\one) &0&0&12^*4&12^*4^*&0&0&-11^*3&0\\
2^*(-\one)1(\one) &12^*4&12^*4^*&0&0&0&0&-22^*3&0\\
\hline
3(-\one)4^*(\one) &11^*3&0&22^*3&0&0&0&0&0\\
4^*(-\one)3(\one) &0&0&0&0&-344^*&0&0&12^*4^*\\
\hline
3(-\one)4(\one) &0&11^*3&0&22^*3&0&0&0&0\\
4(-\one)3(\one) &0&0&0&0&0&344^*&0&12^*4\\
\hline
\hline
\color{red} X_{-\theta}(1) &&&&&&&&\\
\hline
2^*(\one)1^*(\one) &0&0&0&0&0&0&0&0\\
\hline
\end{array}
$$
\caption{Action of the $G_2^{(1)}$ operators on the basis vectors discussed in \ref{thm:hwvs}} \label{tab:rest}
\end{table}

\begin{table}[h!]
$$
\begin{array}{|c|c|c|c|c|c|c|c|}
\hline
\color{red} X_{\beta_1}(0)  &11^*44^* & 22^*44^*& 33^*44^* & 1^*234 & 1^*234^* & 12^*3^*4 & 12^*3^*4^*\\
\hline
2(-\one)3^*(\one) &0&0&23^*44^*&0&0&0&0\\
3^*(-\one)2(\one) &0&23^*44^*&0&0&0&0&0\\
\hline
\hline
\color{red} X_{\beta_2}(0)  &&&&&&&\\
\hline
1(-\one)2^*(\one) &0&12^*44^*&0&-11^*34&-11^*34^*&0&0\\
2^*(-\one)1(\one)  &12^*44^*&0&0&-22^*34&-22^*34^*&0&0\\
\hline
3(-\one)4^*(\one)  &11^*34^*&22^*34^*&0&0&0&-12^*33^*&0\\
4^*(-\one)3(\one) &0&0&0&0&0&-12^*44^*&0\\
\hline
3(-\one)4(\one)  &-11^*34&-22^*34&0&0&0&0&-12^*33^*\\
4(-\one)3(\one)  &0&0&0&0&0&0&12^*44^*\\
\hline
\hline
\color{red} X_{-\theta}(1) &&&&&&&\\
\hline
2^*(\one)1^*(\one) &0&0&0&0&0&0&0\\
\hline
\end{array}
$$
\caption{Action of the $G_2^{(1)}$ operators on the basis vectors discussed in \ref{thm:hwvs}} \label{tab:rest}
\end{table}

\newpage
\subsection{Table Relevant to $L_2^c$}\label{sub:unusedL2}

\begin{table}[h!]
$$
\begin{array}{|c|c|c|c|c|c|c|c|}
\hline
&11^*44^* & 22^*44^*& 33^*44^* & 1^*234 & 1^*234^* & 12^*3^*4 & 12^*3^*4^*\\
\hline
11^*22^*  &0&0&0&0&0&0&0\\  
11^*33^*  &0&0&0&0&0&0&0\\   
22^*33^*  &0&0&0&0&0&0&0\\  
1^*234  &0&0&0&0&0&0&\bf1\\  
1^*234^*  &0&0&0&0&0&\bf1&0\\  
12^*3^*4  &0&0&0&0&\bf1&0&0\\  
12^*3^*4^*  &0&0&0&\bf1&0&0&0\\  
\hline

\end{array}
$$
\caption{Action of the $L_2^{c}$ operators on the basis vectors discussed in \ref{thm:hwvs}} \label{tab:rest}
\end{table}

\begin{landscape}
\subsection{Tables Relevant to $L_1^c$}\label{sub:unusedL1}

\begin{table}[h!]
$$
\begin{array}{|c|c|c|c|c|c|c|c|c|}
\hline
&11^*4 & 11^*4^* & 22^*4 & 22^*4^* & 33^*4 & 33^*4^* &1^*23 &12^*3^*\\
\hline
1(\one)1^*(\one) &4(-\one){\bf1}&4^*(-\one){\bf1}&0&0&0&0&0&0\\
1^*(\one)1(\one)&-4(-\one){\bf1}&-4^*(-\one){\bf1}&0&0&0&0&0&0\\
2(\one)2^*(\one) &0&0&4(-\one){\bf1}&4^*(-\one){\bf1}&0&0&0&0\\
2^*(\one)2(\one)&0&0&-4(-\one){\bf1}&-4^*(-\one){\bf1}&0&0&0&0\\
3(\one)3^*(\one) &0&0&0&0&4(-\one){\bf1}&4^*(-\one){\bf1}&0&0\\
3^*(\one)3(\one)&0&0&0&0&-4(-\one){\bf1}&-4^*(-\one){\bf1}&0&0\\
4(\one)4^*(\one) &0&0&0&0&0&0&0&0\\
4^*(\one)4(\one)&0&0&0&0&0&0&0&0\\
\hline
4(\one)4(\one) &0&0&0&0&0&0&0&0\\
4^*(\one)4^*(\one)&0&0&0&0&0&0&0&0\\
\hline
\hline
11^*22^* &0&0&0&0&0&0&0&0\\
11^*33^* &0&0&0&0&0&0&0&0\\
22^*33^*  &0&0&0&0&0&0&0&0\\
1^*234&0&0&0&0&0&0&0&4(-\one){\bf1}\\
12^*3^*4^*&0&0&0&0&0&0&4^*(-\one){\bf1}&0\\
1^*234^*&0&0&0&0&0&0&0&4^*(-\one){\bf1}\\
12^*3^*4&0&0&0&0&0&0&4(-\one){\bf1}&0\\
\hline
\end{array}
$$
\caption{Action of the $L_1^c$ operators on the basis vectors discussed in \ref{thm:hwvs}} \label{tab:rest}
\end{table}

\begin{table}[h!]

$$
\begin{array}{|c|c|c|c|c|c|c|c|}
\hline
&11^*44^* & 22^*44^*& 33^*44^* & 1^*234 \\
\hline
1(\one)1^*(\one) &4(-\one)4^*(-\one){\bf1}&0&0&0\\
1^*(\one)1(\one)&-4(-\one)4^*(-\one){\bf1}&0&0&0\\
2(\one)2^*(\one) &0&4(-\one)4^*(-\one){\bf1}&0&0\\
2^*(\one)2(\one) &0&-4(-\one)4^*(-\one){\bf1}&0&0\\
3(\one)3^*(\one) &0&0&4(-\one)4^*(-\one){\bf1}&0\\
3^*(\one)3(\one) &0&0&-4(-\one)4^*(-\one){\bf1}&0\\
4(\one)4^*(\one) &1(-\one)1^*(-\one){\bf1}&2(-\one)2^*(-\one){\bf1}&3(-\one)3^*(-\one){\bf1}&0\\
4^*(\one)4(\one) &-1(-\one)1^*(-\one){\bf1}&-2(-\one)2^*(-\one){\bf1}&-3(-\one)3^*(-\one){\bf1}&0\\
\hline
4(\one)4(\one) &0&0&0&0\\
4^*(\one)4^*(\one)&0&0&0&0\\
\hline
\hline
11^*22^* &0&0&0&0\\
11^*33^* &0&0&0&0\\
22^*33^*  &0&0&0&0\\
1^*234&0&0&0&0\\
12^*3^*4^*&0&0&0& 1(-\one)1^*(-\one){\bf1}-2(-\one)2^*(-\one){\bf1} \\
&&&& -3(-\one)3^*(-\one){\bf1}-4(-\one)4^*(-\one){\bf1} \\
1^*234^*&0&0&0&0\\
12^*3^*4&0&0&0&0\\
\hline
\end{array}
$$
\caption{Action of the $L_1^c$ operators on the basis vectors discussed in \ref{thm:hwvs} and continued below} \label{tab:rest}
\end{table}

\begin{table}[h!]
$$
\begin{array}{|c|c|c|c|c|c|c|c|}
\hline
& 1^*234^* & 12^*3^*4 & 12^*3^*4^*\\
\hline
1(\one)1^*(\one) &0&0&0\\
1^*(\one)1(\one)&0&0&0\\
2(\one)2^*(\one) &0&0&0\\
2^*(\one)2(\one) &0&0&0\\
3(\one)3^*(\one) &0&0&0\\
3^*(\one)3(\one) &0&0&0\\
4(\one)4^*(\one) &0&0&0\\
4^*(\one)4(\one) &0&0&0\\
\hline
4(\one)4(\one)&0&0&0\\
4^*(\one)4^*(\one)&0&0&0\\
\hline
\hline
11^*22^* &0&0&0\\
11^*33^* &0&0&0\\
22^*33^*  &0&0&0\\
1^*234&0&0&-1(-\one)1^*(-\one){\bf1}+2(-\one)2^*(-\one){\bf1}\\
&&&3(-\one)3^*(-\one){\bf1}+4(-\one)4^*(-\one){\bf1} \\
12^*3^*4^* &0&0&0\\
1^*234^*&0&-1(-\one)1^*(-\one){\bf1}+2(-\one)2^*(-\one){\bf1}&0\\
&&3(-\one)3^*(-\one){\bf1}-4(-\one)4^*(-\one){\bf1}&\\
12^*3^*4& 1(-\one)1^*(-\one){\bf1}&0&0\\
& -2(-\one)2^*(-\one){\bf1}&&\\
& -3(-\one)3^*(-\one){\bf1}&&\\
& +4(-\one)4^*(-\one){\bf1}&&\\
\hline
\end{array}
$$
\caption{Action of the $L_1^c$ operators on the basis vectors discussed in \ref{thm:hwvs} continuation} \label{tab:rest}
\end{table}
\end{landscape}

\subsection{Tables Relevant to $L_0^c$}\label{sub:unusedL0}

\begin{table}[h!]
$$
\begin{array}{|c|c|c|c||c|c|c|}
\hline
&4(-\one)\bf1&4^*(-\one)\bf1 & 1(-\one)\bf1 & 1(-\one)4(-\one)\bf1  & 1(-\one)4^*(-\one)\bf1 & 2(-\one)3(-\one)\bf1 \\
\hline
4(-\one)4(\one) & 0&4(-\one)\bf1  & 0  & 0 & 1(-\one)4(-\one)\bf1   & 0 \\
4(-\one)4^*(\one) &4(-\one)\bf1 & 0& 0  &  1(-\one)4(-\one)\bf1& 0  & 0 \\
4^*(-\one)4(\one) &0 & 4^*(-\one)\bf1 & 0  &  0&   1(-\one)4^*(-\one)\bf1 & 0 \\
4^*(-\one)4^*(\one) & 4^*(-\one)\bf1 &0 & 0  &  1(-\one)4^*(-\one)\bf1 &   0& 0 \\
1^*(-\one)1(\one) &0 & 0&0   &0  & 0  & 0 \\
2^*(-\one)2(\one) & 0&0 &  0 &0  & 0  & 0 \\
3^*(-\one)3(\one) & 0&0 &  0 &0  & 0  & 0 \\
1(-\one)1^*(\one) & 0&0 &1(-\one)\bf1   &  1(-\one)4(-\one)\bf1 & 1(-\one)4^*(-\one)\bf1   & 0 \\
2(-\one)2^*(\one) & 0&0 &  0 & 0 & 0  &2(-\one)3(-\one)\bf1  \\
3(-\one)3^*(\one) & 0&0 &  0 & 0 & 0  & 2(-\one)3(-\one)\bf1  \\
\hline
\hline
11^*22^* & 0&0 &0   & 0 & 0  &0  \\
11^*33^* & 0&0 &0   &  0&  0 & 0 \\
22^*33^* & 0&0 &0   &0  &  0 & 2(-\one)3(-\one)\bf1  \\
1^* 234 & 0&0 &0   & 0 & -2(-\one)3(-\one)\bf1  &  0\\
12^* 3^* 4^* &0 &0 &0   & 0 &  0 & -1(-\one)4^*(-\one)\bf1  \\
1^* 234^*  &0 &0 &0   & 2(-\one)3(-\one)\bf1  &  0 &0  \\
12^* 3^* 4 &0 &0 &0   & 0 &  0 &  -1(-\one)4(-\one)\bf1\\
\hline
\end{array}
$$
\caption{Action of the $L_0^c$ operators on the basis vectors discussed in \ref{thm:hwvs}} \label{tab:rest}
\end{table}
\begin{landscape}
\begin{table}[h!]
$$
\begin{array}{|c|c|c|c|c|c|c|c|c|}
\hline
 &11^*4 & -11^*4^* &- 22^*4 & 22^*4^* & -33^*4 & 33^*4^* &1^*23 &12^*3^*\\
\hline
4(-\one)4(\one) &0& -11^*4^*&0 &  22^*4^* &  0&   33^*4^*&0&0  \\
4(-\one)4^*(\one)&11^*4 &  0&-22^*4   &  0& -33^*4&0& 0 & 0 \\
4^*(-\one)4(\one)&0 &-11^*4^* &  0 & 22^*4^* &0 &33^*4^*& 0 &  0\\
4^*(-\one)4^*(\one)&11^*4^* &0 & -22^*4  & 0 &-33^*4 &0&0  &0  \\
1^*(-\one)1(\one) &-11^*4 & 11^*4^*& 0  &0  & 0&0&-1^*23  & 0 \\
2^*(-\one)2(\one) & 0& 0& 22^*4  & -22^*4^*&0   &0&0 &-12^*3^*  \\
3^*(-\one)3(\one) &0 &0 & 0  & 0 &33^*4  &-33^*4^*&0 &-12^*3^*  \\
1(-\one)1^*(\one) &-11^*4 & 11^*4^*& 0  &0  & 0&0& 0 &  -12^*3^*\\
2(-\one)2^*(\one) & 0& 0& 22^*4  & -22^*4^*&0   &0&-1^*23 &  0\\
3(-\one)3^*(\one) &0 &0 & 0  & 0 &33^*4  &-33^*4^*& -1^*23&  0\\
\hline
\hline
11^*22^* & 22^*4&-22^*4^* &-11^*4   &11^*4^*  &0 &0&1^*23  &12^*3^*  \\
11^*33^* & 33^*4&-33^*4^* &0&0&-11^*4   &11^*4^*  &-1^*23  &-12^*3^*  \\
22^*33^* &0&0& -33^*4&33^*4^* &-22^*4   &22^*4^*  &1^*23  &12^*3^*  \\
1^* 234 & 0&-1^*23 & 0  & -1^*23 &0&-1^*23&0   & -11^*4+22^*4\\
&&&&&&&  &+ 33^*4 \\
12^* 3^* 4^* &-12^*3^* & 0  & -12^*3^* &0&-12^*3^*&0   & 11^*4^*-22^*4^*&0 \\
&&& &&&  & -33^*4^*& \\
1^* 234^*  & 1^*23& 0& 1^*23  & 0 &1^*23 &0& 0 & -11^*4^*+22^*4^* \\
&&&&&&&  & +33^*4^* \\
12^* 3^* 4 &0 & 12^*3^*&0   &12^*3^*  &0&12^*3^*&11^*4-22^*4   &0  \\
&&&&&&  & -33^*4 &\\
\hline
\end{array}
$$
\caption{Action of the $L_0^c$ operators on the basis vectors discussed in \ref{thm:hwvs}} \label{tab:rest}
\end{table}
\end{landscape}
\begin{table}[h!]
$$
\begin{array}{|c|c|c|c|c|c|c|c|}
\hline
  &11^*44^* &- 22^*44^*&- 33^*44^* & 1^*234 & 1^*234^* & 12^*3^*4 & 12^*3^*4^*\\
\hline
4(-\one)4(\one) &0 &0 &  0 &  0&0&  1^*234 & 1^*234^* \\
4(-\one)4^*(\one) &11^*44^* &-22^*44^* & -33^*44^*  &12^*3^*4  &12^*3^*4^*&0   &0  \\
4^*(-\one)4(\one) &11^*44^* &-22^*44^* & -33^*44^*  &0  &0& 1^*234^*  & 1^*234 \\
4^*(-\one)4^*(\one) &0 & 0&  0&12^*3^*4^* & 12^*3^*4 & 0  & 0 \\
1^*(-\one)1(\one) &-11^*44^* & 0& 0  & -1^*234 &  0 & -1^*234^*&  0\\
2^*(-\one)2(\one) &0 &22^*44^* & 0  & 0 & -12^*3^*4&  0 & -12^*3^*4^* \\
3^*(-\one)3(\one) &0 &0 & 33^*44^*  & 0 &-12^*3^*4 &  0& -12^*3^*4^* \\
1(-\one)1^*(\one) & -11^*44^*&0 & 0  & 0 &-12^*3^*4 &  0 &-12^*3^*4^*  \\
2(-\one)2^*(\one) &0 & 22^*44^*&  0 & -1^*234 & 0 &   -1^*234^*& 0 \\
3(-\one)3^*(\one) &0 & 0& 33^*44^*  & -1^*234 & 0 &  -1^*234^* & 0\\
\hline
\hline
11^*22^* &22^*44^* & -11^*44^*& 0  & -1^*234 &-12^*3^*4  & -1^*234^* &-12^*3^*4^*  \\
11^*33^* & 33^*44^*&0 & -11^*44^*  & -1^*234 &-12^*3^*4 &  -1^*234^* &-12^*3^*4^*  \\
22^*33^* & 0&-33^*44^* & -22^*44^*  & 1^*234 &12^*3^*4 &  1^*234^* &12^*3^*4^*  \\
1^* 234 & -1^*234& -1^*234& -1^*234  & 0 & 0  & 0& -\fbo\\
&&&&&&&+\fbt\\
12^* 3^* 4^* & -12^*3^*4^*&-12^*3^*4^* &  -12^*3^*4^* & -\fbo & 0&  0 & 0 \\
&&&&+\fbt&&&\\
1^* 234^*  &1^*234^* & 1^*234^*& 1^*234^*  &0 & -\fbo  &0& 0 \\
&&&&&-\fbt  &&\\
12^* 3^* 4 &12^*3^*4 &12^*3^*4 &12^*3^*4   &0  & 0  &-\fbo& 0 \\
&&&&&&-\fbt &\\
\hline
\end{array}
$$
\caption{Action of the $L_0^c$ operators on the basis vectors discussed in \ref{thm:hwvs}} \label{tab:rest}
\end{table}

\newpage
\addcontentsline{toc}{chapter}{Bibliography} 


\fontsize{12}{12pt} \selectfont


\bibliographystyle{alphanum}
\bibliography{references}
\nocite{*}
\clearpage
\addcontentsline{toc}{chapter}{Index}
\printindex

\end{document}